# KOLMOGOROV EQUATIONS IN INFINITE DIMENSIONS: WELL-POSEDNESS AND REGULARITY OF SOLUTIONS, WITH APPLICATIONS TO STOCHASTIC GENERALIZED BURGERS EQUATIONS[1]

By Michael Röckner and Zeev Sobol

*Universität Bielefeld and University of Wales Swansea*

We develop a new method to uniquely solve a large class of heat equations, so-called Kolmogorov equations in infinitely many variables. The equations are analyzed in spaces of sequentially weakly continuous functions weighted by proper (Lyapunov type) functions. This way for the first time the solutions are constructed everywhere without exceptional sets for equations with possibly nonlocally Lipschitz drifts. Apart from general analytic interest, the main motivation is to apply this to uniquely solve martingale problems in the sense of Stroock–Varadhan given by stochastic partial differential equations from hydrodynamics, such as the stochastic Navier–Stokes equations. In this paper this is done in the case of the stochastic generalized Burgers equation. Uniqueness is shown in the sense of Markov flows.

**1. Introduction.** In this paper we develop a new technique to uniquely solve generalized heat equations, so-called Kolmogorov equations, in infinitely many variables of type

$$\frac{du}{dt} = Lu$$

for a large class of elliptic operators $L$. The main new idea is to study $L$ on weighted function spaces consisting of sequentially weakly continuous functions on the underlying infinite-dimensional Banach space $X$ (e.g., a classical $L^p$-space). These function spaces are chosen appropriately for the

Received January 2004; revised March 2005.

[1]Supported by the BiBoS-Research Centre and the DFG-Research Group "Spectral Analysis, Asymptotic Distributions and Stochastic Dynamics."

specifically given operator $L$. More precisely, the function space on which $L$ acts is weighted by a properly chosen Lyapunov function $V$ of $L$ and the image space by a function $\Theta$ bounding its image $LV$. Apart from general analytic interest, the motivaton for this work comes from the study of concrete stochastic partial differential equations (SPDEs), such as, for example, those occuring in hydrodynamics (stochastic Navier–Stokes or Burgers equations, etc.). Transition probabilities of their solutions satisfy such Kolmogorov equations in infinitely many variables. To be more specific, below we shall describe a concrete case, to which we restrict in this paper, to explain the method in detail.

Consider the following stochastic partial differential equation on

$$X := L^2(0,1) = L^2((0,1), dr)$$

(where $dr$ denotes Lebesgue measure):

$$\begin{aligned}
(1.1) \qquad dx_t &= (\Delta x_t + F(x_t))\, dt + \sqrt{A}\, dw_t \\
x_0 &= x \in X.
\end{aligned}$$

Here $A \colon X \to X$ is a nonnegative definite symmetric operator of trace class, $(w_t)_{t \geq 0}$ a cylindrical Brownian motion on $X$, $\Delta$ denotes the Dirichlet Laplacian (i.e., with Dirichlet boundary conditions) on $(0,1)$, and $F \colon H_0^1 \to X$ is a measurable vector field of type

$$F(x)(r) := \frac{d}{dr}(\Psi \circ x)(r) + \Phi(r, x(r)), \qquad x \in H_0^1(0,1), r \in (0,1).$$

$H_0^1 := H_0^1(0,1)$ denotes the Sobolev space of order 1 in $L^2(0,1)$ with Dirichlet boundary conditions and $\Psi \colon \mathbb{R} \to \mathbb{R}$, $\Phi \colon (0,1) \times \mathbb{R} \to \mathbb{R}$ are functions satisfying certain conditions specified below. In case $\Psi(x) = \frac{1}{2}x^2$, $\Phi \equiv 0$, SPDE (1.1) is just the classical stochastic Burgers equation, and if $\Psi \equiv 0$ and, for example, $\Phi(r,x) = -x^3$, we are in the situation of a classical stochastic reaction diffusion equation of Ginsburg–Landau type. Therefore, we call (1.1) "stochastic generalized Burgers equation."

Stochastic generalized Burgers equations have been studied in several papers. In fact, the first who included both a "hydrodynamic part" (i.e., $\Psi$ above) and a "reaction diffusion part" (i.e., $\Phi$ above) was Gyöngy in [29], where, as we do in this paper, he also considered the case where the underlying domain is $D = (0,1)$. Later jointly with Rovira in [31] he generalized his results to the case where $\Psi$ is allowed to have polynomial growth; $\Phi$ is still assumed to have linear growth and is locally Lipschitz with at most linearly growing Lipschitz constant. A further generalization to $d$-dimensional domains was done by the same two authors in [32]. Contrary to us, these authors purely concentrated on solving SPDE of type (1.1) directly and did not analyze the corresponding Kolmogorov equations. In fact, they can allow nonconstant (but globally Lipschitz) $\sqrt{A}$ and also explicitly time dependent coefficients. We refer to [29, 31, 32] for the exact conditions, but emphasize that always the reaction diffusion part is assumed to be locally Lipschitz and of at most linear growth. As we shall see below, for the solution of the Kolmogorov equations, our method allows the reaction diffusion part to be of polynomial growth (so Ginsburg–Landau is in fact included) and also the locally Lipschitz condition can be replaced by a much weaker condition of dissipative type [see conditions $(\Phi1)$–$(\Phi3)$ in Section 2 below].

SPDE of type (1.1) with either $\Psi \equiv 0$ or $\Phi \equiv 0$ have been studied extensively. For the case $\Psi \equiv 0$, the literature is so enormous that we cannot record it here, but instead refer, for example, to the monographs [24] and [13] and the references therein. For the case $\Phi \equiv 0$, we refer, for example, to [6, 10, 12, 18, 19, 30, 38, 39, 55], and for the classical deterministic case, for example, to [11, 33, 37, 41, 44]. References concerning the Kolmogorov equations for SPDE will be given below.

The motivation of handling both the hydrodynamic and reaction diffusion part in SPDE of type (1.1) together was already laid out in [29]. It is well known that the mathematical analysis is then much harder, standard theory has to be modified and new techniques must be developed. It is, however, somewhat imaginable that this, with some effort, can be done if as in [29, 31, 32] $\Phi$ has at most linear growth (see, e.g., Remark 8.2 in [35], where this is shown in a finite-dimensional situation). The case of $\Phi$ with polynomial growth treated in this paper seems, however, much harder. In contrast to [29, 31, 32], our methods require, on the other hand, that $\Psi$ grows less than $|x|^{5/2}$ for large $x$ [cf. condition $(\Psi)$ in Section 2].

Showing the range of our method by handling $\Phi$ and $\Psi$ together has the disadvantage that it makes the analysis technically quite hard. Therefore, the reader who only wants to understand the basic ideas of our new general approach is advised to read the paper under the assumption that $\Phi$ does not explicitly depend on $r$ and has polynomial growth strictly less than 5. This simplifies the analysis substantially [e.g., in definition (2.4) of the Lyapunov function below we can take $p = 2$, so the simpler weight functions in (2.3) below suffice].

But now let us turn back to the Kolmogorov equations corresponding to SPDE (1.1).

A heuristic (i.e., not worrying about existence of solutions) application of Itô's formula to (1.1) implies that the corresponding generator or Kolmogorov operator $L$ on smooth cylinder functions $u : X \to \mathbb{R}$, that is,

$$u \in \mathcal{D} := \mathcal{F}C_b^2 := \{u = g \circ P_N | N \in \mathbb{N}, g \in C_b^2(E_N)\} \qquad \text{(cf. below)},$$

is of the following form:

$$Lu(x) := \tfrac{1}{2}\operatorname{Tr}(AD^2u(x)) + (\Delta x + F(x), Du(x))$$

$$(1.2) \qquad = \tfrac{1}{2} \sum_{i,j=1}^{\infty} A_{ij} \partial_{ij}^2 u(x) + \sum_{k=1}^{\infty} (\Delta x + F(x), \eta_k) \partial_k u(x), \qquad x \in H_0^1.$$

Here $\eta_k(r) := \sqrt{2} \sin(\pi k r)$, $k \in \mathbb{N}$, is the eigenbasis of $\Delta$ in $L^2(0,1)$, equipped with the usual inner product $(\cdot, \cdot)$, $E_N := \mathrm{span}\{\eta_k | 1 \le k \le N\}$, $P_N$ is the corresponding orthogonal projection, and $A_{ij} := (\eta_i, A\eta_j)$, $i,j \in \mathbb{N}$. Finally, $Du$, $D^2u$ denote the first and second Fréchet derivatives, $\partial_k := \partial_{\eta_k}$, $\partial_{ij}^2 := \partial_{\eta_i} \partial_{\eta_j}$ with $\partial_y :=$ directional derivative in direction $y \in X$ and $(\Delta x, \eta_k) := (x, \Delta \eta_k)$ for $x \in X$.

Hence, the Kolmogorov equations corresponding to SPDE (1.1) are given by

$$(1.3) \qquad \begin{aligned} \frac{dv}{dt}(t,x) &= \bar{L} v(t,x), \qquad x \in X, \\ v(0, \cdot) &= f, \end{aligned}$$

where the function $f : X \to \mathbb{R}$ is a given initial condition for this parabolic PDE with variables in the infinite-dimensional space $X$. We emphasize that (1.3) is only reasonable for some extension $\bar{L}$ of $L$ (whose construction is an essential part of the entire problem) since even for $f \in \mathcal{D}$, it will essentially never be true that $v(t, \cdot) \in \mathcal{D}$.

Because of the lack of techniques to solve PDE in infinite dimensions, in situations as described above the "classical" approach to solve (1.3) was to first solve (1.1) and then show in what sense the transition probabilities of the solution solve (1.3) (cf., e.g., [3, 13, 17, 24, 26, 27, 45, 50] and the references therein). Since about 1998, however, a substantial part of recent work in this area (cf., e.g., [20, 52, 53] and one of the initiating papers, [46]) is based on the attempt to solve Kolmogorov equations in infinitely many variables [as (1.3) above] directly and, reversing strategies, use the solution to construct weak solutions, that is, solutions in the sense of a martingale problem as formulated by Stroock and Varadhan (cf. [54]) of SPDE as (1.1) above, even for very singular coefficients (naturally appearing in many applications). In the above quoted papers, as in several other works (e.g., [1, 4, 15, 16, 22, 23, 42]), the approach to solve (1.3) directly was, however, based on $L^p(\mu)$-techniques where $\mu$ is a suitably chosen measure depending on $L$, for example, $\mu$ is taken to be an infinitesimally invariant measure of $L$ (see below). So, only solutions to (1.3) in an $L^p(\mu)$-sense were obtained, in particular, allowing $\mu$-zero sets of $x \in X$ for which (1.3) does not hold or where (1.3) only holds for $x$ in the topological support of $\mu$ (cf. [20]).

In this paper we shall present a new method to solve (1.3) for all $x \in X$ (or an explicitly described subset thereof) not using any reference measure. It is based on finite-dimensional approximation, obtaining a solution which, despite the lack of (elliptic and) parabolic regularity results on infinite-dimensional spaces, will, nevertheless, have regularity properties. More precisely, setting $X_p := L^p((0,1), dr)$, we shall construct a semigroup of Markov probability kernels $p_t(x, dy)$, $x \in X_p$, $t > 0$, on $X_p$ such that, for all $u \in \mathcal{D}$, we have $t \mapsto p_t(|Lu|)(x)$ is locally Lebesgue integrable on $[0, \infty)$ and

$$(1.4) \qquad p_t u(x) - u(x) = \int_0^t p_s(Lu)(x)\, ds \qquad \forall x \in X_p.$$

Here, as usual for a measurable function $f : X_p \to \mathbb{R}$, we set

$$(1.5) \qquad p_t f(x) := \int f(y) p_t(x, dy), \qquad x \in X_p, t > 0,$$

if this integral exists. $p$ has to be large enough compared to the growth of $\Phi$ (cf. Theorem 2.2 below). Furthermore, $p_t$ for each $t > 0$ maps a class of sequentially weakly continuous (resp. a class of locally Lipschitz functions) growing at most exponentially into itself. That $p_t$, for $t > 0$, has the property to map the test function space $\mathcal{D}$ (consisting of finitely based, hence, sequentially weakly continuous functions) into itself (as is the case in finite dimensions at least if the coefficients are sufficiently regular) cannot be true in our case since $F$ depends on all coordinates of $x = \sum_{k=1}^{\infty} (x, \eta_k) \eta_k$ and not merely finitely many. So, the regularity property of $p_t$, $t > 0$, to leave the space of exponentially bounded (and, since it is Markov, hence, also the bounded) sequentially weakly continuous functions fixed is the next best possible.

As a second step, we shall construct a conservative strong Markov process with weakly continuous paths, which is unique under a mild growth condition and which solves the martingale problem given by $L$, as in (1.2) and, hence, also (1.1) weakly, for every starting point $x \in X_p$. We also construct an invariant measure for this process.

The precise formulation of these results require more preparations and are therefore postponed to the next section (cf. Theorems 2.2–2.4), where we also collect our precise assumptions. Now we would like to indicate the main ideas of the proof and the main concepts. First of all, we emphasize that these concepts are of a general nature and work in other situations as well (cf., e.g., the companion paper [47] on the 2D-stochastic Navier–Stokes equations). We restrict ourselves to the case described above, so in particular to the (one dimensional) interval $(0, 1)$ for the underlying state space $X_p = L^p((0,1), dr)$, in order to avoid additional complications.

The general strategy is to construct the semigroup solving (1.4) through its corresponding resolvent, that is, we have to solve the equation

$$(\lambda - L)u = f$$

for all $f$ in a function space and $\lambda$ large enough, so that all $u \in \mathcal{D}$ appear as solutions. The proper function spaces turn out to be weighted spaces of sequentially weakly continuous functions on $X$. Such spaces are useful since their dual spaces are spaces of measures, so despite the nonlocal compactness of the state space $X$, positive linear functionals on such function spaces over $X$ are automatically measures (hence, positive operators on it are automatically kernels of positive measures). To choose exponential weights is natural to make these function spaces, which will remain invariant under the to be constructed resolvents and semigroups, as large as possible. More precisely, one chooses a Lyapunov function $V_{p,\kappa}$ of $L$ with weakly compact level sets so that

$$(\lambda - L)V_{p,\kappa} \geq \Theta_{p,\kappa},$$

and so that $\Theta_{p,\kappa}$ is a "large" positive function of (weakly) compact level sets [cf. (2.3), (2.4) below for the precise definitions]. $\Theta_{p,\kappa}$ "measures" the coercivity of $L$ [or of SPDE (1.1)]. Then one considers the corresponding spaces $WC_{p,\kappa}$ and $W_1C_{p,\kappa}$ of sequentially weakly continuous functions over $X$, weighted by $V_{p,\kappa}$ and $\Theta_{p,\kappa}$, respectively, with the corresponding weighted supnorms [cf. (2.2) below]. Then for $\lambda$ large, we consider the operator

$$\lambda - L : \mathcal{D} \subset WC_{p,\kappa} \to W_1C_{p,\kappa}$$

and prove by an approximative maximum principle that, for some $m > 0$,

$$\|(\lambda - L)u\|_{W_1C_{p,\kappa}} \geq m\|u\|_{WC_{p,\kappa}}$$

(cf. Proposition 6.1). So we obtain dissipativity of this operator between these two different spaces and the existence of its continuous inverse $G_\lambda := (\lambda - L)^{-1}$. Considering a finite-dimensional approximation by operators $L_N$ on $E_N$, $N \in \mathbb{N}$, with nice coefficients, more precisely, considering their associated resolvents $(G_\lambda^N)_{\lambda>0}$, we show that $(\lambda - L)(\mathcal{D})$ has dense range and that the continuous extension of $G_\lambda$ to all of $W_1C_{p,\kappa}$ is still one-to-one ("essential maximal dissipativity"). Furthermore, $\lambda G_\lambda^N$ (lifted to all of $X$) converges uniformly in $\lambda$ to $\lambda G_\lambda$ which, hence, turns out to be strongly continuous, but only after restricting $G_\lambda$ to $WC_{p,\kappa}$, which is continuously embedded into $W_1C_{p,\kappa}$, so has a stronger topology (cf. Theorem 6.4). Altogether $(G_\lambda)_{\lambda \geq \lambda_0}$, $\lambda_0$ large, is a strongly continuous resolvent on $WC_{p,\kappa}$, so we can consider its inverse under the Laplace transform (Hille–Yosida theorem) to obtain the desired semigroup $(p_t)_{t>0}$ of operators which are automatically given by probability kernels as explained above. Then one checks that $p_t$, $t > 0$, solves (1.4) and is unique under a mild "growth condition" [cf. (2.17) and Proposition 6.7 below]. Subsequently, we construct a strong Markov process on $X_p$ with weakly continuous paths with transition semigroup $(p_t)_{t>0}$. By general theory, it then solves the Stroock–Varadhan martingale problem corresponding to $(L, \mathcal{D})$, hence, it weakly solves SPDE (1.1). We also prove its uniqueness in the set of all Markov processes satisfying the mild "growth condition" (2.18) below (cf. Theorem 7.1).

In comparison to other constructions of semigroups on weighted function spaces using locally convex topologies and the concept of bicontinuous semigroups (cf. [36] and the references therein), we emphasize that our spaces are (separable) Banach spaces so, as spaces with one single norm, are easier to handle.

In comparison to other constructions of infinite-dimensional Markov processes (see, e.g., [43, 52]) where capacitory methods were employed, we would like to point out that instead of proving the tightness of capacities, we construct Lyapunov functions (which are excessive functions in the sense of potential theory) with compact level sets. The advantage is that we obtain pointwise statements for all points in $X_p$, not just outside a set of zero capacity. Quite a lot is known about the approximating semigroups $(p_t^N)_{t>0}$, that is, the ones corresponding to the $(G_\lambda^{(N)})_{\lambda>0}$, $N \in \mathbb{N}$, mentioned above, since they solve classical finite-dimensional Kolmogorov equations with regular coefficients. So, our construction also leads to a way to "calculate" the solution $(p_t)_{t>0}$ of the infinite-dimensional Kolmogorov equation (1.3).

The organization of this paper is as follows: as already mentioned, in Section 2 we formulate the precise conditions (A) and (F1) on the diffusion coefficient $A$ and the drift $F$, respectively, and state our main results precisely. In Section 3 we prove the necessary estimates on $\mathbb{R}^N$, uniformly in $N$, which are needed for the finite-dimensional approximation. In Section 4 we introduce another assumption (F2) on $F$ which is the one we exactly need in the proof, and we show that it is weaker than (F1). In Section 5 we collect a few essential properties of our weighted function spaces on $X_p$. In particular, we identify their dual spaces which is crucial for our analysis. This part was inspired by [34]. The semigroup of kernels $p_t(x, dy)$, $t > 0$, $x \in X_p$, is constructed in Section 6, and its uniqueness is proved. Here we also prove further regularity properties of $p_t$, $t > 0$. The latter part is not used subsequently in this paper. Section 7 is devoted to constructing the process, respectively showing that it is the solution of the martingale problem given by $L$ as in (1.2), hence, a weak solution to SPDE (1.1), and that it is unique in the mentioned class of Markov processes (see also Lemma A.1 in the Appendix). In deterministic language the latter means that we have uniqueness of the flows given by solutions of (1.1). The invariant measure $\mu$ for $(p_t)_{t>0}$ is constructed in the Appendix by solving the equation $L^*\mu = 0$. As a consequence of the results in the main part of the paper, we get that the closure $(\bar{L}^\mu, \mathrm{Dom}(\bar{L}^\mu))$ of $(L, \mathcal{D})$ is maximal dissipative on $L^s(X, \mu)$, $s \in [1, \infty)$ (cf. Remark A.3), that is, strong uniqueness holds for $(L, \mathcal{D})$ on $L^s(X, \mu)$. In particular, the differential form (1.3) of (1.4) holds with $\bar{L}^\mu$ replacig $\bar{L}$ and the time derivative taken in $L^s(X, \mu)$.

**2. Notation, conditions and main results.** For a $\sigma$-algebra $\mathcal{B}$ on an arbitrary set $E$, we denote the space of all bounded (resp. positive) real-valued $\mathcal{B}$-measurable functions by $\mathcal{B}_b$, $\mathcal{B}^+$, respectively. If $E$ is equipped with a topology, then $\mathcal{B}(E)$ denotes the corresponding Borel $\sigma$-algebra. The spaces $X = L^2(0,1)$ and $H_0^1$ are as in the Introduction and they are equipped with their usual norms $|\cdot|_2$ and $|\cdot|_{1,2}$; so we define, for $x \colon (0,1) \to \mathbb{R}$, measurable,

$$|x|_p := \left( \int_0^1 |x(r)|^p \, dr \right)^{1/p} \ (\in [0,\infty]), \qquad p \in [1,\infty),$$

$$|x|_\infty := \operatorname*{ess\,sup}_{r \in (0,1)} |x(r)|,$$

and define $X_p := L^p((0,1), dr)$, $p \in [1,\infty]$, so $X = X_2$. If $x, y \in H_0^1$, set

$$|x|_{1,2} := |x'|_2, \qquad (x,y)_{1,2} := (x', y'),$$

where $x' := \frac{d}{dr} x$ is the weak derivative of $x$. We shall use this notation from now on and we also write $x'' := \frac{d^2}{dr^2} x = \Delta x$.

Let $H^{-1}$ with norm $|\cdot|_{-1,2}$ be the dual space of $H_0^1$. We always use the continuous and dense embeddings

$$(2.1) \qquad\qquad H_0^1 \subset X \equiv X' \subset H^{-1},$$

so ${}_{H_0^1}\langle x, y \rangle_{H^{-1}} = (x,y)$ if $x \in H_0^1$, $y \in X$. The terms "Borel-measurable" or "measure on $X$, $H_0^1$, $H^{-1}$ resp." will below always refer to their respective Borel $\sigma$-algebras, if it is clear on which space we work. We note that since $H_0^1 \subset X \subset H^{-1}$ continuously, by Kuratowski's theorem, $H_0^1 \in \mathcal{B}(X)$, $X \in \mathcal{B}(H^{-1})$ and $\mathcal{B}(X) \cap H_0^1 = \mathcal{B}(H_0^1)$, $\mathcal{B}(H^{-1}) \cap X = \mathcal{B}(X)$. Furthermore, the Borel $\sigma$-algebras on $X$ and $H_0^1$ corresponding to the respective weak topologies coincide with $\mathcal{B}(X)$, $\mathcal{B}(H_0^1)$, respectively.

For a function $V \colon X \to (0,\infty]$ having weakly compact level sets $\{V \le c\}$, $c \in \mathbb{R}_+$, we define

$$WC_V := \Big\{ f \colon \{V < \infty\} \to \mathbb{R} \Big| f \text{ is continuous on each } \{V \le R\}, R \in \mathbb{R},$$

$$(2.2) \qquad\qquad\qquad\qquad \text{in the weak topology inherited from } X,$$

$$\text{and } \lim_{R \to \infty} \sup_{\{V \ge R\}} \frac{|f|}{V} = 0 \Big\},$$

equipped with the norm $\|f\|_V := \sup_{\{V < \infty\}} V^{-1}|f|$. Obviously, $WC_V$ is a Banach space with this norm. We are going to consider various choices of $V$, distinguished by respective subindices, namely, we define, for $\kappa \in (0,\infty)$,

$$(2.3) \qquad \begin{aligned} V_\kappa(x) &:= e^{\kappa|x|_2^2}, & x \in X, \\ \Theta_\kappa(x) &:= V_\kappa(x)(1 + |x'|_2^2), & x \in H_0^1, \end{aligned}$$

and for $p > 2$,

$$\text{(2.4)} \quad \begin{aligned} V_{p,\kappa}(x) &:= e^{\kappa|x|_2^2}(1 + |x|_p^p), & x \in X, \\ \Theta_{p,\kappa}(x) &:= V_{p,\kappa}(x)(1 + |x'|_2^2) + V_\kappa(x)|(|x|^{p/2})'|_2^2, & x \in H_0^1. \end{aligned}$$

Clearly, $\{V_{p,\kappa} < \infty\} = X_p$ and $\{\Theta_{p,\kappa} < \infty\} = H_0^1$. Each $\Theta_{p,\kappa}$ is extended to a function on $X$ by defining it to be equal to $+\infty$ on $X \setminus H_0^1$. Abusing notation, for $p = 2$, we also set $V_{2,\kappa} := V_\kappa$ and $\Theta_{2,\kappa} := \Theta_\kappa$. For abbreviation, for $\kappa \in (0, \infty)$, $p \in [2, \infty)$, we set

$$\text{(2.5)} \quad WC_{p,\kappa} := WC_{V_{p,\kappa}}, \qquad W_1 C_{p,\kappa} := WC_{\Theta_{p,\kappa}},$$

and we also abbreviate the norms correspondingly,

$$\text{(2.6)} \quad \|\cdot\|_{p,\kappa} := \|\cdot\|_{V_{p,\kappa}}, \qquad \|\cdot\|_\kappa := \|\cdot\|_{0,\kappa} \quad \text{and} \quad \|\cdot\|_{1,p,\kappa} := \|\cdot\|_{\Theta_{p,\kappa}}.$$

All these norms are, of course, well defined for any function on $X$ with values in $[-\infty, \infty]$. And therefore we shall apply them below not just for functions in $WC_{p,\kappa}$ or $W_1 C_{p,\kappa}$. For $p' \geq p$ and $\kappa' \geq \kappa$, by restriction, $WC_{p,\kappa}$ is continuously and densely embedded into $WC_{p',\kappa'}$ and into $W_1 C_{p,\kappa}$ (see Corollary 5.6 below), as well is the latter into $W_1 C_{p',\kappa'}$. $V_{p,\kappa}$ will serve as convenient Lyapunov functions for $L$. Furthermore, $\Theta_{p,\kappa}$ bounds $(\lambda - L)V_{p,\kappa}$ from below for large enough $\lambda$, thus, $\Theta_{p,\kappa}$ measures the coercivity of $L$ (cf. Lemma 4.6 below). Note that the level sets of $\Theta_{p,\kappa}$ are even strongly compact in $X$.

We recall that, for $P_N$ as in the Introduction, there exists $\alpha_p \in [1, \infty)$ such that

$$\text{(2.7)} \quad |P_N x|_p \leq \alpha_p |x|_p \qquad \text{for all } x \in X_p, N \in \mathbb{N}$$

(cf. [40], Section 2c16), of course, with $\alpha_2 = 1$. In particular,

$$\text{(2.8)} \quad V_{\kappa,p} \circ P_N \leq \alpha_p^p V_{\kappa,p}.$$

For a function $V : X \to (1, \infty]$, we also define spaces $\text{Lip}_{l,p,\kappa}$, $p \geq 2$, $\kappa > 0$, consisting of functions on $X$ which are locally Lipschitz continuous in the norm $|(-\Delta)^{-l/2} \cdot|_2$, $l \in \mathbb{Z}_+$. The respective seminorms are defined as follows:

$$\text{(2.9)} \quad (f)_{l,p,\kappa} := \sup_{y_1, y_2 \in X_p} (V_{p,\kappa}(y_1) \vee V_{p,\kappa}(y_2))^{-1} \frac{|f(y_1) - f(y_2)|}{|(-\Delta)^{-l/2}(y_1 - y_2)|_2}$$

$$(\in [0, \infty]).$$

For $l \in \mathbb{Z}_+$, we define

$$\text{(2.10)} \quad \text{Lip}_{l,p,\kappa} := \{f : X_p \to \mathbb{R} | \|f\|_{\text{Lip}_{l,p,\kappa}} < \infty\},$$

where $\|f\|_{\text{Lip}_{l,p,\kappa}} := \|f\|_{p,\kappa} + (f)_{l,p,\kappa}$. When $X$ is of finite dimension, $(f)_{l,p,\kappa}$ is a weighted norm of the generalized gradient of $f$ (cf. Lemma 3.6 below).

Also, $(\mathrm{Lip}_{l,p,\kappa}, \|\cdot\|_{\mathrm{Lip}_{l,p,\kappa}})$ is a Banach space (cf. Lemma 5.7 below) and $\mathrm{Lip}_{l,p,\kappa} \subset \mathrm{Lip}_{l',p',\kappa'}$ for $l' \le l$, $p' \ge p$ and $\kappa' \ge \kappa$. In this paper we shall mostly deal with the case $l \in \{0, 1\}$.

Obviously, each $f \in \mathrm{Lip}_{l,p,\kappa}$ is uniformly $|(-\Delta)^{-l/2} \cdot |_2$-Lipschitz continuous on every $|\cdot|_p$-bounded set. In particular, any $f \in \mathrm{Lip}_{1,p,\kappa}$ is sequentially weakly continuous on $X_p$, consequently weakly continuous on bounded subsets of $X_p$. Hence, for all $p' \in [p, \infty)$, $\kappa' \in [\kappa, \infty)$,

$$(2.11) \qquad \mathcal{B}_b(X_p) \cap \mathrm{Lip}_{1,p,\kappa} \subset WC_{p',\kappa'}$$

and obviously, by restriction,

$$(2.12) \qquad \mathcal{B}_b(X_p) \cap \mathrm{Lip}_{0,p,\kappa} \subset W_1 C_{p',\kappa'}.$$

Further properties of these function spaces will be studied in Section 5 below.

Besides the space $\mathcal{D} := \mathcal{F}C_b^2$ defined in the Introduction, other test function spaces $\mathcal{D}_{p,\kappa}$ on $X$ will turn out to be convenient. They are for $p \in [2, \infty)$, $\kappa \in (0, \infty)$ defined as follows:

$$(2.13) \qquad \mathcal{D}_{p,\kappa} := \{u = g \circ P_N | N \in \mathbb{N}, g \in C^2(\mathbb{R}^N),$$
$$\|u\|_{p,\kappa} + \||Du|_2\|_{p,\kappa} + \|\operatorname{Tr}(AD^2u)\|_{p,\kappa} < \infty\}.$$

Again we set $\mathcal{D}_\kappa := \mathcal{D}_{2,\kappa}$. Obviously, $\mathcal{D}_{p,\kappa} \subset WC_{p,\kappa}$ and $\mathcal{D}_{p,\kappa} \subset \mathcal{D}_{p',\kappa'}$ if $p' \in [p, \infty)$ and $\kappa' \in [\kappa, \infty)$. We extend the definition (1.2) of the Kolmogorov operator $L$ for all $u \in \mathcal{F}C^2 := \{u = g \circ P_N | N \in \mathbb{N}, g \in C^2(\mathbb{R}^\mathbb{N})\}$. So, $L$ can be considered with domain $\mathcal{D}_{p,\kappa}$.

Now let us collect our precise hypotheses on the terms in SPDE (1.1), respectively the Kolmogorov operator (1.2). First, we recall that in the entire paper $\Delta = x''$ is the Dirichlet Laplacian on $(0, 1)$ and $(W_t)_{t \ge 0}$ is a cylindrical Browninan motion on $X$. Consider the following condition on the map $A : X \to X$:

(A)  $A$ is a nonnegative symmetric linear operator from $X$ to $X$ of trace class such that $A_N := P_N A P_N$ is an invertible operator represented by a diagonal matrix on $E_N$ for all $N \in \mathbb{N}$.

Here $E_N$, $P_N$ are as defined in the Introduction. Furthermore, we set

$$(2.14) \qquad a_0 := \sup_{x \in H_0^1 \setminus \{0\}} \frac{(x, Ax)}{|x'|_2^2} = |A|_{H_0^1 \to H^{-1}},$$

where $|\cdot|_{H_0^1 \to H^{-1}}$ denotes the usual operator norm on bounded linear operators from $H_0^1$ into its dual $H^{-1}$.

Consider the following conditiions on the map $F : H_0^1 \to X$:

(F1)

$$F(x) = \frac{d}{dr}(\Psi \circ x)(r) + \Phi(r, x(r)), \qquad x \in H_0^1(0,1), r \in (0,1), \tag{2.15}$$

where $\Psi \colon \mathbb{R} \to \mathbb{R}$, $\Phi \colon (0,1) \times \mathbb{R} \to \mathbb{R}$ satisfy the following conditions:

($\Psi$)  $\Psi \in C^{1,1}(\mathbb{R})$ (i.e., $\Psi$ is differentiable with locally Lipschitz derivative) and there exist $C \in [0, \infty)$ and a bounded, Borel-measurable function $\omega \colon [0, \infty) \to [0, \infty)$ vanishing at infinity such that

$$|\Psi_{xx}|(x) \leq C + \sqrt{|x|}\,\omega(|x|) \qquad \text{for } dx\text{-a.e. } x \in \mathbb{R}.$$

($\Phi$1)  $\Phi$ is Borel-measurable in the first and continuous in the second variable and there exists $g \in L^{q_1}(0,1)$ with $q_1 \in [2, \infty]$ and $q_2 \in [1, \infty)$ such that

$$|\Phi(r, x)| \leq g(r)(1 + |x|^{q_2}) \qquad \text{for all } r \in (0,1), x \in \mathbb{R}.$$

($\Phi$2)  There exist $h_0, h_1 \in L_+^1(0,1)$, $|h_1|_1 < 2$, such that for a.e. $r \in (0,1)$

$$\Phi(r, x)\,\text{sign}\,x \leq h_0(r) + h_1(r)|x| \qquad \text{for all } x \in \mathbb{R}.$$

($\Phi$3)  There exist $\rho_0 \in (0,1]$, $g_0 \in L_+^1(0,1)$, $g_1 \in L_+^{p_1}(0,1)$ for some $p_1 \in [2, \infty]$, and a function $\omega \colon [0, \infty) \to [0, \infty)$ as in ($\Psi$) such that with $\sigma \colon (0,1) \times \mathbb{R} \to \mathbb{R}$, $\sigma(r, x) := \frac{|x|}{\sqrt{r(1-r)}}$ for a.e. $r \in (0,1)$

$$\Phi(r, y) - \Phi(r, x) \leq [g_0(r) + g_1(r)|\sigma(r, x)|^{2 - 1/p_1}\omega(\sigma(r, x))](y - x)$$

for all $x, y \in \mathbb{R}$, $0 \leq y - x \leq \rho_0$.

Furthermore, we say that condition (F1+) holds if, in addition to (F1), we have

($\Phi$4)  $\Phi$ is twice continuously differentiable and there exist $g_2, g_3 \in L_+^2(0,1)$, $g_4, g_5 \in L_+^1(0,1)$, and $\omega \colon [0, \infty) \to [0, \infty)$ as in ($\Psi$) such that, for their partial derivatives $\Phi_{xx}$, $\Phi_{xr}$, $\Phi_x$, $\Phi_r$, and with $\sigma$ as in ($\Phi$3),

$$|\Phi_{xx}| + \frac{|\Phi_x|^2}{|\Phi| + 1} \leq g_2 + g_3\sqrt{\sigma}\,\omega(\sigma)$$

and

$$|\Phi_{xr}| + \frac{|\Phi_x \Phi_r|}{|\Phi| + 1} \leq g_4 + g_5\sigma^{3/2}\omega(\sigma).$$

REMARK 2.1.  (i) Integrating the inequality in ($\Psi$) twice, one immediately sees that ($\Psi$) implies that there exist a bounded Borel-measurable function $\hat\omega \colon \mathbb{R}_+ \to \mathbb{R}_+$, $\hat\omega(r) \to 0$ as $r \to \infty$, and $C \in (0, \infty)$ such that

$$|\Psi'(x)| \leq C + |x|^{3/2}\hat\omega(|x|), \qquad |\Psi(x)| \leq C + |x|^{5/2}\hat\omega(|x|) \qquad \text{for all } x \in \mathbb{R}.$$

(ii) We emphasize that conditions ($\Phi 2$), ($\Phi 3$) are one-sided estimates, so that ($\Phi 1$)–($\Phi 3$) is satisfied if $\Phi(r, x) = P(x), r \in (0, 1), x \in \mathbb{R}$, where $P$ is a polynomial of odd degree with strictly negative leading coefficient.

(iii) Under the assumptions in (F1), SPDE (1.1) will not have a strong solution in general for all $x \in X$.

(iv) If ($\Phi 1$) holds, ($\Phi 2$) only needs to be checked for $x \in \mathbb{R}$ such that $|x| \geq R$ for some $R \in (0, \infty)$. And replacing $\omega$ [in ($\Psi$) and ($\Phi 3$)] by $\tilde{\omega}(r) := \sup_{s \geq r} \omega(s)$, we may assume that $\omega$ is decreasing.

(v) ($\Phi 4$) implies that there exists a bounded measurable function $\hat{\omega} : \mathbb{R}_+ \to \mathbb{R}_+$, $\hat{\omega}(r) \to 0$ as $r \to \infty$, such that

$$|\Phi_x| \leq C + \sigma^{3/2} \hat{\omega}(\sigma) \quad \text{and} \quad |\Phi| \leq g_1 + C + \sigma^{5/2} \hat{\omega}(\sigma).$$

In particular, ($\Phi 4$) implies ($\Phi 3$) with $p_1 = 2$, $g_0(r) = g_1(r) = const$. Indeed, we have, for $x \in \mathbb{R}$, $r \in (0, \frac{1}{2})$,

$$\Phi_x(r, x) = \Phi_x(0, 0) + \int_0^r \Phi_{xx}\left(s, \frac{x}{r}s\right)\frac{x}{r}\,ds + \int_0^r \Phi_{xr}\left(s, \frac{x}{r}s\right)ds.$$

As shown in the previous item, we may assume $\omega$ decreasing. Then it follows from ($\Phi 4$) and Hölder's inequality that

$$
\begin{aligned}
|\Phi_x|(r, x) \leq\ & C + \frac{|x|}{r}\int_0^r g_2\,ds + \left(\frac{|x|}{r}\right)^{3/2}\int_0^r g_3(s)\omega\left(\frac{|x|}{r}\sqrt{s}\right)s^{1/4}\,ds + \int_0^r g_4\,ds \\
& + \left(\frac{|x|}{r}\right)^{3/2}\int_0^r g_5(s)\omega\left(\frac{|x|}{r}\sqrt{s}\right)s^{3/4}\,ds \\
\leq\ & C + \frac{|x|}{\sqrt{r}}|g_2|_2 + \left(\frac{|x|}{\sqrt{r}}\right)^{3/2}|g_3|_2\left(\int_0^1 \omega^2\left(\frac{|x|}{\sqrt{r}}\sqrt{\tau}\right)\sqrt{\tau}\,d\tau\right)^{1/2} \\
& + |g_4|_1 + \left(\frac{|x|}{\sqrt{r}}\right)^{3/2}\int_0^1 g_5(s)\omega\left(\frac{|x|}{\sqrt{r}}\sqrt{s}\right)ds.
\end{aligned}
$$

Now observe that

$$\tilde{\omega}(\sigma) := \left(2\int_0^1 \omega^2(\sqrt{2}\sigma\tau)\tau\,d\tau\right)^{1/2} + \int_0^1 g_5(s)\omega(\sqrt{2}\sigma\sqrt{s})\,ds$$

is a bounded measurable function and $\tilde{\omega}(r) \to 0$ as $r \to \infty$. So the first assertion follows for $r \in (0, \frac{1}{2})$. For the case $r \in (\frac{1}{2}, 1)$, the assertion is proved by the change of variables $r' = 1 - r$. The second assertion is proved similarly.

In the rest of this paper hypothesis (A) (though repeated in each statement to make partial reading possible) will always be assumed. As it is already said in the Introduction, all of our results are proved for general $F : H_0^1 \to X$ under condition (F2) [resp. (F2+), or parts thereof], which is introduced in Section 4 and which is weaker than (F1) [resp. (F1+)]. For the convenience of the reader, we now, however, formulate our results for the concrete $F$ given in (2.15), under condition (F1) [(F1+) resp.]. For their proofs, we refer to the respective more general results, stated and proved in one of the subsequent sections.

THEOREM 2.2 ("Pointwise solutions of the Kolmogorov equations"). *Suppose* (A) *and* (F1) *hold. Let* $\kappa_0 := \frac{2-|h_1|_1}{8a_0}$ *(with $a_0$ as in* (2.15) *and $h_1$ as in* (Φ2), $0 < \kappa_1 \in \kappa^* < \kappa_0$, *and let* $p \in [2, \infty) \cap (q_2 - 3 + \frac{2}{q_1}, \infty)$ *[with $q_1, q_2$ as in* (Φ1)]*). Then there exists a semigroup* $(p_t)_{t>0}$ *of probability kernels on* $X_p$, *independent of $\kappa^*$, having the following properties:*

(i) ("Existence") *Let* $u \in \mathcal{D}_{\kappa_1}$. *Then* $t \mapsto p_t(|Lu|)(x)$ *is locally Lebesgue integrable on* $[0, \infty)$ *and*

$$(2.16) \qquad p_t u(x) - u(x) = \int_0^t p_s(Lu)(x)\, ds \qquad \text{for all } x \in X_p.$$

*In particular, for all* $s \in [0, \infty)$,

$$\lim_{t \to 0} p_{s+t} u(x) = p_s u(x) \qquad \text{for all } x \in X_p.$$

(ii) *There exists* $\lambda_{\kappa^*} \in (0, \infty)$ *such that*

$$(2.17) \qquad \int_0^\infty e^{-\lambda_{\kappa^*} s} p_s(\Theta_{p,\kappa^*})(x)\, ds < \infty \qquad \text{for all } x \in X_p.$$

(iii) ("Uniqueness") *Let* $(q_t)_{t>0}$ *be a semigroup of probability kernels on* $X_p$ *satisfying* (i) *with* $(p_t)_{t>0}$ *replaced by* $(q_t)_{t>0}$ *and* $\mathcal{D}_{\kappa_1}$ *by* $\mathcal{D}$. *If, in addition,* (2.17) *holds with* $(q_t)_{t>0}$ *replacing* $(p_t)_{t>0}$ *for some* $\kappa \in (0, \kappa_0)$ *replacing* $\kappa^*$, *then* $p_t(x, dy) = q_t(x, dy)$ *for all* $t > 0, x \in X_p$.

(iv) ("Regularity") *Let* $t \in (0, \infty)$. *Then* $p_t(W_{p,\kappa^*}) \subset W_{p,\kappa^*}$. *Furthermore, let* $f \in \text{Lip}_{0,2,\kappa_1} \cap \mathcal{B}_b(X) \cap W_{p,\kappa^*} (\supset \mathcal{D})$. *Then* $p_t f$ *uniquely extends to a continuous function on* $X$, *again denoted by* $p_t f$, *which is in* $\text{Lip}_{0,2,\kappa_1} \cap \mathcal{B}_b(X)$. *Let* $q \in [2, \infty)$, $\kappa \in [\kappa_1, \kappa^*]$. *Then there exists* $\lambda_{q,\kappa} \in (0, \infty)$, *independent of $t$ and $f$, such that*

$$\|p_t f\|_{q,\kappa} \le e^{\lambda_{q,\kappa} t} \|f\|_{q,\kappa}$$

*and*

$$(p_t f)_{0,q,\kappa} \le e^{\lambda_{q,\kappa} t} (f)_{0,q,\kappa}.$$

*If moreover,* (F1+) *holds, then there exists* $\lambda'_{q,\kappa} \in (0, \infty)$, *independent of $t$, such that, for all* $f \in \text{Lip}_{1,2,\kappa_1} \cap \mathcal{B}_b(X)$,

$$(p_t f)_{1,q,\kappa} \le e^{\lambda'_{q,\kappa} t} (f)_{1,q,\kappa}.$$

PROOF.   The assertions follow from Corollary 4.2, Remark 6.6, Propositions 6.7, 6.9 and 6.11(iii).   □

THEOREM 2.3 ["Martingale and weak solutions to SPDE (1.1)"].   *Assume that* (A) *and* (F1) *hold, and let* $p, \kappa^*$ *be as in Theorem* 2.2.

(i) *There exists a conservative strong Markov process* $\mathbb{M} := (\Omega, \mathcal{F}, (\mathcal{F}_t)_{t \geq 0}, (x_t)_{t \geq 0}, (\mathbb{P}_x)_{x \in X_p})$ *in* $X_p$ *with continuous sample paths in the weak topology whose transition semigroup is given by* $(p_t)_{t > 0}$ *from Theorem* 2.2. *In particular, for* $\lambda_{\kappa^*}$ *as in Theorem* 2.2(ii),

$$\mathbb{E}_x \left[ \int_0^\infty e^{-\lambda_{\kappa^*} s} \Theta_{p, \kappa^*}(x_s) \, ds \right] < \infty \qquad \text{for all } x \in X_p.$$

(ii) ("Existence") *Let* $\kappa_1 \in (0, \kappa_0 - \kappa^*)$. *Then* $\mathbb{M}$ *satisfies the martingale problem for* $(L, \mathcal{D}_{\kappa_1})$, *that is, for all* $u \in \mathcal{D}_{\kappa_1}$ *and all* $x \in X_p$, *the function* $t \mapsto |Lu(x_t)|$ *is locally Lebesgue integrable on* $[0, \infty)$ $\mathbb{P}_x$-*a.s. and under* $\mathbb{P}_x$,

$$u(x_t) - u(x) - \int_0^t Lu(x_s) \, ds, \qquad t \geq 0,$$

*is an* $(\mathcal{F}_t)_{t \geq 0}$-*martingale starting at* 0 (*cf.* [54]).

(iii) ("Uniqueness") $\mathbb{M}$ *is unique among all conservative (not necessarily strong) Markov processes* $\mathbb{M}' := (\Omega', \mathcal{F}', (\mathcal{F}'_t)_{t \geq 0}, (x'_t)_{t \geq 0}, (\mathbb{P}'_x)_{x \in X_p})$ *with weakly continuous sample paths in* $X_p$ *satisfying the martingale problem for* $(L, \mathcal{D})$ [*as specified in* (ii) *with* $\mathcal{D}$ *replacing* $\mathcal{D}_{\kappa_1}$] *and having the additional property that, for some* $\kappa \in (0, \kappa_0)$, *there exists* $\lambda_\kappa \in (0, \infty)$ *such that*

$$(2.18) \qquad \mathbb{E}'_x \left[ \int_0^\infty e^{-\lambda_\kappa s} (\Theta_{p, \kappa})(x'_s) \, ds \right] < \infty \qquad \text{for all } x \in X_p.$$

(iv) *If* $p \geq 2q_2 - 6 + 4/q_1$, *then* $\mathbb{M}$ *weakly solves SPDE* (1.1).

PROOF.   Corollary 4.2, Remark 6.6, Theorem 7.1 and Remark 7.2 below.   □

THEOREM 2.4 ("Invariant measure").   *Assume that* (A) *and* (F1) *hold. Let* $p, \kappa^*$ *be as in Theorem* 2.2.

(i) *There exists a probability measure* $\mu$ *on* $H_0^1$ *which is "L-infinitesimally invariant," that is,* $Lu \in L^1(H_0^1, \mu)$ *and*

$$(2.19) \qquad \int Lu \, d\mu = 0 \qquad \text{for all } u \in \mathcal{D}$$

($L^* \mu = 0$ *for short*). *Furthermore,*

$$(2.20) \qquad \int \Theta_{p, \kappa^*} \, d\mu < \infty.$$

(ii) $\mu$, extended by zero to all of $X_p$, is $(p_t)_{t>0}$-invariant, that is, for all $f\colon X\to\mathbb{R}$, bounded, measurable, and all $t>0$,

$$\int p_t f\,d\mu = \int f\,d\mu$$

[with $(p_t)_{t>0}$ from Theorem 2.2]. In particular, $\mu$ is a stationary measure for the Markov process $\mathbb{M}$ from Theorem 2.3.

PROOF. See the Appendix. □

**3. Finite-dimensional approximation: uniform estimates.** In this section we study finite-dimensional approximation of (1.2)–(1.3). The results will be used in an essential way below.

The main result of this section is Proposition 3.4, giving estimates on the resolvent, including its gradients associated with the approximation $L_N$ of our operator $L$ on $E_N$ [cf. (3.3) below], but these estimates are uniform with respect to $N$. As a preparation, we need several results of which the second (i.e., an appropriate version of a weak maximum principle) is completely standard. Nevertheless, we include the proof for the convenience of the reader.

Below, the background space is the Euclidean space $\mathbb{R}^N$, $N\in\mathbb{N}$, with the Euclidean inner product denoted by $(\cdot,\cdot)$, $dx$ denotes the Lebesgue measure on $\mathbb{R}^N$ and $L^p(\mathbb{R}^N)$, $W^{r,p}_{\mathrm{loc}}(\mathbb{R}^N)$, $r\in\mathbb{N}\cup\{0\}$, $p\in[1,\infty]$ the corresponding $L^p$ and local Sobolev spaces, respectively.

PROPOSITION 3.1. *Let* $A\colon\mathbb{R}^N\to\mathbb{R}^N$ *be a symmetric strictly positive definite linear operator (matrix),* $F\colon\mathbb{R}^N\to\mathbb{R}^N$ *be a bounded measurable vector field,* $\lambda_* := \sup_{x\in\mathbb{R}^N}\frac{(F(x),A^{-1}F(x))}{4}$, $\rho\in L^1(\mathbb{R}^N)$ *be strictly positive and locally Lipschitz and* $W\in L^\infty_{\mathrm{loc}}(\mathbb{R}^N)$, $W\geq 0$. *Let*

$$Lu := \rho^{-1}\operatorname{div}(A\rho Du) + (F,Du) = \operatorname{Tr} AD^2 u + (\rho^{-1}AD\rho + F,Du) - Wu,$$

$$u\in W^{2,1}_{\mathrm{loc}}(\mathbb{R}^N).$$

*Then there exists a unique sub-Markovian pseudo-resolvent* $(\mathcal{R}_\lambda)_{\lambda>0}$ *on* $L^\infty(\mathbb{R}^N)$, *that is, a family of operators satisfying the first resolvent equation, which is Markovian if* $W=0$, *such that:*

(a) $\operatorname{Range}(\mathcal{R}_\lambda)\subset\operatorname{Dom}:=\{u\in\bigcap_{p<\infty}W^{2,p}_{\mathrm{loc}}(\mathbb{R}^N)\,|\,u, Lu\in L^\infty(\mathbb{R}^N)\}$ *and*

$$(\lambda-L)\mathcal{R}_\lambda = \operatorname{id}\qquad \text{for all } \lambda>0.$$

(b) *For all* $\lambda>\lambda_*$ *and* $f\in L^\infty(\mathbb{R}^N)$, *one has* $|D\mathcal{R}_\lambda f|\in L^2(\mathbb{R}^N,\rho\,dx)$.

(c) *For all* $f\in L^\infty(\mathbb{R}^N)$, *one has* $\lim_{\lambda\to\infty}\lambda\mathcal{R}_\lambda f = f$ *in* $L^2(\mathbb{R}^N,\rho\,dx)$.

*Hence, in particular, $\mathcal{R}_\lambda f$ for $f \in L^\infty(\mathbb{R}^N)$ has a continuous $dx$-version, as have its first weak derivatives, and for the continuous versions of $\mathcal{R}_\lambda f$, $\lambda > 0$, the resolvent equation holds pointwise on all of $\mathbb{R}^N$. If both $f$ and $F$ above are in addition locally Lipschitz, then $\mathcal{R}_\lambda f \in \bigcap_{p<\infty} W_{loc}^{3,p}(\mathbb{R}^N)$ for every $\lambda > 0$, hence, its continuous $dx$-version is in $C^2(\mathbb{R}^N)$.*

PROOF. Consider the following bi-linear form $(\mathcal{E}, D(\mathcal{E}))$ in $L^2(\mathbb{R}^N, \rho\, dx)$:

$$\mathcal{E}(u,v) := \int_{\mathbb{R}^N} [(Du, ADv) - (F, Du)v + Wuv]\rho\, dx,$$

$$D(\mathcal{E}) := \left\{ u \in W_{loc}^{1,2}(\mathbb{R}^N) \,\Big|\, \int_{\mathbb{R}^N} [u^2 + |Du|^2 + Wu^2]\rho\, dx < \infty \right\}.$$

Since, for all $u, v \in D(\mathcal{E})$,

(3.1)             $|u(F, Dv)| \leq |(Dv, ADv)| + \lambda_*|u|^2,$

it follows that $\mathcal{E} \geq -\lambda_*$. Then it is easy to show that $(\mathcal{E} + \lambda_*(\cdot, \cdot), D(\mathcal{E}))$ is a Dirichlet form (cf. [43], Section I.4., i.e., a closed sectorial Markovian from) on $L^2(\mathbb{R}^N, \rho\, dx)$. Hence, there exists an associated sub-Markovian strongly continuous resolvent $(R_\lambda)_{\lambda > \lambda_*}$ and semigroup $(P_t)_{t \geq 0}$ on $L^2(\mathbb{R}^N, \rho\, dx)$ (cf. ibid.). Note that $1 \in D(\mathcal{E})$ and $\mathcal{E}(1, v) = 0$ for all $v \in D(\mathcal{E})$, provided $W = 0$, so $(R_\lambda)_{\lambda > \lambda_*}$ and $(P_t)_{t>0}$ are even Markovian in this case. In particular, assertion (b) holds. Note that, for a bounded $f \in L^2(\mathbb{R}^N, \rho\, dx)$, we can define

$$\mathcal{R}_\lambda f := \int_0^\infty e^{-\lambda t} P_t f\, dt$$

even for all $\lambda > 0$ instead of $\lambda > \lambda_*$. Here, the $L^2(\mathbb{R}^N, \rho\, dx)$-valued intregral is taken in the sense of Bochner. Then $\lambda R_\lambda f = \lambda \mathcal{R}_\lambda f \xrightarrow{\lambda \to \infty} f$ in $L^2(\mathbb{R}^N, \rho\, dx)$ and $(\mathcal{R}_\lambda)_{\lambda > 0}$ is a sub-Markovian pseudo-resolvent on $L^\infty(\mathbb{R}^N)$. In particular, the first resolvent equation and assertion (c) hold.

To show (a), we first note that, for $\lambda > \lambda_*$ and $f \in L^\infty(\mathbb{R}^N)$, the bounded function $u := \mathcal{R}_\lambda f$ is a weak solution to the equation

$$\lambda u - Lu = \lambda u - \rho^{-1} \operatorname{div}(A\rho Du) - (F, Du) + Wu = f \qquad \text{in } \mathbb{R}^N.$$

Hence, it follows from [28], Theorem 8.8, that $u \in W_{loc}^{2,2}(\mathbb{R}^N)$. Then [28], Lemma 9.16, yields that $u \in \text{Dom}$. Thus, $\text{Range}(\mathcal{R}_\lambda) \subseteq \text{Dom}$, provided $\lambda > \lambda_*$. Now let $\lambda \in (0, \lambda_*]$. Then for all $\lambda' > \lambda_*$, $\mathcal{R}_\lambda f = \mathcal{R}_{\lambda'} f + (\lambda' - \lambda)\mathcal{R}_{\lambda'}\mathcal{R}_\lambda f$. Hence, $\mathcal{R}_\lambda f \in \text{Dom}$ and $(\lambda' - L)\mathcal{R}_\lambda f = f + (\lambda' - \lambda)\mathcal{R}_\lambda f$. So, $(\lambda - L)\mathcal{R}_\lambda f = f$. The last part follows by Sobolev embedding. $\square$

LEMMA 3.2. *Let $A: \mathbb{R}^N \to \mathbb{R}^N$ be a symmetric strictly positive definite linear operator (matrix), $F: \mathbb{R}^N \to \mathbb{R}^N$ be a bounded measurable vector field,* $\lambda_* := \sup_{x \in \mathbb{R}^N} \frac{(F(x), A^{-1}F(x))}{4}$, $\rho > 0$ be locally Lipschitz and $W \in L^\infty_{\mathrm{loc}}(\mathbb{R}^N)$, $W \geq 0$.

For $\lambda > \lambda_*$, let $u \in W^{1,2}_{\mathrm{loc}}(\mathbb{R}^N) \cap L^2(\mathbb{R}^N, \rho\,dx)$ be a weak super-solution to the equation

$$\lambda u - \rho^{-1} \operatorname{div}(A\rho Du) - (F, Du) + Wu = 0 \qquad on \ \mathbb{R}^N$$

[i.e., a weak solution to the inequality $\lambda u - \rho^{-1} \operatorname{div}(A\rho Du) - (F, Du) + Wu \geq 0$].

Then $u \geq 0$.

PROOF.    For $\theta \in C^1_c(\mathbb{R}^N)$, choose $u^- \theta^2 \rho$ as a test function. Then, using the fact that $u^+ \wedge u^- = 0$, we obtain that, for all $\varepsilon > 0$,

$$0 \leq -\int [(\lambda + W)(u^- \theta)^2 + (D(u^- \theta^2), ADu^-) - u^- \theta^2 (F, Du^-)]\rho\,dx$$

$$= -\int [(\lambda + W)(u^- \theta)^2 + (D(u^- \theta), AD(u^- \theta)) - u^- \theta (F, D(u^- \theta))]\,\rho\,dx$$

$$\quad + \int (u^-)^2 [(D\theta, AD\theta) - \theta(F, D\theta)]\rho\,dx$$

$$\leq -\int (\lambda - (1+\varepsilon)\lambda_*)(u^- \theta)^2 \rho\,dx + \int (u^-)^2 \left(1 + \frac{1}{\varepsilon}\right)(D\theta, AD\theta)\rho\,dx,$$

where we used the fact that $W \geq 0$, $\mathcal{E} \geq -\lambda_*$ and we applied (3.1) with $\varepsilon\theta$, $\frac{1}{\varepsilon}\theta$ replacing $u$, $v$, respectively. Hence, for all $\varepsilon > 0$,

$$(\lambda - (1+\varepsilon)\lambda_*) \int (u^- \theta)^2\,dx \leq \left(1 + \frac{1}{\varepsilon}\right) \int (u^-)^2 (D\theta, AD\theta)\rho\,dx.$$

Now we choose $\varepsilon > 0$ such that $\lambda > (1+\varepsilon)\lambda_*$ and let $\theta \nearrow 1$ and $D\theta \to 0$ such that $(D\theta, AD\theta) \leq C_A$. Then the dominated convergence theorem yields $u^- = 0$.  □

COROLLARY 3.3.    Let $A : \mathbb{R}^N \to \mathbb{R}^N$ be a symmetric strictly positive definite linear operator (matrix), $F : \mathbb{R}^N \to \mathbb{R}^N$ be a bounded measurable vector field, $\lambda_* := \sup_{x \in \mathbb{R}^N} \frac{(F(x), A^{-1}F(x))}{4}$, $\rho > 0$ be locally Lipschitz and $W \in L^\infty_{\mathrm{loc}}(\mathbb{R}^N)$.

Let $V \in C^2(\mathbb{R}^N)$, $V \geq 1$ be such that, for some $\lambda_V \in \mathbb{R}$,

$$(3.2) \qquad \lambda_V V - \rho^{-1} \operatorname{div}(A\rho DV) - (F, DV) + WV \geq 0.$$

Let $f \in L^2(\mathbb{R}^N, \rho\,dx)$, $V^{-1}f \in L^\infty(\mathbb{R}^N)$, $\lambda > \lambda_* + \lambda_V$ and $u \in W^{1,2}_{\mathrm{loc}}(\mathbb{R}^N) \cap L^2(\mathbb{R}^N, \rho\,dx)$ be a weak sub-solution to the equation

$$\lambda u - \rho^{-1} \operatorname{div}(A\rho Du) - (F, Du) + Wu = f \qquad on \ \mathbb{R}^N$$

[*i.e., a weak solution to the inequality* $\lambda u - \rho^{-1} \operatorname{div}(A\rho Du) - (F, Du) + Wu \leq f$]. *Then*

$$\|V^{-1}u\|_\infty \leq \frac{1}{\lambda - \lambda_V}\|V^{-1}f\|_\infty.$$

PROOF.   Let $\tilde{W} := V^{-1}[\lambda_V V - \rho^{-1}\operatorname{div}(A\rho DV) - (F, DV) + WV]$ and $v := V^{-1}u$. It is easy to see that $v \in W^{1,2}_{\text{loc}}(\mathbb{R}^N) \cap L^2(\mathbb{R}^N, V^2\rho\,dx)$ and it is a weak sub-solution to the equation

$$(\lambda - \lambda_V)v - \frac{1}{V^2\rho}\operatorname{div}(AV^2\rho Dv) - (F, Dv) + \tilde{W}v = V^{-1}f \qquad \text{on } \mathbb{R}^N.$$

Note that $V^{-1}f \in L^2(\mathbb{R}^N, V^2\rho\,dx)$. Since $\tilde{W} \geq 0$, the result now follows from Lemma 3.2 and the fact that the resolvent associated on $L^2(\mathbb{R}^N, V^2\rho\,dx)$ with the bi-linear form

$$\mathcal{E}(g, h) := \int_{\mathbb{R}^N}[(Dg, ADh) - (F, Dg)h + \tilde{W}gh]V^2\rho\,dx,$$

$$D(\mathcal{E}) := \left\{ g \in W^{1,2}_{\text{loc}}(\mathbb{R}^N)\,\Big|\, \int_{\mathbb{R}^N}[g^2 + |Dg|^2 + \tilde{W}g^2]V^2\rho\,dx < \infty \right\},$$

is sub-Markovian.   □

PROPOSITION 3.4.   *Let* $A, H : \mathbb{R}^N \to \mathbb{R}^N$ *be symmetric strictly positive definite linear operators (matrices) such that* $AH = HA$. *Let* $F : \mathbb{R}^N \to \mathbb{R}^N$ *be a bounded locally Lipschitz vector field. Let*

$$Lu(x) := \operatorname{Tr}(AD^2u)(x) + (-Hx + F(x), Du(x)),$$

(3.3)
$$u \in W^{2,1}_{\text{loc}}(\mathbb{R}^N), x \in \mathbb{R}^N.$$

*Let* $\Gamma : \mathbb{R}^N \to \mathbb{R}^N$ *be a symmetric nondegenerate linear operator (matrix) such that* $\Gamma H = H\Gamma$. *Assume the following:*

(i)  *there exists* $V_0 \in C^2(\mathbb{R}^N)$, $V_0 \geq 1$ *and* $\lambda_{V_0} \in \mathbb{R}$ *such that*

(3.4)
$$(\lambda_{V_0} - L)V_0 \geq 0;$$

(ii)  *there exists* $V_1 \in C^2(\mathbb{R}^N)$, $V_1 \geq 1$ *and* $\lambda_{V_1} \in \mathbb{R}$ *such that*

$$(\lambda_{V_1} - L - W)V_1 \geq 0$$

(3.5)
$$\text{with } W(x) := \sup_{|y|=1}[(DF(x)\Gamma y, \Gamma^{-1}y) - |H^{1/2}y|^2], x \in \mathbb{R}^N.$$

*Then:*

(i) *there exists a unique Markovian pseudo-resolvent* $(\mathcal{R}_\lambda)_{\lambda>0}$ *on* $L^\infty(\mathbb{R}^N)$ *such that*

$$\operatorname{Range}(\mathcal{R}_\lambda) \subset \left\{ u \in \bigcap_{p<\infty} W^{2,p}_{\mathrm{loc}}(\mathbb{R}^N) \,|\, u, Lu \in L^\infty(\mathbb{R}^N) \right\},$$

$(\lambda - L)\mathcal{R}_\lambda = \mathrm{id}$ *for all* $\lambda > 0$, *and* $\lambda\mathcal{R}_\lambda f \to f$ *as* $\lambda \to \infty$ *pointwise on* $\mathbb{R}^N$ *for bounded locally Lipschitz* $f$;

(ii) *for a bounded locally Lipschitz* $f$, *we have*

$$\tag{3.6} \|V_0^{-1}\mathcal{R}_\lambda f\|_\infty \leq \frac{1}{\lambda - \lambda_{V_0}} \|V_0^{-1}f\|_\infty$$

*for all* $\lambda > \lambda_{V_0}$; *and*

$$\tag{3.7} \sup_x V_1^{-1}|\Gamma D\mathcal{R}_\lambda f|(x) \leq \frac{1}{\lambda - \lambda_{V_1}} \operatorname*{essup}_x V_1^{-1}|\Gamma Df|(x)$$

*for all* $\lambda > \lambda_{V_1}$ *provided* $V_1^{-1}|Df| \in L^\infty(\mathbb{R}^N)$ *and* $|Df| \in L^2(\mathbb{R}^N, \rho\,dx)$. *Here* $D\mathcal{R}_\lambda f$ *and* $\mathcal{R}_\lambda f$ *denote the (unique) continuous* $dx$-*versions of* $D\mathcal{R}_\lambda f$, $\mathcal{R}_\lambda f$, *respectively, which exist by assertion* (i) *and* $\rho(x) := \exp\{-\frac{1}{2}(x, A^{-1}Hx)\}$, $x \in \mathbb{R}^N$.

To prove Proposition 3.4, we need another lemma.

LEMMA 3.5. *Let* $A, H : \mathbb{R}^N \to \mathbb{R}^N$ *be symmetric strictly positive definite linear operators (matrices) such that* $AH = HA$. *Let* $F : \mathbb{R}^N \to \mathbb{R}^N$ *be a bounded locally Lipschitz vector field. Let* $L$ *be defined as in* (3.3).

*Let* $\Gamma : \mathbb{R}^N \to \mathbb{R}^N$ *be a symmetric nondegenerate linear operator (matrix) such that* $\Gamma H = H\Gamma$.

*Let* $\lambda \in \mathbb{R}$, $f$ *be locally Lipschitz and* $u \in W^{1,2}_{\mathrm{loc}}(\mathbb{R}^N)$ *be a weak solution to the equation* $(\lambda - L)u = f$ *on* $\mathbb{R}^N$.

*Then* $u \in \bigcap_{p<\infty} W^{2,p}_{\mathrm{loc}}(\mathbb{R}^N)$ *and* $v := |\Gamma Du|$ *is a weak sub-solution to the equation*

$$(\lambda - L - W)v = |\Gamma Df|$$

*with* $W(x) := \sup_{|y|=1} [(DF(x)\Gamma y, \Gamma^{-1}y) - |H^{1/2}y|^2], x \in \mathbb{R}^N.$

PROOF. Throughout the proof let $\langle f, g \rangle$ stand for $\int_{\mathbb{R}^N} f(x)g(x)\,dx$ or $\int_{\mathbb{R}^N} (f(x), g(x))\,dx$ whenever $fg \in L^1(\mathbb{R}^N, dx)$ or $(f, g) \in L^1(\mathbb{R}^N, dx)$, for $f, g : \mathbb{R}^N \to \mathbb{R}$ or $f, g : \mathbb{R}^N \to \mathbb{R}^N$ measurable, $\eta_m$, $m = 1, \ldots, N$ be the (common) orthonormal eigenbasis for $H$ and $\Gamma$, $\Gamma\eta_m = \gamma_m\eta_m$, $m = 1, \ldots, N$.

By [28], Theorem 8.8 and Lemma 9.16, $u \in \bigcap_{p<\infty} W^{2,p}_{\mathrm{loc}}(\mathbb{R}^N)$. For $m = 1, \ldots, N$, let $u_m := \partial_m u$. Then, for a bounded $\phi \in W^{1,2}_c(\mathbb{R}^N)$, integration by parts yields

$$\langle Du_m, AD\phi \rangle = -\langle \lambda u_m - \partial_m f - (-H + F, Du_m) - (-H\eta_m + \partial_m F, Du), \phi \rangle.$$

Set $[Du] := |\Gamma Du|$ and, for $\varepsilon > 0$, $[Du]_\varepsilon := \sqrt{|\Gamma Du|^2 + \varepsilon}$. For $\theta \in C^\infty_c(\mathbb{R}^N)$ and $m = 1, \ldots, \dim E$, choose $\phi_m := \frac{|\gamma_m|^2 u_m}{[Du]_\varepsilon} \theta$. Then $\phi_m$ is bounded and

$$D\phi_m = \frac{|\gamma_m|^2 u_m}{[Du]_\varepsilon} D\theta + \frac{|\gamma_m|^2 Du_m}{[Du]_\varepsilon} \theta - \frac{|\gamma_m|^2 u_m D^2 u \Gamma^2 Du}{[Du]^3_\varepsilon} \theta$$

$$\text{with } |D\phi_m| \in \bigcap_{p<\infty} L^p(\mathbb{R}^N, dx).$$

Hence, a.e. on $\mathbb{R}^N$,

$$\sum_m (Du_m, AD\phi_m)$$
$$= ([Du]_\varepsilon, AD\theta)$$
$$+ \left[ \mathrm{Tr}\{\Gamma D^2 u A D^2 u \Gamma\} - \left( \frac{\Gamma Du}{[Du]_\varepsilon}, \Gamma D^2 u A D^2 u \Gamma \frac{\Gamma Du}{[Du]_\varepsilon} \right) \right] \frac{\theta}{[Du]_\varepsilon}.$$

Since $D([Du]_\varepsilon) = \frac{D^2 u \Gamma^2 Du}{[Du]_\varepsilon}$, it follows that $v_\varepsilon := [Du]_\varepsilon$ is a weak solution of the equation $(\lambda - L - W_\varepsilon)v = G_\varepsilon$, where

$$W_\varepsilon = \frac{1}{[Du]^2_\varepsilon}(-|H^{1/2}\Gamma Du|^2_2 + (\Gamma Du, \Gamma(DF)^t Du))$$

and

$$G_\varepsilon := \lambda \frac{\varepsilon}{[Du]_\varepsilon} + \left( \frac{\Gamma Du}{[Du]_\varepsilon}, \Gamma Df \right)$$
$$- \frac{1}{[Du]_\varepsilon}\left[ \mathrm{Tr}\{\Gamma D^2 u A D^2 u \Gamma\} - \left( \frac{\Gamma Du}{[Du]_\varepsilon}, \Gamma D^2 u A D^2 u \Gamma \frac{\Gamma Du}{[Du]_\varepsilon} \right) \right].$$

We have $W_\varepsilon \leq W$ a.e. so $v_\varepsilon$ is a weak sub-solution to the equation $(\lambda - L - W)v = G_\varepsilon$. Passing to the limit as $\varepsilon \to 0$, we see that $v_\varepsilon = [Du]_\varepsilon$ converges to $v = [Du]$ in $W^{1,2}_{\mathrm{loc}}(\mathbb{R}^N)$ and, thus, the assertion follows. $\square$

PROOF OF PROPOSITION 3.4. Note that, provided $AH = HA$, we have

$$Lu = \rho^{-1} \mathrm{div}(A\rho Du) + (F, Du),$$

where $\rho(x) = \exp\{-\frac{1}{2}(x, A^{-1}Hx)\}$. Hence, Proposition 3.1(a) implies assertion (i) except for the fact that $\lambda \mathcal{R}_\lambda f \to f$ pointwise as $\lambda \to \infty$, which we shall prove at the end.

Let $f : \mathbb{R}^N \to \mathbb{R}$ be bounded and locally Lipschitz and such that

$$V_1^{-1} |\Gamma Df| \in L^\infty(\mathbb{R}^N) \quad \text{and} \quad u := \mathcal{R}_\lambda f.$$

By assertion (i), $u$ is a weak solution of the equation $(\lambda - L)u = f$ on $\mathbb{R}^N$ and, by Lemma 3.5, $v := |\Gamma Du|$ is a weak sub-solution to the equation $(\lambda - L - W)v = |\Gamma Df|$ on $\mathbb{R}^N$ with $W$ as in Lemma 3.5. Let first $\lambda > \lambda_* + \lambda_{V_0} \vee \lambda_{V_1}$. Note that $u \in L^2(\mathbb{R}^N, \rho\, dx)$ and $v \in L^2(\mathbb{R}^N, \rho\, dx)$ by Proposition 3.1(b). Then (3.6)–(3.7) follow from assumptions (i) and (ii) and Corollary 3.3, since $f, |\Gamma Df| \in L^2(\mathbb{R}^N, \rho\, dx)$.

By density, for $\lambda > \lambda_* + \lambda_{V_0} \vee \lambda_{V_1}$, the operator $\mathcal{R}_\lambda$ can be continuously extended to the completion of the bounded locally Lipschitz functions on $\mathbb{R}^N$ with respect to $\|V_0^{-1} \cdot\|_\infty$, preserving the resolvent identity and estimate (3.6). Moreover, for a locally Lipschitz $f$ such that $V_0^{-1} f, V_1^{-1}|\Gamma Df| \in L^\infty(\mathbb{R}^N)$ and $|Df| \in L^2(\mathbb{R}^N, \rho\, dx)$, estimate (3.7) holds. This is easy to see by replacing $f$ by $(f \vee (-v)) \wedge n$ and letting $n \to \infty$. Now, for $\lambda \in (\lambda_{V_0}, \lambda_* + \lambda_{V_0} \vee \lambda_{V_1}]$, one can define

$$\tag{3.8} \mathcal{R}_\lambda = \sum_{k=1}^\infty (\lambda_0 - \lambda)^{k-1} \mathcal{R}_{\lambda_0}^k$$

with some $\lambda_0 > \lambda_* + \lambda_{V_0} \vee \lambda_{V_1}$. The series converges in operator norm due to (3.6) and (3.6) is preserved:

$$\|V_0^{-1} \mathcal{R}_\lambda f\|_\infty \leq \sum_{k=1}^\infty (\lambda_0 - \lambda)^{k-1} \|V_0^{-1} \mathcal{R}_{\lambda_0}^k f\|_\infty \leq \sum_{k=1}^\infty \frac{(\lambda_0 - \lambda)^{k-1}}{(\lambda_0 - \lambda_{V_0})^k} \|V_0^{-1} f\|_\infty$$

$$= \frac{1}{\lambda - \lambda_{V_0}} \|V_0^{-1} f\|_\infty.$$

On $L^\infty(\mathbb{R}^N)$, obviously $\mathcal{R}_\lambda$ defined in (3.8) coincides with $\mathcal{R}_\lambda$ defined in Proposition 3.1 with $W = 0$. So, $\lambda \mathcal{R}_\lambda$ remains Markovian for $\lambda \in (\lambda_{V_0}, \lambda_* + \lambda_{V_0} \vee \lambda_{V_1}]$. By similar arguments, using the closability of $\Gamma D$, we prove that (3.7) is preserved for $\lambda \in (\lambda_{V_1}, \lambda_* + \lambda_{V_0} \vee \lambda_{V_1}]$.

We are left to prove that $\lambda \mathcal{R}_\lambda f \to f$ pointwise on $\mathbb{R}^N$ as $\lambda \to \infty$, for any bounded locally Lipschitz $f$. The proof is by contradiction. Let $x_0 \in \mathbb{R}^N$ and some subsequence $\lambda_n \to \infty$ and some $\varepsilon \in (0, 1]$,

$$\tag{3.9} |\lambda_n \mathcal{R}_{\lambda_n} f(x_0) - f(x_0)| > \varepsilon \qquad \forall n \in \mathbb{N}.$$

Selecting another subsequence if necessary, by Proposition 3.1(c), we may assume that the complement of the set

$$M := \left\{ x \in \mathbb{R}^N \,\Big|\, \lim_{n \to \infty} \lambda_n \mathcal{R}_{\lambda_n} f(x) = f(x) \right\}$$

in $\mathbb{R}^N$ has Lebesgue measure zero, so $M$ is dense in $\mathbb{R}^N$. By (3.7), the sequence $(\lambda_n \mathcal{R}_{\lambda_n} f)_{n \in \mathbb{N}}$ is equicontinuous and converges on the dense set $M$ to the continuous function $f$, hence, it must converge everywhere on $\mathbb{R}^N$ to $f$. This contradicts (3.9).  □

LEMMA 3.6.  *Let $V\colon\mathbb{R}^N\to[1,\infty)$ be convex (hence, continuous) and let $\Gamma\colon\mathbb{R}^N\to\mathbb{R}^N$ be a symmetric invertible linear operator (matrix) and $f\colon\mathbb{R}^N\to\mathbb{R}$ be locally Lipschitz. Then*

$$(3.10) \qquad \|V^{-1}|\Gamma Df|\|_\infty = \sup_{y_1,y_2\in\mathbb{R}^N} \frac{1}{V(y_1)\vee V(y_2)} \frac{|f(y_1)-f(y_2)|}{|\Gamma^{-1}(y_1-y_2)|}.$$

PROOF.  We may assume that $f\in C^1(\mathbb{R}^N)$. The general case follows by approximation. Let $x\in\mathbb{R}^N$. Then we have

$$\frac{1}{V(x)}|\Gamma Df(x)| = \lim_{\substack{y_1,y_2\to x \\ y_1,y_2\in\mathbb{R}^N}} \frac{1}{V(y_1)\vee V(y_2)} \frac{|f(y_1)-f(y_2)|}{|\Gamma^{-1}(y_1-y_2)|}.$$

On the other hand, for $y_1,y_2\in\mathbb{R}^N$,

$$\frac{1}{V(y_1)\vee V(y_2)} \frac{|f(y_1)-f(y_2)|}{|\Gamma^{-1}(y_1-y_2)|}$$
$$= \frac{1}{V(y_1)\vee V(y_2)}$$
$$\times \left| \int_0^1 (\Gamma Df(\tau y_1+(1-\tau)y_2), \Gamma^{-1}(y_2-y_1)|\Gamma^{-1}(y_2-y_1)|^{-1})\,d\tau \right|$$
$$\le \|V^{-1}|\Gamma Df|\|_\infty,$$

where we used that $V(\tau y_1+(1-\tau)y_2)\le V(y_1)\vee V(y_2)$, since $V$ is convex. Hence, the assertion follows.  □

REMARK 3.7.  We note that if the right-hand side of (3.10) is finite, then $f$ is Lipschitz on the level sets of $V$.

## 4. Approximation and condition (F2).

In this section we construct a sequence $F_N\colon E_N\to E_N$, $N\in\mathbb{N}$, of bounded locally Lipschitz continuous vector fields approximating the nonlinear drift $F$. The corresponding operators $L_N$, $N\in\mathbb{N}$, are of the form

$$(4.1) \qquad L_N u(x) := \tfrac{1}{2}\operatorname{Tr}(A_N D^2 u)(x) + (x''+F_N(x),Du(x)),$$
$$u\in W^{2,1}_{\mathrm{loc}}(E_N), x\in E_N, N\in\mathbb{N},$$

whose resolvents $(G_\lambda^{(N)})_{\lambda>0}$, lifted to $X_p$, will be shown in Section 6 to converge weakly to the resolvent of $L$.

We introduce the following condition for a map $F\colon H_0^1\to X$:

(F2) For every $k \in \mathbb{N}$, the map $F^{(k)} := (F, \eta_k) : H_0^1 \to \mathbb{R}$ is $|\cdot|_2$-continuous on $|\cdot|_{1,2}$-balls and there exists a sequence $F_N : E_N \to E_N$, $N \in \mathbb{N}$, of bounded locally Lipschitz continuous vector fields satisfying the following conditions:

  (F2a) There exist $\kappa_0 \in (0, \frac{1}{4a_0}]$ and a set $Q_{\text{reg}} \subset [2, \infty)$ such that $2 \in Q_{\text{reg}}$ and for all $\kappa \in (0, \kappa_0)$, $q \in Q_{\text{reg}}$, there exist $m_{q,\kappa} > 0$ and $\lambda_{q,\kappa} \in \mathbb{R}$ such that for all $N \in \mathbb{N}$,

$$(4.2) \qquad L_N V_{q,\kappa} := L_N(V_{q,\kappa}\restriction_{E_N}) \leq \lambda_{q,\kappa} V_{q,\kappa} - m_{q,\kappa} \Theta_{q,\kappa} \qquad \text{on } E_N.$$

  (F2b) For all $\varepsilon \in (0,1)$, there exists $C_\varepsilon \in (0, \infty)$ such that for all $N \in \mathbb{N}$ and $dx$-a.e. $x \in E_N$ (where $dx$ denotes Lebesgue measure on $E_N$)

$$(DF(x)y, y) \leq |y'|_2^2 + (\varepsilon|x'|_2^2 + C_\varepsilon)|y|_2^2 \qquad \forall y \in E_N.$$

  (F2c) $\lim_{N \to \infty} |P_N F - F_N \circ P_N|_2(x) = 0 \ \forall x \in H_0^1$.

  (F2d) For $\kappa_0$ and $Q_{\text{reg}}$ as in (F2a), there exist $\kappa \in (0, \kappa_0)$, $p \in Q_{\text{reg}}$ such that, for some $C_{p,\kappa} > 0$ and some $\omega : [0, \infty) \to [0, 1]$ vanishing at infinity,

$$|F_N \circ P_N|_2(x) \leq C_{p,\kappa} \Theta_{p,\kappa}(x) \omega(\Theta_{p,\kappa}(x)) \qquad \forall x \in H_0^1, N \in \mathbb{N}.$$

Furthermore, we say that condition (F2+) holds if, in addition, to (F2) we have:

(F2e) For all $\varepsilon \in (0, 1)$, there exists $C_\varepsilon \in (0, \infty)$ such that, for all $N \in \mathbb{N}$ and $dx$-a.e. $x \in E_N$,

$$(DF_N(x)(-\Delta)^{1/2}y, (-\Delta)^{-1/2}y) \leq |y'|_2^2 + (\varepsilon|x'|_2^2 + C_\varepsilon)|y|_2^2$$
$$\forall y \in E_N.$$

The main result of this section is the following:

PROPOSITION 4.1. *Let $F$ be as in* (2.15) *and let assumptions* $(\Psi)$, $(\Phi 1)$–$(\Phi 3)$ *be satisfied. Then* (F2) *holds. More precisely,* (F2a) *holds with* $\kappa_0 := \frac{2 - |h_1|_1}{8 a_0}$, $Q_{\text{reg}} := [2, \infty)$, (F2c) *holds uniformly on $H_0^1$-balls, and* (F2d) *holds with $p \in [2, \infty) \cap (q_2 - 3 + \frac{2}{q_1}, \infty)$ and any $\kappa \in (0, \kappa_0)$. If, in addition,* $(\Phi 4)$ *is satisfied, then* (F2+) *holds.*

To prove our main results formulated in Section 2, we shall only use conditions (F2), (F2+), respectively. Before we prove Proposition 4.1, as a motivation, we shall prove that (F2) [in fact, even only (F2a)–(F2c)] and (F2e) will imply regularity and convergence (see also Theorem 6.4 below) of the above mentioned resolvents $(G_\lambda^N)_{\lambda > 0}$.

COROLLARY 4.2. *Let* (A) *and* (F2a)–(F2c) *hold and let* $L_N$ *be as in* (4.1) *with* $F_N$ *as in* (F2). *Let* $(\mathcal{R}_\lambda^{(N)})_{\lambda>0}$ *be the corresponding Markovian pseudo-resolvent on* $L^\infty(E_N)$ *from Proposition* 3.1. *For a bounded Borel measurable* $f : X \to \mathbb{R}$, *we define*

$$G_\lambda^{(N)} f := (\mathcal{R}_\lambda^{(N)}(f{\restriction_{E_N}})) \circ P_N.$$

*Then* $\lambda G_\lambda^{(N)}$ *is Markovian and* $\lambda G_\lambda^{(N)} f \to f \circ P_N$ *pointwise as* $\lambda \to \infty$ *for all bounded* $f$ *which are locally Lipschitz on* $E_N$.

*Let* $\kappa_0$, $Q_{\mathrm{reg}}$ *be as in* (F2a) *and let* $\kappa \in (0, \kappa_0)$, $q \in Q_{\mathrm{reg}}$ *with* $\lambda_{q,\kappa}$ *as in* (F2a). *Set* $\lambda'_{q,\kappa} := \lambda_{q,\kappa} + C_{m_{q,\kappa}}$, *with* $m_{q,\kappa}$ *as in* (F2a) *and function* $\varepsilon \mapsto C_\varepsilon$ *as in* (F2b). *Let* $N \in \mathbb{N}$ *and* $f \in \mathrm{Lip}_{0,q,\kappa}$, $f$ *bounded. Then*

$$(4.3) \quad |G_\lambda^{(N)} f(x)| \le \frac{1}{\lambda - \lambda_{q,\kappa}} V_{q,\kappa}(P_N x) \|f\|_{q,\kappa}, \qquad x \in X_q, \lambda > \lambda_{q,\kappa},$$

*and for* $y_1, y_2 \in X_q$,

$$(4.4) \quad \frac{|G_\lambda^{(N)} f(y_1) - G_\lambda^{(N)} f(y_2)|}{|y_1 - y_2|_2} \le \frac{|G_\lambda^{(N)} f(y_1) - G_\lambda^{(N)} f(y_2)|}{|P_N(y_1 - y_2)|_2}$$
$$\le \frac{V_{q,\kappa}(P_N y_1) \vee V_{q,\kappa}(P_N y_2)}{\lambda - \lambda'_{q,\kappa}}(f)_{0,q,\kappa},$$
$$\lambda > \lambda'_{q,\kappa}.$$

*In particular, if* $\lambda > \lambda'_{q,\kappa} \vee \lambda_{q,\kappa}$, *then* $G_\lambda^{(N)} f \in \bigcap_{\varepsilon>0} \mathcal{D}_{q,\kappa+\varepsilon}$ *and, provided* $f \in \mathcal{D}$, $G_\lambda^{(N)} f \in \bigcap_{\varepsilon>0} \mathcal{D}_\varepsilon$. *Furthermore, for all* $x \in H_0^1, \lambda > \lambda'_{q,\kappa}$,

$$(4.5) \quad |(\lambda - L)G_\lambda^{(N)} f(x) - (f \circ P_N)(x)|$$
$$\le \frac{1}{\lambda - \lambda'_{q,\kappa}} |P_N F - F_N \circ P_N|_2(x)\alpha_q^q V_{q,\kappa}(x)(f)_{0,q,\kappa}.$$

*In particular, for all* $\lambda_* > \lambda'_{q,\kappa}$,

$$(4.6) \quad \lim_{m \to \infty} \sup_{\lambda \ge \lambda_*} \lambda|(\lambda - L)G_\lambda^{(m)} f - f|(x) = 0 \qquad \forall x \in H_0^1.$$

*If, moreover,* (F2e) *holds, let* $\lambda''_{q,\kappa} := \lambda_{q,\kappa} + C_{m_{q,\kappa}}$, *with* $m_{q,\kappa}$ *as in* (F2a) *and function* $\varepsilon \mapsto C_\varepsilon$ *as in* (F2e). *Then, for* $N \in \mathbb{N}$ *and* $f \in \mathrm{Lip}_{1,q,\kappa}$, $f$ *bounded, we have, for* $y_1, y_2 \in X_q$,

$$(4.7) \quad \frac{|G_\lambda^{(N)} f(y_1) - G_\lambda^{(N)} f(y_2)|}{|(-\Delta)^{-1/2}(y_1 - y_1)|_2} \le \frac{V_{q,\kappa}(P_N y_1) \vee V_{q,\kappa}(P_N y_2)}{\lambda - \lambda''_{q,\kappa}}(f)_{1,q,\kappa}.$$

PROOF. To prove (4.3), (4.4) and (4.7), fix $x \in X_q$. By (F2a), we can apply Proposition 3.4 with $V_0 := V_{q,\kappa} \upharpoonright_{E_N}$ to conclude that, for $\lambda > \lambda_{q,\kappa}$,

$$|G_\lambda^{(N)} f(x)| = |\mathcal{R}_\lambda^{(N)}(f \upharpoonright_{E_N})(P_N x)|$$

$$\leq \frac{1}{\lambda - \lambda_{q,\kappa}} V_{q,\kappa}(P_N x) \sup_{y \in E_N} V_{q,\kappa}^{-1}(y) |f(y)|$$

$$\leq \frac{1}{\lambda - \lambda_{q,\kappa}} V_{q,\kappa}(P_N x) \sup_{y \in X_q} V_{q,\kappa}^{-1}(y) |f(y)|,$$

which proves (4.3). By (F2a), (F2b), respectively, (F2a), (F2e), we can apply Proposition 3.4 with $V_1 := V_{q,\kappa} \upharpoonright_{E_N}$ to conclude that, for $\lambda > \lambda_0 := \lambda'_{q,\kappa}$ or $\lambda''_{q,\kappa}$ if $l := 0$, respectively, $l := 1$ and all $y_1, y_2 \in X_q$,

$$\frac{|G_\lambda^{(N)} f(y_1) - G_\lambda^{(N)} f(y_2)|}{|(-\Delta)^{-l/2}(y_1 - y_2)|_2}$$

$$\leq \frac{|\mathcal{R}_\lambda^{(N)}(f \upharpoonright_{E_N})(P_N y_1) - \mathcal{R}_\lambda^{(N)}(f \upharpoonright_{E_N})(P_N y_2)|}{|(-\Delta)^{-l/2}(P_N y_1 - P_N y_2)|_2}$$

$$\leq V_{q,\kappa}(P_N y_1) \vee V_{q,\kappa}(P_N y_2) \sup_{y \in E_N} V_{q,\kappa}^{-1}(y) |(-\Delta)^{l/2}(D\mathcal{R}_\lambda^{(N)}(f \upharpoonright_{E_N})(y))|_2$$

$$\leq \frac{V_{q,\kappa}(P_N y_1) \vee V_{q,\kappa}(P_N y_2)}{\lambda - \lambda_0} (f)_{l,q,\kappa},$$

where we used both Proposition 3.4 and Lemma 3.6 in the last two steps. We note that, by our assumption on $\kappa_0$ in (F2a), we really have that $|Df \upharpoonright_{E_N}| \in L^2(E_N, \rho \, dx)$, so the conditions to have (3.7) are indeed fulfilled.

By the last part of Proposition 3.1, we have that $u := \mathcal{R}_\lambda^{(N)} f \upharpoonright_{E_N} \in C^2(E_N)$ and that

$$(4.8) \qquad \lambda u(x) - L_N u(x) = f(x) \qquad \forall \, x \in E_N.$$

Hence, it follows from (4.3), (4.4), Lemma 3.6 and (2.8) that $G_\lambda^{(N)} f \in \bigcap_{\varepsilon > 0} \mathcal{D}_{p,\kappa+\varepsilon}$ and, provided $f \in \mathcal{D}$, that $G_\lambda^{(N)} f \in \bigcap_{\varepsilon > 0} \mathcal{D}_\varepsilon$. Furthermore, (4.8) implies that, on $H_0^1$,

$$|(\lambda - L)((\mathcal{R}_\lambda^{(N)} f \upharpoonright_{E_N}) \circ P_N) - f \circ P_N|$$

$$= |(P_N F - F_N \circ P_N, D(\mathcal{R}_\lambda^{(N)} f \upharpoonright_{E_N}) \circ P_N)|$$

$$\leq \frac{1}{\lambda - \lambda'_{q,\kappa}} |P_N F - F_N \circ P_N|_2 (f)_{0,q,\kappa} V_{q,\kappa} \circ P_N,$$

where we used (4.4) and Lemma 3.6. Now (4.5) follows by (2.8) and (4.6) follows by (F2c). $\quad \square$

Now we turn to the proof of Proposition 4.1, which will be the consequence of a number of lemmas which we state and prove first.

In the rest of this section, $\phi\colon(0,1)\times\mathbb{R}\to\mathbb{R}$ will be a function square integrable in the first variable locally uniformly in the second and continuous in the second variable, and $\psi\in C^1(\mathbb{R})$. For such functions, we define

$$(4.9)\qquad F_\phi(x):=\phi(\cdot,x(\cdot)),\qquad G_\psi(x):=x'\psi'\circ x,\qquad x,y\in H_0^1.$$

Note that $F_\phi\colon H_0^1\to X$ and $G_\psi\colon H_0^1\to X$.

LEMMA 4.3. *Let $\psi$ satisfy* $(\Psi)$, *and* $\theta\in C_c^\infty(-1,1)$, $0\le\theta\le1$, $\theta\restriction_{[-1/2,1/2]}\equiv1$. *For $N\in\mathbb{N}$, let* $\psi^{(N)}(x):=\psi(x)\theta(\frac{x}{N})$, $x\in\mathbb{R}$.

*Then for $N\in\mathbb{N}$, $\psi^{(N)}\in C_c^{1,1}(\mathbb{R})$ satisfying* $(\Psi)$ *uniformly in $N$, that is, with some $\hat{C}\ge0$ and $\hat{\omega}\colon\mathbb{R}_+\to\mathbb{R}_+$, $\hat{\omega}(r)\to0$ as $r\to\infty$, independent of $N$. Moreover, $|G_{\psi^{(N)}}-G_\psi|_2\to0$ as $N\to\infty$ uniformly on balls in $H_0^1$.*

PROOF. Let, for $x\in\mathbb{R}$, $\theta_1(x):=x\theta'(x)$ and $\theta_2(x):=x^2\theta''(x)$. Then $\psi_{xx}^{(N)}(x)=\psi_{xx}(x)\theta(\frac{x}{N})+2\frac{\psi_x(x)}{x}\theta_1(\frac{x}{N})+\frac{\psi(x)}{x^2}\theta_2(\frac{x}{N})$. Hence, the first assertion follows from Remark 2.1(i).

Note that $\psi^{(N)}(x)=\psi(x)$ whenever $|x|\le\frac{N}{2}$. Hence, the second assertion follows. $\square$

LEMMA 4.4. *Let $\theta\in C^\infty(\mathbb{R})$, odd, $0\le\theta'\le1$, $\theta(x)=x$ for $x\in[-1,1]$ and $\theta(x)=\frac{3}{2}\operatorname{sign}(x)$ for $x\in\mathbb{R}\setminus[-2,2]$.*

*For $N\in\mathbb{N}$, let $\theta_N(x):=N\theta(N^{-1}x)$, $x\in\mathbb{R}$ and $\phi_N:=\theta_N\circ\phi$.*

*Then for all $N\in\mathbb{N}$, $\phi_N$ is a bounded function.*

*If $\phi$ satisfies* $(\Phi1)$–$(\Phi3)$, *then so does $\phi_N$, $N\in\mathbb{N}$, with the same $q_2\ge1$ and functions $g$, $h_0$, $h_1$, $g_0$, $g_1$ and $\omega$. Moreover, $|F_\phi-F_{\phi_N}|_2\to0$ as $N\to\infty$ uniformly on balls in $H_0^1$.*

*If, in addition, $\phi$ satisfies* $(\Phi4)$, *then $\phi_N$ is twice continuously differentiable and*

$$|\partial_{xx}^2\phi_N(r,x)|\le c_\theta g_2(r)+c_\theta g_3(r)\left|\frac{x}{\sqrt{r(1-r)}}\right|^{1/2}\omega\left(\frac{|x|}{\sqrt{r(1-r)}}\right),$$
$$r\in(0,1),x\in\mathbb{R},$$

*and*

$$|\partial_{xr}^2\phi_N(r,x)|\le c_\theta g_4(r)+c_\theta g_5(r)\left|\frac{x}{\sqrt{r(1-r)}}\right|^{3/2}\omega\left(\frac{|x|}{\sqrt{r(1-r)}}\right),$$
$$r\in(0,1),x\in\mathbb{R},$$

*with $c_\theta:=1\vee\sup_\xi\xi^2|\theta''(\xi)|$.*

PROOF. The first assertion is obvious. Then, given that $\phi$ satisfies ($\Phi 1$), ($\Phi 2$), so does $\phi_N$ since $\theta_N$ is an odd contraction. Note that $\theta_N(\eta) - \theta_N(\xi) < 0$ whenever $\eta < \xi$ and $0 \le \theta_N(\eta) - \theta_N(\xi) \le \eta - \xi$ for $\eta \ge \xi$. So, ($\Phi 3$) holds also for $\phi_N$ if it holds for $\phi$. To prove the next assertion, we note that, since $\theta_N(x) = x$ if $|x| \le N$, for $x \in H_0^1$, condition ($\Phi 1$) implies that $\{g(1 + |x|_\infty^{q_2}) \le N\} \subset \{\phi(\cdot, x(\cdot)) = \phi_N(\cdot, x(\cdot))\}$. Hence, again by ($\Phi 1$),

$$|F_\phi - F_{\phi_N}|_2^2(x) = \int |\phi(r, x(r)) - \phi_N(r, x(r))|^2 \, dr$$

$$\le 4(1 + |x'|_2^{q_2})^2 \int \mathbb{1}_{\{g \ge N/(1+|x|_\infty^{q_2})\}} g^2(r) \, dr,$$

which converges to zero as $N \to \infty$ uniformly for $x$ in any ball in $H_0^1$.

Finally, the last assertion follows from the following identities: with $\theta_{(2)}(\xi) := \xi \theta''(\xi)$,

$$\partial_{xx}^2 \phi_N = (\theta_N' \circ \phi) \partial_{xx}^2 \phi + \mathbb{1}_{\{|\phi| \ge N\}} \left( \theta_{(2)} \circ \frac{\phi}{N} \right) \frac{(\partial_x \phi)^2}{\phi},$$

$$\partial_{xr}^2 \phi_N = (\theta_N' \circ \phi) \partial_{xr}^2 \phi + \mathbb{1}_{\{|\phi| \ge N\}} \left( \theta_{(2)} \circ \frac{\phi}{N} \right) \frac{\partial_x \phi \, \partial_r \phi}{\phi}. \qquad \square$$

LEMMA 4.5. Let $\delta \in C_c^\infty((-1,1))$, nonnegative, even, and $\int \delta(x) \, dx = 1$. For $\beta \in (0,1)$, $x \in \mathbb{R}$, $r \in (0,1)$, let

$$\delta_\beta(r, x) := \frac{1}{\beta\sqrt{r(1-r)}} \delta \left( \frac{x}{\beta\sqrt{r(1-r)}} \right)$$

and

$$\phi_\beta(r, x) := \int_{\mathbb{R}} \phi(r, x - y) \delta_\beta(r, y) \, dy.$$

Then $\phi_\beta(r, \cdot) \in C^\infty(\mathbb{R})$ for all $r \in (0,1)$.

If $\phi$ is bounded, then, for $\beta \in (0,1)$, $n = 0, 1, 2, \ldots$, $x \in \mathbb{R}$ and $r \in (0,1)$,

$$\left| \frac{\partial^n}{\partial x^n} \phi_\beta \right| (r, x) \le |\phi|_\infty \frac{\int_{\mathbb{R}} |\delta^{(n)}|(y) \, dy}{(\beta\sqrt{r(1-r)})^n}.$$

If $\phi$ satisfies ($\Phi 1$)–($\Phi 3$), then $\phi_\beta$, $\beta \in (0,1)$, does so, with the same $q_1 \in [2, \infty]$ and $q_2 \in [1, \infty)$ and functions $h_1$ and $g_1$ and $g' = 2^{q_2+1} g$, $h_0' = h_1 + h_0 + 2^{q_2+2} g$, $g_0' = g_0 + 9(\sup_r \omega(r)) g_1$, and $\omega'(r) := \frac{9}{4} \sup\{\omega(s) | s > \frac{r}{2}\}$, $r > 0$.

Moreover, $|F_\phi - F_{\phi_\beta}|_2(x) \to 0$ as $\beta \to 0$ uniformly on balls in $H_0^1$.

PROOF. The first two assertions are well-known properties of the convolution. By ($\Phi 1$), for all $\beta \in (0,1)$, $x \in \mathbb{R}$ and $r \in (0,1)$,

$$|\phi_\beta(r,x)| \le g(r) \int_\mathbb{R} (1+|x-y|^{q_2}) \delta_\beta(r,y)\, dy$$

$$\le 2^{q_2} g(r) \left(1 + |x|^{q_2} + (\beta\sqrt{r(1-r)}\,)^{q_2} \int |y|^{q_2} \delta(y)\, dy\right).$$

So, all $\phi_\beta$, $\beta \in (0,1)$, satisfy (Φ1) with $g' = 2^{q_2+1} g$.

By Remark 2.1(iv), since $\phi_\beta$ satisfy (Φ1) uniformly in $\beta \in (0,1)$, it suffices to verify (Φ2) for all $x \in \mathbb{R}$, $|x| > 1$. Then $\text{sign}(x-y) = \text{sign}(x)$ for all $y \in \bigcup_{\beta, r \in (0,1)} \text{supp}\, \delta_\beta(r,\cdot) \subset (-1,1)$, $\beta \in (0,1)$. Since $\phi$ satisfies (Φ2), for a.e. $r \in (0,1)$, all $x \in \mathbb{R}$, $|x| > 1$, $\beta \in (0,1)$, we obtain

$$\phi_\beta(r,x)\,\text{sign}(x) = \int_\mathbb{R} \text{sign}(x-y) \phi(r, x-y) \delta_\beta(r,y)\, dy$$

$$\le h_0(r) + h_1(r) \int_\mathbb{R} |x-y| \delta_\beta(r,y)\, dy$$

$$\le h_0(r) + h_1(r) \left(|x| + \beta\sqrt{r(1-r)} \int_\mathbb{R} |y| \delta(y)\, dy\right).$$

Hence, $\phi_\beta$, $\beta \in (0,1)$, satisfy (Φ2) with the same $h_1$ as $\phi$ does and with $h_0' = h_1 + h_0 + 2^{q_2+2} g$.

Set $\xi(r,x) := \frac{x}{\sqrt{r(1-r)}}$, $x \in \mathbb{R}$, $r \in (0,1)$. By (Φ3), for all $\rho \in (0, \rho_0)$, $x \in \mathbb{R}$, $N \in \mathbb{N}$, $\beta \in (0,1)$, $r \in (0,1)$,

$$\frac{1}{\rho}(\phi_\beta(r, x+\rho) - \phi_\beta(r,x))$$

$$= \frac{1}{\rho} \int_\mathbb{R} (\phi(r, x+\rho-y) - \phi(r, x-y)) \delta_\beta(r,y)\, dy$$

$$\le g_0(r) + g_1(r) \int_\mathbb{R} |\xi(r, x-y)|^{2-1/p_1} \omega(|\xi(r,x-y)|) \delta_\beta(r,y)\, dy$$

$$= g_0(r) + g_1(r) \int_\mathbb{R} |\xi(r,x) - \beta y|^{2-1/p_1} \omega(|\xi(r,x) - \beta y|) \delta(y)\, dy.$$

By Remark 2.1(iv), we may assume $\omega$ nonincreasing, by replacing $\omega$ with $\tilde\omega(r) := \sup_{s>r} \omega(s)$. Then, for $|\xi(r,x)| \le 2$,

$$\int_\mathbb{R} |\xi(r,x) - \beta y|^{2-1/p_1} \omega(|\xi(r,x) - \beta y|) \delta(y)\, dy \le 9\omega(0),$$

and, for $|\xi(r,x)| > 2$, $\frac{1}{2}|\xi(r,x)| \le |\xi(r,x) - \beta y| \le \frac{3}{2}|\xi(r,x)|$, provided $|y| \le 1$, hence,

$$\int_\mathbb{R} |\xi(r,x) - \beta y|^{2-1/p_1} \omega(|\xi(r,x) - \beta y|) \delta(y)\, dy \le (\tfrac{3}{2}|\xi(r,x)|)^{2-1/p_1} \omega(\tfrac{1}{2}|\xi(r,x)|).$$

Thus, $\phi_\beta$, $\beta \in (0,1)$, satisfy ($\Phi 3$) with the same $g_1$ as $\phi$ does, and with $g'_0 = g_0 + 9\omega(0)g_1$ and $\omega'(r) := \frac{9}{4}\tilde\omega(\frac{r}{2})$, $r \in \mathbb{R}_+$.

Finally, to prove the last assertion, we first note that, for all $x \in H_0^1$ and $\beta \in (0,1)$,

$$|F_\phi - F_{\phi_\beta}|_2^2(x) = \int_0^1 |\phi(r,x(r)) - \phi_\beta(r,x(r))|^2 \, dr$$

$$\leq \int_0^1 \sup_{\substack{y \in \mathbb{R} \\ |y| \leq |x'|_2}} |\phi(r,y) - \phi_\beta(r,y)|^2 \, dr.$$

But $\phi_\beta(r,y) \to \phi(r,y)$ as $\beta \to 0$ locally uniformly in $y$ for all $r \in (0,1)$ and, since we have seen that each $\phi_\beta$ satisfies ($\Phi 1$) with $2^{q_2+1}g$ and $q_2$, we also have that the integrand is bounded by

$$2^{2q_2+4}g(r)^2(1 + (|x'|_2 + 1)^{q_2})^2.$$

Therefore, the last assertion follows by Lebesgue's dominated convergence theorem. $\square$

LEMMA 4.6. *Define, for $N \in \mathbb{N}$, $u \in W_{\mathrm{loc}}^{2,1}(E_N)$,*

$$L_{\phi,\psi}u(x) := \frac{1}{2}\mathrm{Tr}(A_N D^2 u)(x) + (x'' + F_\phi(x) + G_\psi(x), Du(x)), \qquad x \in E_N.$$

*Assume that ($\Phi 2$) holds. Let $\kappa_0 := \frac{2 - |h_1|_1}{8a_0}$. For $\kappa \in (0, \kappa_0)$, let $\lambda_\kappa := 2\kappa \, \mathrm{Tr} \, A + \frac{|h_0|_1^2 \kappa}{4 - 2|h_1|_1 - 8\kappa a_0}$. Then*

$$(4.10) \qquad L_{\phi,\psi}V_\kappa := L_{\phi,\psi}(V_\kappa \restriction_{E_N}) \leq \lambda_\kappa V_\kappa \qquad on \ E_N,$$

*and, for all $\lambda > 2\lambda_\kappa$,*

$$(4.11) \qquad L_{\phi,\psi}V_\kappa \leq \lambda V_\kappa - m_{\kappa,\lambda}\Theta_\kappa \qquad on \ E_N,$$

*with*

$$(4.12) \qquad m_{\kappa,\lambda} := \min\left(\frac{\lambda}{2}, 2\kappa - |h_1|_1\kappa - \frac{|h_0|_1^2\kappa^2}{\lambda - 4\kappa\,\mathrm{Tr}\,A} - 4a_0\kappa^2\right) \ (>0).$$

*Moreover, for all $q \in [2,\infty)$ and $\kappa \in (0,\kappa_0)$, there exist $\lambda_{q,\kappa} > 2\lambda_\kappa$ and $m_{q,\kappa} < \min\{q(q-1), m_{\kappa,\lambda}\}$ depending only on $q$, $\kappa$, $|h_0|_1$, $|h_1|_1$, $|A|_{X \to X}$ and $\mathrm{Tr} \, A$ such that*

$$(4.13) \quad L_{\phi,\psi}V_{q,\kappa} := L_{\phi,\psi}(V_{q,\kappa}\restriction_{E_N}) \leq \lambda_{q,\kappa}V_{q,\kappa} - m_{q,\kappa}\Theta_{q,\kappa} \qquad on \ E_N.$$

PROOF. First observe that, due to ($\Phi 2$), for all $q \in [2,\infty)$ and $x \in H_0^1$,

$$(4.14) \quad (F_\phi(x), x|x|^{q-2}) \leq \int_0^1 (h_1|x|^q + h_0|x|^{q-1}) \, dr \leq |h_1|_1|x|_\infty^q + |h_0|_1|x|_\infty^{q-1}$$

and

$$(4.15) \quad \begin{aligned} (G_\psi(x), x|x|^{q-2}) &= -(q-1)\int_0^1 x'|x|^{q-2}\psi \circ x\, dr \\ &= -(q-1)\int_{x(0)}^{x(1)} \psi(\tau)|\tau|^{q-2}\, d\tau = 0, \end{aligned}$$

since $x(1) = x(0) = 0$.

To prove the first assertion, note that, for $x \in E_N$, $i, j = 1, \dots, N$,

$$(4.16) \qquad \partial_i |x|_2^2 = 2(x, \eta_i) \quad \text{and} \quad \partial_{ij}^2 |x|_2^2 = 2(\eta_i, \eta_j) = 2\delta_{ij}.$$

So, we have, for $x \in E_N$ by (4.15) with $q = 2$,

$$(4.17) \quad L_{\phi,\psi} V_\kappa(x) = 2\kappa e^{\kappa|x|_2^2}(\mathrm{Tr}\, A_N + (F_\phi(x), x) + 2\kappa(x, Ax) - |x'|_2^2).$$

Now (4.14) for $q = 2$, together with the estimates $|x|_\infty \leq \frac{1}{\sqrt{2}}|x'|_2$ and the inequality $ab \leq 2\varepsilon a^2 + \frac{b^2}{8\varepsilon}$, $a, b, \varepsilon > 0$, imply that, for all $\varepsilon > 0$ and $x \in H_0^1$,

$$(F_\phi(x), x) \leq (\tfrac{1}{2}|h_1|_1 + \varepsilon)|x'|_2^2 + \frac{|h_0|_1^2}{8\varepsilon},$$

hence,

$$\mathrm{Tr}\, A_N + (F_\phi(x), x) + 2\kappa(x, Ax) - |x'|_2^2$$

$$\leq \mathrm{Tr}\, A + \frac{|h_0|_1^2}{8\varepsilon} - \left(1 - \frac{1}{2}|h_1|_1 - \varepsilon - 2\kappa a_0\right)|x'|_2^2.$$

So, (4.10) follows by choosing $\varepsilon > 0$ so that the last term in brackets is equal to zero. Equation (4.11) follows by choosing $\varepsilon > 0$ so that

$$2\kappa\left(\mathrm{Tr}\, A + \frac{|h_0|_1^2}{8\varepsilon}\right) = \frac{\lambda}{2}.$$

To prove the second assertion, observe that, for $x \in E_N$, $i, j = 1, \dots, N$,

$$(4.18) \quad \begin{aligned} \partial_i |x|_q^q &= q(x|x|^{q-2}, \eta_i), \\ \partial_{ij}^2 |x|_q^q &= q(q-1)(|x|^{q-2}\eta_i, \eta_j), \\ \partial_j(x|x|^{q-2}, \eta_i) &= (q-1)(|x|^{q-2}, \eta_i\eta_j), \\ (x|x|^{q-2}, x'') &= -(q-1)|x'|x|^{q/2-1}|_2^2. \end{aligned}$$

So by (4.15), we have, for $x \in E_N$,

$$(4.19) \quad \begin{aligned} L_{\phi,\psi} V_{q,\kappa}(x) &= (1 + |x|_q^q) L_{\phi,\psi} V_\kappa(x) \\ &\quad + q e^{\kappa|x|_2^2}[(F_\phi(x), |x|^{q-2}x) + 4\kappa(Ax, |x|^{q-2}x)] \\ &\quad + q(q-1)e^{\kappa|x|_2^2}\left[\left(|x|^{q-2}, \sum_{i=1}^N A_{ii}\eta_i^2\right) - |x'|x|^{q/2-1}|_2^2\right]. \end{aligned}$$

It follows from (4.11) that, for all $\lambda > 2\lambda_\kappa$, $x \in E_N$,

$$(4.20) \qquad (1 + |x|_q^q)L_{\phi,\psi}V_\kappa(x) \le V_{q,\kappa}(x)(\lambda - m_{\kappa,\lambda}(|x'|_2^2 + 1)).$$

Below we shall use the following consequence of the inequality $|z|_\infty^2 \le 2|z'|_2|z|_2$, $z \in H_0^1$: For $x \in H_0^1$ and $q \ge 2$,

$$(4.21) \qquad \begin{aligned} |x|_\infty^q &= |x|x|^{q/2-1}|_\infty^2 \le 2|(x|x|^{q/2-1})'|_2|x|x|^{q/2-1}|_2 \\ &= q|x'|x|^{q/2-1}|_2|x|_q^{q/2}. \end{aligned}$$

It follows by (4.14) and (4.21), together with Young's inequality, that there exists $c_1(q) > 0$ depending only on $q$, such that, for all $\varepsilon > 0$,

$$(4.22) \qquad \begin{aligned} (F_\phi(x), |x|^{q-2}x) \\ &\le |h_1|_1|x|_\infty^q + |h_0|_1|x|_\infty^{q-1} \\ &\le q|h_1|_1|x'|x|^{q/2-1}|_2|x|_q^{q/2} \\ &\quad + q^{(q-1)/q}|h_0|_1|x'|x|^{q/2-1}|_2^{(q-1)/q}|x|_q^{(q-1)/2} \\ &\le \varepsilon|x'|x|^{q/2-1}|_2^2 \\ &\quad + c_1(q)(|h_1|_1^2\varepsilon^{-1} + |h_0|_1^{2q/(q+1)}\varepsilon^{-(q-1)/(q+1)})(1 + |x|_q^q). \end{aligned}$$

It follows from the estimate $|z|_p \le |z|_\infty$, (4.21) and Young's inequality that, for every $\varepsilon > 0$,

$$(4.23) \qquad \begin{aligned} |(Ax, |x|^{q-2}x)| &\le |A|_{X \to X}|x|_2|x|_{2q-2}^{q-1} \le |A|_{X \to X}|x|_\infty^q \\ &\le q|A|_{X \to X}|x'|x|^{q/2-1}|_2|x|_q^{q/2} \\ &\le \varepsilon|x'|x|^{q/2-1}|_2^2 + \frac{q^2}{4\varepsilon}|A|_{X \to X}^2|x|_q^q. \end{aligned}$$

Next, observe that $\sum_{i=1}^N A_{ii}\eta_i^2(r) \ge 0$ for all $r \in (0,1)$. Hence, it follows by (4.21) and Young's inequality that there exists $c_2(q) > 0$ depending only on $q$, such that, for every $\varepsilon > 0$,

$$(4.24) \qquad \begin{aligned} \left(|x|^{q-2}, \sum_{i=1}^N A_{ii}\eta_i^2\right) &\le |x|_\infty^{q-2}\sum_{i=1}^N A_{ii} \\ &\le q^{(q-2)/q}|x'|x|^{q/2-1}|_2^{(q-2)/q}|x|_q^{(q-2)/2}\operatorname{Tr} A \\ &\le \varepsilon|x'|x|^{q/2-1}|_2^2 \\ &\quad + c_2(q)(\operatorname{Tr} A)^{2q/(q+2)}\varepsilon^{-(q-2)/(q+2)}(1 + |x|_q^q). \end{aligned}$$

Collecting (4.22), (4.23) and (4.24), we conclude that there exists $c_q > 0$ depending only on $q$, such that, for every $\varepsilon \in (0,1)$,

$$qe^{\kappa|x|_2^2}[(F_\phi(x), |x|^{q-2}x) + 4\kappa(Ax, |x|^{q-2}x)]$$

$$+ q(q-1)e^{\kappa|x|_2^2}\left[\left(|x|^{q-2}, \sum_{i=1}^{N} A_{ii}\eta_i^2\right) - |x'|\,|x|^{q/2-1}|_2^2\right]$$

$$\leq c_q\varepsilon^{-1}(|h_1|_1^2 + |h_0|_1^{2q/(q+1)} + \kappa|A|_{X\to X}^2 + (\operatorname{Tr} A)^{2q/(q+2)})V_{q,\kappa}(x)$$

$$- q(q-1-(4\kappa+q)\varepsilon)V_\kappa(x)|x'|\,|x|^{q/2-1}|_2^2.$$

This together with (4.20) and (4.19) implies (4.13).   $\square$

LEMMA 4.7.   *Let $\phi$ be continuously differentiable in the second variable such that $\sup_{|\xi|\leq R}\phi_x(\cdot,\xi) \in L^1(0,1)$ for all $R > 0$ and let $\phi$ satisfy* ($\Phi$3).

*Then there exists a nonnegative function $\varepsilon \mapsto C(\varepsilon)$ depending only on $\omega$, $p_1$, $|g_0|_1$ and $|g_1|_{p_1}$ such that, for all $\varepsilon > 0$ and $x, y \in H_0^1$,*

$$\partial_y(F_\phi, y)(x) \leq \tfrac{1}{2}|y'|_2^2 + (\varepsilon|x'|_2^2 + C(\varepsilon))|y|_2^2.$$

*If, moreover, $\phi$ is twice continuously differentiable and there exist $g_2, g_3 \in L_+^2(0,1)$, $g_4, g_5 \in L_+^1(0,1)$ and a bounded Borel-measurable function $\omega:\mathbb{R}_+ \to \mathbb{R}_+$, $\omega(r) \to 0$ as $r \to \infty$, such that*

(4.25)
$$|\phi_{xx}(r,x)| \leq g_2(r) + g_3(r)\left|\frac{x}{\sqrt{r(1-r)}}\right|^{1/2}\omega\left(\frac{|x|}{\sqrt{r(1-r)}}\right),$$
$$r \in (0,1), x \in \mathbb{R},$$

*and*

$$|\phi_{xr}(r,x)| \leq g_4(r) + g_5(r)\left|\frac{x}{\sqrt{r(1-r)}}\right|^{3/2}\omega\left(\frac{|x|}{\sqrt{r(1-r)}}\right), \qquad r \in (0,1), x \in \mathbb{R}$$

*[which is the case, if $\phi$ satisfies* ($\Phi$4)] *then there exists a nonnegative function $\varepsilon \mapsto C(\varepsilon)$ depending only on $\omega$, $p_1$, $|g_0|_1$, $|g_1|_{p_1}$ $|g_2|_2$, $|g_3|_2$, $|g_4|_1$ and $|g_5|_1$ such that, for all $\varepsilon > 0$ and $x, y \in H_0^1$,*

$$\partial_{(-\Delta)^{1/2}y}(F_\phi, (-\Delta)^{-1/2}y)(x) \leq \tfrac{1}{2}|y'|_2^2 + (\varepsilon|x'|_2^2 + C_\varepsilon)|y|_2^2.$$

PROOF.   As before, we set $\sigma(r,x) := \frac{|x|}{\sqrt{r(1-r)}}$. Since $\phi$ is continuously differentiable in the second variable, ($\Phi$3) implies that, for all $x \in \mathbb{R}$ and $r \in (0,1)$,

(4.26)        $\phi_x(r,x) \leq g_0(r) + g_1(r)|\sigma(r,x)|^{2-1/p_1}\omega(\sigma(r,x)).$

Fix $x \in H_0^1$. Note that, for $\xi, \eta \in H_0^1$, since $\sup_{|\xi| \le R} \phi_x(\cdot, \xi) \in L^1(0, 1)$ for all $R > 0$, we have

$$\partial_\xi(F_\phi, \eta)(x) = \int_0^1 \xi(r) \eta(r) \phi_x(r, x(r)) \, dr.$$

Hence, (4.26) implies that, for $y \in H_0^1$,

$$\begin{aligned}
\partial_y(F_\phi, y)(x) &= \int_0^1 (y^2)(r) \phi_x(r, x(r)) \, dr \\
&\le |y|_\infty^2 |g_0|_1 + |y|_{2p_1/(p_1-1)}^2 |g_1|_{p_1} |\sigma^{2-1/p_1} \omega \circ \sigma|_\infty(x),
\end{aligned}$$

where, for $\alpha, \beta \ge 0$, we set

$$|\sigma^\alpha \omega^\beta \circ \sigma|_\infty(x) := \sup_{r \in (0,1)} |\sigma^\alpha(r, x(r)) \omega^\beta(\sigma(r, x(r)))|.$$

Note that, for $y \in H_0^1$, $|y|_\infty^2 \le 2|y'|_2 |y|_2$ and, hence,

$$|y|_{2p_1/(p_1-1)}^2 \le |y|_\infty^{2/p_1} |y|_2^{2(p_1-1)/p_1} \le 2|y'|_2^{1/p_1} |y|_2^{(2p_1-1)/p_1}.$$

Hence, by Young's inequality, there exists $\hat{c}_{p_1} > 0$ such that

$$\partial_y(F_\phi, y)(x) \le \tfrac{1}{2}|y'|_2^2 + \hat{c}_{p_1}|y|_2^2[|g_0|_1^2 + |g_1|_{p_1}^{2p_1/(2p_1-1)} |\sigma^2 \omega^{2p_1/(2p_1-1)} \circ \sigma|_\infty(x)].$$

Observe now that, for all $\varepsilon > 0$,

$$\hat{c}_{p_1}|g_1|_{p_1}^{2p_1/(2p_1-1)} |\sigma^2 \omega^{2p_1/(2p_1-1)} \circ \sigma|_\infty(x) \le \varepsilon|\sigma|_\infty^2(x) + \hat{C}(\varepsilon),$$

with $\hat{C}(\varepsilon) := \sup\{\hat{c}_{p_1}|g_1|_{p_1}^{2p_1/(2p_1-1)} s^2 \omega^{2p_1/(2p_1-1)}(s)|s \ge 0$ such that $\hat{c}_{p_1}|g_1|_{p_1}^{2p_1/(2p_1-1)} \omega^{2p_1/(2p_1-1)}(s) > \varepsilon\}$. Now the first assertion follows from the inequality $|\sigma|_\infty(x) = \sup_r \frac{|x|(r)}{\sqrt{r(1-r)}} \le \sqrt{2}|x'|_2$, $x \in H_0^1$, which is a consequence of the fundamental theorem of calculus (or of Sobolev embedding).

To prove the second assertion, let $z := (-\Delta)^{-1/2} y$, $y \in H_0^1$. Then $(-\Delta)^{1/2} y = -z''$, $|z'|_2 = |y|_2$ and $|z''|_2 = |y'|_2$. Moreover,

$$\begin{aligned}
\partial_{(-\Delta)^{1/2} y}(F_\phi, (-\Delta)^{-1/2} y)(x) &= -\int_0^1 z''(r) z(r) \phi_x(r, x(r)) \, dr \\
&= \int_0^1 |z'|^2(r) \phi_x(r, x(r)) \, dr \\
(4.27) \qquad\qquad & \quad + \int_0^1 z'(r) z(r) x'(r) \phi_{xx}(r, x(r)) \, dr \\
& \quad + \int_0^1 z'(r) z(r) \phi_{xr}(r, x(r)) \, dr.
\end{aligned}$$

We can estimate the first term in the right-hand side of (4.27) in the same way as above. Indeed, note that (4.25) was shown in the proof of Remark 2.1(v) to imply (4.26). So, as above, we obtain that there exists a nonnegative function $\varepsilon \mapsto C_1(\varepsilon)$ depending only on $\omega$, $p_1$, $|g_0|_1$ and $|g_1|_{p_1}$ such that, for all $\varepsilon > 0$,

$$
\begin{aligned}
(4.28) \quad \int_0^1 |z'|^2(r)\phi_x(r, x(r))\, dr &\le \tfrac{1}{4}|z''|_2^2 + (\varepsilon|x'|_2^2 + C_1(\varepsilon))|z'|_2^2 \\
&\le \tfrac{1}{4}|y'|_2^2 + (\varepsilon|x'|_2^2 + C_1(\varepsilon))|y|_2^2.
\end{aligned}
$$

To estimate the second and the last terms in the right-hand side of (4.27), we note that

$$
|z'|_\infty \le (2|z''|_2|z'|_2)^{1/2} = (2|y'|_2|y|_2)^{1/2}, \qquad |z|_\infty \le 2^{-1/2}|z'|_2 = 2^{-1/2}|y|_2.
$$

By (4.25) and the estimate $|\sigma|_\infty(x) \le \sqrt{2}|x'|_2$, we conclude that, for all $\varepsilon > 0$,

$$
\begin{aligned}
|\phi_{xx}(\cdot, x)|_2 &\le |g_2|_2 + |g_3|_2|\sigma^{1/2}\omega \circ \sigma|_\infty(x) \\
&\le \frac{1}{6 \cdot 2^{1/4}}\varepsilon|\sigma|_\infty^{1/2}(x) + C_2(\varepsilon) \le \frac{1}{6}\varepsilon|x'|_2^{1/2} + C_2(\varepsilon)
\end{aligned}
$$

and

$$
\begin{aligned}
|\phi_{xr}(\cdot, x)|_1 &\le |g_4|_1 + |g_5|_1|\sigma^{3/2}\omega \circ \sigma|_\infty(x) \\
&\le \frac{1}{6 \cdot 2^{3/4}}\varepsilon|\sigma|_\infty^{3/2}(x) + C_3(\varepsilon) \le \frac{1}{6}\varepsilon|x'|_2^{3/2} + C_3(\varepsilon),
\end{aligned}
$$

with

$$
C_2(\varepsilon) := |g_2|_2 + \sup\left\{|g_3|_2 s^{1/2}\omega(s)\,|\, s \ge 0 \text{ such that } |g_3|_2\omega(s) > \frac{1}{6 \cdot 2^{1/4}}\varepsilon\right\},
$$

$$
C_3(\varepsilon) := |g_4|_1 + \sup\left\{|g_5|_1 s^{3/2}\omega(s)\,|\, s \ge 0 \text{ such that } |g_5|_1\omega(s) > \frac{1}{6 \cdot 2^{3/4}}\varepsilon\right\}.
$$

Thus, it follows from Young's inequality that there exists a nonnegative function $\varepsilon \mapsto \tilde{C}(\varepsilon)$ dependent on $\omega$, $|g_2|_2$, $|g_3|_2$, $|g_4|_1$ and $|g_5|_1$ only such that, for all $\varepsilon \in (0, 1)$,

$$
\begin{aligned}
\int_0^1 &z'(r)z(r)x'(r)\phi_{xx}(r, x(r))\, dr + \int_0^1 z'(r)z(r)\phi_{xr}(r, x(r))\, dr \\
&\le |y'|_2^{1/2}|y|_2^{3/2}[\tfrac{1}{6}\varepsilon|x'|_2^{3/2} + |x'|_2 C_2(\varepsilon) + \tfrac{1}{6}\varepsilon|x'|_2^{3/2} + C_3(\varepsilon)] \\
&\le \tfrac{1}{4}|y'|_2^2 + (\varepsilon|x'|_2^2 + \tilde{C}(\varepsilon))|y|_2^2.
\end{aligned}
$$

Now the second assertion follows from (4.28). $\quad\square$

LEMMA 4.8. *Let $\psi$ satisfy* ($\Psi$).

*Then there exists a nonnegative function $\varepsilon \mapsto C(\varepsilon)$ depending on $\omega$ and $C$ such that, for all $\varepsilon > 0$ and $x, y \in H_0^1$,*

$$(4.29) \qquad \begin{aligned} \partial_y(G_\psi, y)(x) &\leq \tfrac{1}{2}|y'|_2^2 + (\varepsilon|x'|_2^2 + C_\varepsilon)|y|_2^2, \\ \partial_{(-\Delta)^{1/2}y}(G_\psi, (-\Delta)^{-1/2}y)(x) &\leq \tfrac{1}{2}|y'|_2^2 + (\varepsilon|x'|_2^2 + C_\varepsilon)|y|_2^2. \end{aligned}$$

PROOF. Fix $x \in H_0^1$. Note that, for $\xi, \eta \in H_0^1$, we have

$$\partial_\xi(G_\psi, \eta)(x) = -\int_0^1 \xi\eta'\psi_x \circ x \, dr.$$

Hence, for all $y \in H_0^1$,

$$\begin{aligned} \partial_y(G_\psi, y)(x) &= -\tfrac{1}{2}\int_0^1 (y^2)'\psi_x \circ x \, dr \\ &= \tfrac{1}{2}\int_0^1 y^2 x'\psi_{xx} \circ x \, dr \leq \tfrac{1}{2}|y|_4^2|x'|_2|\psi_{xx} \circ x|_\infty. \end{aligned}$$

Set $z := (-\Delta)^{-1/2}y$ so that $(-\Delta)^{1/2}y = -z''$. Then

$$\partial_{(-\Delta)^{1/2}y}(G_\psi, (-\Delta)^{-1/2}y)(x) = \int_0^1 z'z''\psi_x \circ x \, dr \leq \tfrac{1}{2}|z'|_4^2|x'|_2|\psi_{xx} \circ x|_\infty.$$

Note that $|y|_4^2 \leq |y|_\infty|y|_2 \leq \sqrt{2}|y'|_2^{1/2}|y|_2^{3/2}$. Hence, by Young's inequality, there exists $\hat{c} > 0$ such that

$$\partial_y(G_\psi, y)(x) \leq \tfrac{1}{2}|y'|_2^2 + \hat{c}|y|_2^2|x'|_2^{4/3}|\psi_{xx} \circ x|_\infty^{4/3},$$

$$\partial_{(-\Delta)^{1/2}y}(G_\psi, (-\Delta)^{-1/2}y)(x) \leq \tfrac{1}{2}|y'|_2^2 + \hat{c}|y|_2^2|x'|_2^{4/3}|\psi_{xx} \circ x|_\infty^{4/3}.$$

($\Psi$) implies that, for all $\varepsilon > 0$, $|\psi_{xx}|^{4/3}(x) \leq \varepsilon|x|^{2/3} + \hat{C}(\varepsilon)$ with $\hat{C}(\varepsilon) := \sup\{(C + r^{1/2}\omega(r))^{4/3}|r \geq 0 \text{ such that } Cr^{-1/2} + \omega(r) > \varepsilon^{3/4}\}$. Now the assertion follows from the estimate $|x|_\infty \leq \tfrac{1}{\sqrt{2}}|x'|_2$. □

LEMMA 4.9. *Assume that $\sup_{|x| \leq R}|\phi(\cdot, x)| \in L^2(0,1)$ for all $R > 0$.*

*Then $F_\phi : H_0^1 \to L^2(0,1)$ is $|\cdot|_2$-continuous on $|\cdot|_{H_0^1}$-balls.*

*If, in addition, $\sup_x |\phi(\cdot, x)|_1 < \infty$ and $\phi$ is differentiable in the second variable with $\sup_{|\xi| \leq R}|\phi_x(\cdot, \xi)| \in L^1(0,1)$ for all $R > 0$, then, for all $N \in \mathbb{N}$, $P_N F_\phi \circ P_N : E_N \to E_N$ is bounded and locally Lipschitz continuous.*

*If $\phi$ satisfies* ($\Phi 1$), *then, for all $p \in [2, \infty)$, there exists $c_{p,q_1,q_2} > 0$ such that*

$$|F_\phi|_2(x) \leq c_{p,q_1,q_2}|g|_{q_1}\Theta_{p,\kappa}^{(q_2-1+2/q_1)/(p+2)}(x) \qquad \text{for all } x \in H_0^1, \kappa > 0.$$

PROOF.    Let $(x_n)_{n\in\mathbb{N}}$ be a bounded sequence in $H_0^1$ and $\lim_{n\to\infty} x_n = x$ in $H_0^1$ in the $|\cdot|_{1,2}$-topology. Since a $|\cdot|_{1,2}$-bounded set is compact in $C_0(0,1)$, we conclude that $x_n \to x$ uniformly on $(0,1)$ and, hence, $\phi(r, x_n(r)) \to \phi(r, x(r))$ as $n \to \infty$ for all $r \in (0,1)$ and $\sup_n |\phi|(r, x_n(r)) \le \sup_{|\xi| \le |x|_\infty + 1} |\phi|(r, \xi) \in L^2(0,1)$. Thus, the first assertion follows by the dominated convergence theorem.

Let now the second assumption hold. Then, for all $n \in \mathbb{N}$, $x, y \in H_0^1$,

$$|(F_\phi(x), \eta_n)| \le \sup_\xi |\phi(\cdot, \xi)|_1 |\eta_n|_\infty$$

$$|(F_\phi(x) - F_\phi(y), \eta_n)| \le \Big| \sup_{|\xi| \le |x|_\infty \vee |y|_\infty} |\phi_x(\cdot, \xi)| \Big|_1 |\eta_n|_\infty |x - y|_\infty.$$

Hence, the second assertion follows.

To prove the last assertion, we first note that by ($\Phi 1$), for all $x \in H_0^1$,

$$|F_\phi|_2(x) \le |g|_{q_1} |1 + |x|^{q_2}|_{2q_1/(q_1-2)} \le |g|_{q_1} (1 + |x|_s^{q_2})$$

with $s := \frac{2q_1 q_2}{q_1 - 2}$, and for $p \in [2, \infty)$,

$$(4.30) \qquad |x|_\infty^{1+p/2} \le \frac{p+2}{2} \int_0^1 |x'| |x|^{p/2}\, dr \le \frac{p+2}{2} |x'|_2 |x|_p^{p/2}.$$

Since $|x|_s^s \le |x|_2^2 |x|_\infty^{s-2}$, it follows that

$$|x|_s^{q_2} \le |x|_2^{2q_2/s} |x|_\infty^{q_2(1-2/s)} \le |x|_2^{2q_2/s} \left[ \left( \frac{p+2}{2} \right)^2 |x'|_2^2 |x|_p^p \right]^{(q_2 - 2q_2/s)/(p+2)}.$$

Substituting $s$, we find

$$|F_\phi|_2(x) \le 2|g|_{q_1} \left( \frac{p+2}{2} \right)^{(2q_2 - 2 + 4/q_1)/(p+2)} (1 + |x|_2^{1-2/q_1})$$

$$\times [(1 + |x'|_2^2)(1 + |x|_p^p)]^{(q_2 - 1 + 2/q_1)/(p+2)},$$

which implies the assertion.    $\square$

LEMMA 4.10.    $G_\psi : H_0^1 \to L^2(0,1)$ is continuous.

If, in addition, $\psi$ is bounded, then, for all $N \in \mathbb{N}$, $P_N G_\psi \circ P_N : E_N \to E_N$ is bounded and locally Lipschitz continuous.

If $|\psi'|(x) \le C(1 + |x|^{q_0})$, then, for all $x \in H_0^1$, $\kappa \in (0, \infty)$, $p \in [2, \infty)$,

$$|G_\psi|_2(x) \le 2C \left( \frac{p+2}{2} \right)^{2q_0/(p+2)} \Theta_{p,\kappa}^{1/2 + q_0/(p+2)}(x).$$

In particular, if $\psi$ satisfies ($\Psi$), then

$$|G_\psi|_2(x) \le 2C \left( \frac{p+2}{2} \right)^{3/(p+2)} \Theta_{p,\kappa}^{(1/2)(1+3/(p+2))}(x) \qquad \text{for all } x \in H_0^1.$$

PROOF. Let $(x_n)_{n \in \mathbb{N}}$ be a $|\cdot|_{1,2}$-bounded sequence such that $\lim_{n \to \infty} x_n = x \in H_0^1$ in the $|\cdot|_2$-topology. Since an $|\cdot|_{1,2}$-bounded set is compact in $C_0(0,1)$, we conclude that $x_n \to x$ uniformly on $(0,1)$ and, hence, $\psi' \circ x_n \to \psi' \circ x$ uniformly on $(0,1)$. Thus, the first assertion follows by the definition of $G_\psi$.

Let now $\psi$ be bounded. Then, for all $n \in \mathbb{N}$, $x, y \in H_0^1$,

$$|(G_\psi(x), \eta_n)| \leq |\psi|_\infty |\eta_n'|_1$$

$$|(G_\psi(x) - G_\psi(y), \eta_n)| \leq \operatorname*{ess\,sup}_{|s| \leq |x|_\infty \vee |y|_\infty} |\psi'(s)| |\eta_n'|_2 |x - y|_2.$$

Hence, the second assertion follows.

The third assertion follows from the estimate $|G_\psi|_2(x) \leq C(1 + |x|_\infty^{q_0})|x'|_2$ and (4.30). The last assertion is then clear, because we can take $q_0 = \frac{3}{2}$ by Remark 2.1(i). $\square$

Now we are prepared for the following:

PROOF OF PROPOSITION 4.1. Let $N \in \mathbb{N}$ and let $B_N$ denote the closed ball in $H_0^1$ of radius $N$. By Lemmas 4.4 and 4.5, there exist $\beta_N \in (0,1)$ such that

$$\sup_{x \in B_N} |F_{\Phi_N} - F_{(\Phi_N)_\beta}|_2(x) \leq \frac{1}{N}$$

for all $\beta \leq \beta_N$ and $\beta_{N+1} \leq \beta_N$. Define

(4.31) $$F_N := F_{(\Phi_N)_{\beta_N}} + G_{\Psi_N}, \qquad N \in \mathbb{N}.$$

Then $\lim_{N \to \infty} |F - F_N|_2 = 0$ uniformly on balls in $H_0^1$, where $F$ is as in (2.15), by Lemmas 4.3 and 4.4. Since by Lemmas 4.9 and 4.10, $F$ is $|\cdot|_2$-continuous on $|\cdot|_{1,2}$-balls, and since $P_N x \to x$ in $H_0^1$ as $N \to \infty$ for all $x \in H_0^1$, it follows that (F2c) holds.

By Lemmas 4.4 and 4.5, it follows that Lemma 4.6 applies to $(\Phi_N)_{\beta_N}$ and $\Psi_N$ for all $q \in [2, \infty)$ with $\kappa_0, \lambda_{q,\kappa}$ and $m_{q,\kappa}$ independent of $N$. So, (F2a) holds.

By Lemmas 4.3 and 4.5, we see that Lemmas 4.7 and 4.8 apply to $(\Phi_N)_{\beta_N}$ and $\Psi_N$ with the functions $\varepsilon \to C_\varepsilon$ independent of $N$. So, (F2b) holds.

Since in Lemma 4.9 we have $(q_2 - 1 + \frac{2}{q_1})/(p+2) \leq 1$ if and only if $p \geq q_2 - 3 + \frac{2}{q_1}$, (F2d) follows by Lemmas 4.9 and 4.10.

The boundedness and local Lipschitz continuity of $F_N$ follow by Lemmas 4.5, 4.9 and 4.10. So, (F2) is proved.

If, in addition, ($\Phi 4$) holds, then (F2e) follows from Lemmas 4.7 and 4.8 in the same way as we have derived (F2b). $\square$

**5. Some properties of the function spaces $WC_{p,\kappa}$, $W_1C_{p,\kappa}$ and $\mathrm{Lip}_{l,p,\kappa}$.**
Below for a topological vector space $\mathcal{V}$ over $\mathbb{R}$ let $\mathcal{V}'$ denote its dual space.

The following we formulate for general completely regular topological spaces and recall that our $X = L^2(0,1)$ equipped with the weak topology is such a space.

Let $X$ be a completely regular topological space, $V : X \to [1, \infty]$ a function, and $X_V := \{V < \infty\}$ equipped with the topology induced by $X$. Analogously to (2.2), we define

$$
(5.1) \quad
\begin{aligned}
C_V := \Big\{ & f : X_V \to \mathbb{R} | f \upharpoonright_{\{V \le R\}} \text{ is continuous } \forall R \in \mathbb{R}_+ \text{ and} \\
& \lim_{R \to \infty} \sup_{\{V \ge R\}} V^{-1}|f| = 0 \Big\},
\end{aligned}
$$

equipped with the norm $\|f\|_V := \sup V^{-1}|f|$. Obviously, $C_V$ is a Banach space.

THEOREM 5.1. *Let $X$ be a completely regular topological space. Let $V : X \to [1, \infty]$ be of metrizable compact level sets $\{V \le R\}$, $R \ge 0$, and let $C_V$ be as above. Then $\sigma(C_V) = \mathcal{B}(X_V)$ and*

$$
(5.2) \quad C_V' = \Big\{ \nu | \nu \text{ is a signed Borel measure on } X_V, \int V \, d|\nu| < \infty \Big\},
$$

$\|\nu\|_{C_V'} = \int V \, d|\nu|$. *In particular, $f_n \to f$ weakly in $C_V$ as $n \to \infty$ if and only if $(f_n)$ is bounded in $C_V$, $f \in C_V$, and $f_n \to f$ pointwise on $X_V$ as $n \to \infty$.*

PROOF. Let $\nu$ be a signed Borel measure on $X_V$ such that $\int V d|\nu| < \infty$. Then $f \mapsto \nu(f) := \int f d\nu$ is a linear functional on $C_V$ and, since

$$
\left| \int f \, d\nu \right| \le \int \frac{|f|}{V} V \, d|\nu| \le \|f\|_V \int V \, d|\nu|,
$$

we conclude that $\nu \in C_V'$ and $\|\nu\|_{C_V'} \le \int V \, d|\nu|$.

Now let $l \in C_V'$. Note that, for every $f \in C_V$, there exists $x \in X_V$ such that $\|f\|_V = |f|(x)V^{-1}(x)$. Hence, we can apply [14], Corollary 36.5, to conclude that there exist positive $l_1, l_2 \in C_V'$ such that $l = l_1 - l_2$ and $\|l\|_{C_V'} = \|l_1\|_{C_V'} + \|l_2\|_{C_V'}$. So, we may assume that $l \ge 0$. Let $f_n \in C_V$, $n \in \mathbb{N}$, such that $f_n \downarrow 0$ as $n \to \infty$. Then by Dini's theorem, $f_n \to 0$ as $n \to \infty$ uniformly on all sets $\{V \le R\}$, $R \ge 1$. Hence, $\|f_n\|_V \to 0$ as $n \to \infty$ so $l(f_n) \to 0$. $C_V$ is a Stone-lattice generating the Borel $\sigma$-algebra on $X_V$. Indeed, we first note that $X_V \in \mathcal{B}(X)$ as a $\sigma$-compact set, and if $B \in \mathcal{B}(X_V)$, then $B = \bigcup_{n=1}^{\infty} B_n$ with $B_n \in \mathcal{B}(K_n)$, $K_n := \{V \le n\}$. But since $K_n$ is a metric space, $\mathcal{B}(K_n) = \sigma(C(K_n))$. But $C(K_n) = C_{V \upharpoonright_{K_n}}$ by Tietze's extension theorem (which holds for compact sets in completely regular spaces). Hence, $\mathcal{B}(K_n) = \sigma(C_{V\restriction K_n}) = \sigma(C_V) \cap K_n$. So, $B \in \sigma(C_V)$. We conclude by the Daniell–Stone theorem (cf., e.g., [5], 39.4) that there exists a positive Borel measure $\nu$ on $X_V$ such that

$$\int f \, d\nu = l(f) \qquad \forall f \in C_V.$$

Since $1 \in C_V$, $\nu$ is a finite measure. To calculate $\|l\|_{C_V'}$, let $f_n \uparrow V$ be a sequence of bounded positive continuous functions on $X_V$ increasing to $V$. Such a sequence exists by [51], Lemma II.1.10, since $X_V$ as a union of metrizable compacts is strongly Lindelöf. Then $f_n \in C_V$ and $\|f_n\|_V \leq 1$ for all $n \in \mathbb{N}$ and

$$\|l\|_{C_V'} \geq \int f_n \, d\nu \to \int V \, d\nu \qquad \text{as } n \to \infty.$$

Hence, $\|l\|_{C_V'} = \int V \, d\nu$. The rest of the assertion follows from the dominated convergence theorem. $\quad\square$

COROLLARY 5.2. *Let $X, Y$ be completely regular topological spaces. Let $\Theta \colon Y \to [1, \infty]$ have metrizable compact level sets, and let $X \colon V \to [1, \infty]$ be a function. Let $X_V$ and $Y_\Theta$, $C_V$ and $C_\Theta$ be as above. Let $M \colon C_\Theta \to C_V$ be a positive bounded linear operator. Then there exists a kernel $m(x, dy)$ from $X_V$ to $Y_\Theta$ such that, for all $f \in C_\Theta$, $Mf(x) = \int f(y) m(x, dy)$ and $\int \Theta(y) m(x, dy) \leq \|M\|_{C_\Theta \to C_V} V(x)$.*

COROLLARY 5.3. *An algebra of bounded continuous functions on $X_V$ generating $\mathcal{B}(X_V)$ is dense in $C_V$.*

PROOF. By a simple monotone class argument, it follows that the algebra forms a measure determining class on $X_V$. So by Theorem 5.1, it follows that the algebra is dense in $C_V$ with respect to the weak topology, hence, also with respect to the strong topology since it is a linear space. $\quad\square$

REMARK 5.4. In fact, on $X_V$ there is a generalization of the full Stone–Weierstrass theorem and it can be deduced from the Daniell–Stone theorem, even in more general cases than considered here. In particular, the algebra in Theorem 5.3 generates $\mathcal{B}(X_V)$ if it separates points. We refer to [47].

LEMMA 5.5. *Let $X$ be a completely regular space, let $V, \Theta \colon X \to [1, \infty]$ have metrizable compact level sets, $V \leq c\Theta$ for some $c \in (0, \infty)$, and such that, for all $R > 0$, there exists $R' \geq R$ such $\{V \leq R\}$ is contained in the closure of the set $\{V \leq R'\} \cap X_\Theta$.*
    *Then $C_V \subset C_\Theta$ continuously and densely.*

PROOF.   Note that $X_\Theta \subset X_V$. If $f \in C_V$, then, for $R \in (0, \infty)$,

$$|f| \leq \left(\sup_{\{V \geq \sqrt{R}\}} \frac{|f|}{V}\right)V + \sqrt{R}\|f\|_V,$$

hence,

$$\sup_{\{\Theta \geq R\}} \frac{|f|}{\Theta} \leq c \sup_{\{V \geq \sqrt{R}\}} \frac{|f|}{V} + \frac{1}{\sqrt{R}}\|f\|_V.$$

Letting $R \to \infty$, we conclude that $f\!\restriction_{X_\Theta} \in C_\Theta$. Moreover, the last assumption implies that, if $f \in C_V$ vanishes on $X_\Theta$, then it vanishes on $\{V \leq R\}$ for every $R > 0$, since $f$ is continuous on $\{V \leq R'\}$. Hence, the restriction to $X_\Theta$ is an injection $C_V \to C_\Theta$. Since $V \leq \Theta$, the injection is continuous. The density follows from Corollary 5.3. Indeed, we have seen in its proof that $\sigma(C_V) = \mathcal{B}(X_V)$. But then $\sigma(C_{V\restriction_{X_\Theta}}) = \sigma(C_V) \cap X_\Theta \supset \mathcal{B}(X_V) \cap X_\Theta = \mathcal{B}(X_\Theta)$, since $X_\Theta \in \mathcal{B}(X)$.   □

Now we come to our concrete situation.

COROLLARY 5.6.   *For $p \in [2, \infty)$, $p' \geq p$, and $x \in (0, \infty)$, $\kappa' \geq \kappa$, we have $WC_{p,\kappa} \subset WC_{p',\kappa'}$ and $WC_{p,\kappa} \subset W_1 C_{p,\kappa} \subset W_1 C_{p',\kappa'}$ densely and continuously.*

PROOF.   Note that, for $x \in L^p(0,1)$, $p > 1$, $P_N x \in H_0^1$, $N \in \mathbb{N}$, and $P_m x \to x$ in $L^p(0,1)$ as $m \to \infty$ (see, e.g., [40], Section 2c16). Also by (2.8), $V_{p,\kappa} \circ P_N \leq \alpha_p^p V_{p,\kappa}$ and, hence, $\{P_N x | V_{p,\kappa}(x) \leq R, N \in \mathbb{N}\} \subset \{V_{p,\kappa} \leq \alpha_p^p R\} \cap H_0^1$. Furthermore, since

$$\left(\frac{2}{p}\right)^2 |(|x|^{p/2})'|_2^2 = |x'|x|^{p/2-1}|_2^2$$

$$= ||x'\mathbb{1}_{\{|x|<1\}}|x|^{p/2-1}|_2^2 + ||x'\mathbb{1}_{\{|x|\geq 1\}}|x|^{p/2-1}|_2^2$$

$$\leq |x'|_2^2 + |x'|x|^{p'/2-1}|_2^2,$$

it follows that there exists $c_p \in (0, \infty)$ such that

(5.3)                        $\Theta_{p,\kappa} \leq c_p \Theta_{p',\kappa'}$.

Now the assertion follows from Lemma 5.5.   □

LEMMA 5.7.   *Let $l \in \mathbb{Z}_+$, $p \in [2, \infty)$, $\kappa \in (0, \infty)$, $(f_n)_{n \in \mathbb{N}} \subset \text{Lip}_{l,p,\kappa}$, be such that $f(x) := \lim_{n \to \infty} f_n(x)$ exists for all $x \in X_p$. Then*

$$\|f\|_{p,\kappa} \leq \liminf_{n \to \infty} \|f_n\|_{p,\kappa} \quad \text{and} \quad (f)_{l,p,\kappa} \leq \liminf_{n \to \infty} (f_n)_{l,p,\kappa}.$$

*In particular, $(\text{Lip}_{l,p,\kappa}, \|\cdot\|_{\text{Lip}_{l,p,\kappa}})$ is complete.*

PROOF. The assertion follows from the fact that, for a set $\Omega$ and $\psi_n : \Omega \to \mathbb{R}$, $n \in \mathbb{N}$, we have $\sup_{\omega \in \Omega} \liminf_{n \to \infty} \psi_n(\omega) \leq \liminf_{n \to \infty} \sup_{\omega \in \Omega} \psi_n(\omega)$. □

PROPOSITION 5.8. *Let* $l \in \mathbb{Z}_+$, $p \in [2, \infty)$, *and* $\kappa \in (0, \infty)$. *Let* $(f_n)_{n \in \mathbb{N}}$ *be a bounded sequence in* $\mathrm{Lip}_{l, p, \kappa}$. *Then there exists a subsequence* $(f_{n_k})_{k \in \mathbb{N}}$ *converging pointwise to some* $f \in \mathrm{Lip}_{l, p, \kappa}$.
*If* $l > 0$, *then* $f$ *is sequentially weakly continuous on* $X_p$.

PROOF. Let $Y \subset X_p$ be countable such that $Y \cap \{V_{p,\kappa} < n\}$ is $|\cdot|_p$-dense in $\{V_{p,\kappa} < n\}$ for all $n \in \mathbb{N}$, and let $(f_{n_k})_{k \in \mathbb{N}}$ be a subsequence converging pointwise on $Y$. Since $f_{n_k}$, $k \in \mathbb{N}$, are bounded in $\mathrm{Lip}_{l,p,\kappa}$, they are $|\cdot|_p$-equicontinuous on the $|\cdot|_p$-open sets $\{V_{p,\kappa} < n\}$ for all $n \in \mathbb{N}$. Hence, there exists a $|\cdot|_p$-continuous function $f : X_p \to \mathbb{R}$ such that $f_{n_k}(x) \to f(x)$ as $k \to \infty$ for all $x \in X_p$. By Lemma 5.7, we have $f \in \mathrm{Lip}_{l,p,\kappa}$.

Since $f_{n_k}$, $k \in \mathbb{N}$, are $|(-\Delta)^{-l/2} \cdot |_2$-equicontinuous, $f$ is $|(-\Delta)^{-l/2} \cdot |_2$-continuous, in particular, sequentially weakly continuous on $X_p$. □

## 6. Construction of resolvents and semigroups.
In this section we construct the resolvent and semigroup in the spaces $WC_{p,\kappa}$ associated with the differential operator $L$ defined in (1.2) with $F$ satisfying (F2).

PROPOSITION 6.1. *Let* $F : H_0^1 \to X$ *satisfying* (F2a) *and* (F2c), *and let* $\kappa_0$, $Q_{\mathrm{reg}}$ *and* $\lambda_{q,\kappa}$, $m_{q,\kappa}$ *for* $q \in Q_{\mathrm{reg}}$ *be as in* (F2a). *Assume that* $V_{\kappa_1} F^{(k)} \in W_1 C_{q,\kappa}$ *for all* $k \in \mathbb{N}$ *for some* $q \in Q_{\mathrm{reg}}$ *and* $\kappa \in (0, \kappa_0)$, $\kappa_1 \in [0, \kappa)$. *Then we have*

$$(6.1) \qquad \|u\|_{q,\kappa} \leq \frac{1}{m_{q,\kappa}} \|\lambda u - Lu\|_{1,q,\kappa} \qquad \forall\, u \in \mathcal{D}_{\kappa_1}, \lambda \geq \lambda_{q,\kappa}.$$

For the the proof of this proposition, we need the following two results.

LEMMA 6.2. *Let* $q \in [2, \infty)$, $\kappa \in (0, \infty)$.

(i) $V_{q,\kappa}$ *is Gâteaux differentiable on* $L^q(0,1)$ *with derivative given by*

$$(6.2) \quad DV_{q,\kappa}(x) = V_{q,\kappa}(x) \left( 2\kappa x + \frac{q}{1+|x|_q^q} x |x|^{q-2} \right) \ (\in L^{q/(q-1)}(0,1)).$$

(ii) *On* $H_0^1$ *the function* $V_{q,\kappa}$ *is twice Gâteaux differentiable. Moreover,* $DV_{q,\kappa} : H_0^1 \to H_0^1$ $[\subset H^{-1},$ *see* (2.1)$]$, $D^2 V_{q,\kappa} : H_0^1 \to \mathcal{L}(L^2(0,1))$ $[:=$ *bounded linear operators on* $L^2(0,1)]$. *Furthermore, both maps are continuous and, for* $x, \xi, \eta \in H_0^1$,

$$(\xi, D^2 V_{q,\kappa}(x)\eta) = V_{q,\kappa}(x) \left[ \left( 2\kappa(\xi, x) + q \frac{(\xi, x|x|^{q-2})}{1+|x|_q^q} \right) \right.$$

$$\times \left( 2\kappa(\eta, x) + q\frac{(\eta, x|x|^{q-2})}{1+|x|_q^q} \right)$$

$$(6.3) \qquad + 2\kappa(\xi, \eta) + q(q-1)\frac{(\xi, \eta|x|^{q-2})}{1+|x|_q^q}$$

$$- q^2\frac{(\xi, x|x|^{q-2})(\eta, x|x|^{q-2})}{(1+|x|_q^q)^2} \Big].$$

PROOF.  Identities (6.2) and (6.3) follow from the formulas

$$\frac{\partial}{\partial \eta}|x|_q^q = q(\eta, x|x|^{q-2}), \qquad \frac{\partial}{\partial \xi}(x|x|^{q-2}, \eta) = (q-1)(|x|^{q-2}, \xi\eta) \quad \text{and}$$

$$\frac{\partial^2}{\partial \xi \, \partial \eta}|x|_q^q = q(q-1)(\xi, \eta|x|^{q-2}), \qquad x, \xi, \eta \in H_0^1.$$

The continuity of $DV_{q,\kappa}$ and $D^2V_{q,\kappa}$ in the mentioned topologies follows from the fact that, given $x_n \to x$ in $H_0^1$ as $n \to \infty$, then $x_n' \to x'$ in $L^2(0,1)$ and $x_n \to x$ in $C_0[0,1]$ as $n \to \infty$.  □

LEMMA 6.3.  Let $q \in [2,\infty)$ and $\kappa \in (0,\infty)$. Let $u \in WC_{q,\kappa}$ be such that $u = u \circ P_N$ for some $N \in \mathbb{N}$. Then there exists $x_0 \in (C_0 \cap C_b^1)(0,1)$ such that $\|u\|_{q,\kappa} = \frac{|u|}{V_{q,\kappa}}(x_0)$.

PROOF.  We may assume $u \not\equiv 0$. Since $V_{q,\kappa}^{-1}|u|$ is weakly upper semi-continuous on $X$ and $V_{q,\kappa}$ has weakly compact level sets, there exists $x_0 \in X_q$ such that $\|u\|_{q,\kappa} = |u|(x_0)V_{q,\kappa}^{-1}(x_0)$. Set $x_1 := P_N x_0$ and $x_2 := x_0 - x_1$. Since $u(x_0) = u(P_N x_0)$, we conclude that

$$V_{q,\kappa}(x_0) = \min\{V_{q,\kappa}(x) | x \in X, P_N x = P_N x_0\}.$$

Hence, by Lemma 6.2(i), we have that $(DV_{q,\kappa}(x_0), \eta) = 0$ for all $\eta \in L^q(0,1) \cap E_N^\perp$. Since $\{\eta_k | k \in \mathbb{N}\}$ is a Schauder basis of $L^s(0,1)$ for all $s \in (1,\infty)$ (cf. [40], Section 2c16), it follows that $DV_{q,\kappa}(x_0) \in E_N \subset (C_0 \cap C_b^1)(0,1)$.

Consider  $h \in C^1(\mathbb{R})$,  $h(s) := 2\kappa s + \frac{q}{1+|x_0|_q^q}s|s|^{q-2}$,  $s \in \mathbb{R}$. By (6.2), $DV_{q,\kappa}(x_0) = V_{q,\kappa}(x_0)h \circ x_0$. Hence, $h \circ x_0 \in (C_0 \cap C_b^1)(0,1)$. Since, for $s \in \mathbb{R}$, $h'(s) = 2\kappa + \frac{q(q-1)}{1+|x_0|_q^q}|s|^{q-2} \geq 2\kappa > 0$, the assertion follows, by the inverse function theorem.  □

PROOF OF PROPOSITION 6.1.  For $N \in \mathbb{N}$, we introduce a differential operator $L^{(N)}$ on the space of all continuous functions $v: H_0^1 \to \mathbb{R}$ having continuous partial derivatives up to second order in all directions $\eta_k$, $k \in \mathbb{N}$,

defined by

$$L^{(N)}v(x) \equiv \frac{1}{2}\sum_{i=1}^{N}A_{ii}\partial_{ii}^2 v(x) + \sum_{k=1}^{N}((x,\eta_k'') + (F(x),\eta_k))\partial_k v(x), \qquad x \in H_0^1.$$

Let $\lambda \geq \lambda_{q,\kappa}$, $u \in \mathcal{D}_{\kappa_1}$, $u = u \circ P_N$ for some $N \in \mathbb{N}$. Then, for $m \geq N$ and $x \in H_0^1$,

$$(\lambda - L)u = (\lambda - L^{(m)})u = -V_{q,\kappa}L^{(m)}(uV_{q,\kappa}^{-1}) - 2(A_m DV_{q,\kappa}, D(uV_{q,\kappa}^{-1}))$$
$$+ uV_{q,\kappa}^{-1}(\lambda - L^{(m)})V_{q,\kappa}.$$

Since $u \in \mathcal{D}_{\kappa_1} \subset WC_{q,\kappa}$, Lemma 6.3 implies that there exists $x_0 \in (C_0 \cap C_b^1)(0,1)$ such that $\|u\|_{q,\kappa} = \frac{|u|}{V_{q,\kappa}}(x_0)$. We may assume, without loss of generality, that $u(x_0) \geq 0$. Then $x_0$ is a point, where the function $uV_{q,\kappa}^{-1}$ achieves its maximum. Hence,

$$D(uV_{q,\kappa}^{-1})(x_0) = 0 \quad \text{and} \quad L^{(m)}(uV_{q,\kappa}^{-1})(x_0) \leq 0.$$

Therefore,

$$(\lambda - L)u(x_0) \geq \|u\|_{q,\kappa}\liminf_{m\to\infty}(\lambda - L^{(m)})V_{q,\kappa}(x_0).$$

For $m \in \mathbb{N}$, let now $L_m$ be as in (4.1). Note that

$$|L_m(V_{q,\kappa}\restriction_{E_m}) \circ P_m - L^{(m)}V_{q,\kappa}|(x) \to 0 \qquad \text{as } m \to \infty, \ x \in H_0^1.$$

This is so since $A$ is of trace class, (F2c) holds and, for $x \in H_0^1$, $P_m x \to x$ in $H_0^1$ as $m \to \infty$ and hence, by Lemma 6.2(ii), $DV_{q,\kappa}(P_m x) \to DV_{q,\kappa}(x)$ in $H_0^1$ and $D^2 V_{q,\kappa}(P_m x) \to D^2 V_{q,\kappa}(x)$ in $\mathcal{L}(L^2(0,1))$ as $m \to \infty$. Hence, by (F2a),

$$(\lambda - L)u(x_0) \geq \|u\|_{q,x}\liminf_{m\to\infty}(\lambda - L_m)(V_{q,\kappa}\restriction_{E_m})(P_m x_0)$$
$$\geq m_{q,\kappa}\|u\|_{q,x}\liminf_{m\to\infty}\Theta_{q,\kappa}(P_m x_0) = m_{q,\kappa}\|u\|_{q,x}\Theta_{q,\kappa}(x_0).$$

Since, by assumption, $V_{\kappa_1}F^{(k)} \in W_1 C_{q,\kappa}$, $k \in \mathbb{N}$, it follows that $Lu \in W_1 C_{q,\kappa}$. So, the assertion follows. $\square$

Now we can prove our main existence result on resolvents and semigroups (see also Proposition 6.7 below).

THEOREM 6.4. *Let* (A), (F2) *hold, and let* $\kappa_0$, $Q_{\mathrm{reg}}$ *be as in* (F2a), $\kappa \in (0,\kappa_0)$ *and* $p \in Q_{\mathrm{reg}}$ *be as in* (F2d). *Let* $\kappa^* \in (\kappa,\kappa_0)$, $\kappa_1 \in (0,\kappa^* - \kappa]$, *and let* $\lambda_{p,\kappa^*}$ *and* $\lambda_{2,\kappa_1}'$ *be as in Corollary* 4.2, *with* $\kappa^*$ *and* $\kappa_1$, *respectively, replacing* $\kappa$.

Then for $\lambda > \lambda_{p,\kappa^*} \vee \lambda'_{2,\kappa_1}$, $((\lambda - L), \mathcal{D}_{\kappa_1})$ is one-to-one and has a dense range in $W_1 C_{p,\kappa^*}$. Its inverse $(\lambda - L)^{-1}$ has a unique bounded linear extension $G_\lambda \colon W_1 C_{p,\kappa^*} \to W C_{p,\kappa^*}$, defined by the following limit:

$$\lambda G_\lambda f := \lim_{m \to \infty} \lambda G_\lambda^{(m)} f, \qquad f \in \mathrm{Lip}_{0,2,\kappa_1}, \ f \ bounded, \lambda > \lambda_{p,\kappa^*} \vee \lambda'_{2,\kappa_1},$$

weakly in $W C_{p,\kappa^*}$ (hence, pointwise on $X_p$), uniformly in $\lambda \in [\lambda_*, \infty)$ for all $\lambda_* > \lambda_{p,\kappa^*} \vee \lambda'_{2,\kappa}$. Furthermore,

$$\lim_{m \to \infty} \lambda(\lambda - L) G_\lambda^{(m)} f = \lambda f$$

weakly in $W_1 C_{p,\kappa^*}$ uniformly in $\lambda \in [\lambda^*, \infty)$. $G_\lambda$, $\lambda > \lambda_{p,\kappa^*} \vee \lambda'_{2,\kappa_1}$, is a Markovian pseudo-resolvent on $W_1 C_{p,\kappa^*}$ and a strongly continuous quasi-contractive resolvent on $W C_{p,\kappa^*}$ with $\|G_\lambda\|_{W C_{p,\kappa^*} \to W C_{p,\kappa^*}} \leq (\lambda - \lambda_{p,\kappa^*})^{-1}$. $G_\lambda$ is associated with a Markovian quasi-contractive $C_0$-semigroup $P_t$ on $W C_{p,\kappa^*}$ satisfying

$$\|P_t\|_{W C_{p,\kappa^*} \to W C_{p,\kappa^*}} \leq e^{\lambda_{p,\kappa^*} t}, \qquad t > 0.$$

For the proof of the theorem, we need the following lemma.

LEMMA 6.5.  Let $G_\lambda$, $\lambda > \lambda_0$, be a pseudo-resolvent on a Banach space $\mathbb{F}$, such that $\|\lambda G_\lambda\|_{\mathbb{F} \to \mathbb{F}} \leq M$ for all $\lambda > \lambda_0$. Then the set $\mathbb{F}_G$ of strong continuity of $G$,

$$\mathbb{F}_G := \{f \in \mathbb{F} | \lambda G_\lambda f \to f \ as \ \lambda \to \infty\},$$

is the (weak) closure of $G_\lambda \mathbb{F}$.

PROOF.  First observe that $\mathbb{F}_G$ is a closed linear subspace of $\mathbb{F}$. Indeed, let $f \in \mathbb{F}$, $f_n \in \mathbb{F}_G$, $n \in \mathbb{N}$, such that $f_n \to f$ as $n \to \infty$. Then

$$\lambda G_\lambda f - f = (\lambda G_\lambda f_n - f_n) + (\lambda G_\lambda - \mathrm{id})(f - f_n).$$

The first term in the right-hand side vanishes as $\lambda \to \infty$ for all $n \in \mathbb{N}$ and the second term vanishes as $n \to \infty$ uniformly in $\lambda > \lambda_0$, since $\|\lambda G_\lambda - \mathrm{id}\|_{\mathbb{F} \to \mathbb{F}} \leq M + 1$. So, we conclude that $\lambda G_\lambda f \to f$ as $\lambda \to \infty$.

By the resolvent identity, for $f \in \mathbb{F}$ and $\lambda, \mu > \lambda_0$, we have

$$\lambda G_\lambda G_\mu f = \frac{\lambda}{\lambda - \mu} G_\mu f - \frac{1}{\lambda - \mu} \lambda G_\lambda f \to G_\mu f$$

as $\lambda \to \infty$ since $\|\lambda G_\lambda f\|_{\mathbb{F}} \leq M \|f\|_{\mathbb{F}}$. Thus, $\overline{G_\lambda \mathbb{F}} \subset \mathbb{F}_G$. On the other hand, $\mathbb{F}_G \subset \overline{G_\lambda \mathbb{F}}$ by definition. Finally, since $G_\lambda \mathbb{F}$ is linear, its weak and strong closures coincide, by the Mazur theorem.  $\square$

PROOF OF THEOREM 6.4. We have that (F2d) holds with $\kappa^*$ replacing $\kappa$, and for all $k \in \mathbb{N}$, that $F^{(k)} \in W_1 C_{p,\kappa}$ by (F2c) and (F2d), so $V_{\kappa_1} F^{(k)} \in W_1 C_{p,\kappa^*}$. Therefore, Proposition 6.1 implies that $(\lambda - L) : \mathcal{D}_{\kappa_1} \to W_1 C_{p,\kappa^*}$ is one-to-one with bounded left inverse from $W_1 C_{p,\kappa^*} \supset (\lambda - L)(\mathcal{D}_{\kappa_1})$ to $W C_{p,\kappa^*}$ for all $\lambda > \lambda_{p,\kappa^*}$.

Now we prove that $(\lambda - L)(\mathcal{D}_{\kappa_1})$ is dense in $W_1 C_{p,\kappa^*}$ for $\lambda > \lambda'_{2,\kappa_1}$. Let $m \in \mathbb{N}$, $f \in \mathrm{Lip}_{0,2,\kappa_1}(\subset W_1 C_{p,\kappa^*})$, $f$ bounded, and $\lambda > \lambda'_{2,\kappa_1}$. By Corollary 4.2, $G_\lambda^{(m)} f \in \bigcap_{\varepsilon > 0} \mathcal{D}_{\kappa_1 + \varepsilon}$ and, by (4.6), $(\lambda - L) G_\lambda^{(m)} f(x) \to f(x)$ as $m \to \infty$ for all $x \in H_0^1$, and by (4.5) and (F2d),

$$|(\lambda - L) G_\lambda^{(m)} f(x) - (f \circ P_m)(x)| \le \frac{2}{\lambda - \lambda'_{2,\kappa_1}} \Theta_{p,\kappa^*}(x) V_{2,\kappa_1}(x)(f)_{0,2,\kappa_1}$$

$$= \frac{2}{\lambda - \lambda'_{2,\kappa_1}} \Theta_{p,\kappa^*}(x)(f)_{0,2,\kappa_1}.$$

Hence, $|\lambda(\lambda - L) G_\lambda^{(m)} f - \lambda f| \to 0$ as $m \to \infty$ weakly in $W_1 C_{p,\kappa^*}$, uniformly in $\lambda \in [\lambda_*, \infty)$ for all $\lambda_* > \lambda'_{2,\kappa_1}$, by Theorem 5.1. By Corollary 5.3, $\mathcal{D}(\subset \mathrm{Lip}_{0,2,\kappa_1})$ is dense in $W_1 C_{p,\kappa^*}$. So, taking $f \in \mathcal{D}$ and recalling that $G_\lambda^{(m)} f \in \mathcal{D}_{\kappa_1}$ by Corollary 4.2, we conclude that $(\lambda - L)(\mathcal{D}_{\kappa_1})$ is of (weakly) dense range. Therefore, for $\lambda > \lambda'_{2,\kappa_1} \vee \lambda_{p,\kappa^*}$, the left inverse $(\lambda - L)^{-1}$ can be extended to a bounded linear operator $G_\lambda : W_1 C_{p,\kappa^*} \to W C_{p,\kappa^*}$. Then one has $\lambda G_\lambda^{(m)} f \to \lambda G_\lambda f$ as $m \to \infty$ weakly in $W C_{p,\kappa^*}$ (in particular, pointwise on $X_p$) for all $\lambda > \lambda'_{2,\kappa_1} \vee \lambda_{p,\kappa^*}$ and all $f \in \mathrm{Lip}_{0,2,\kappa_1} \cap \mathcal{B}_b(X)$. So, $\lambda G_\lambda$ is Markovian and $\lambda \mapsto G_\lambda f$ is decreasing if $f \ge 0$ for such $\lambda$, since $G_\lambda^{(m)} f$ has the same properties. In addition, for $\nu \in W C'_{p,\kappa^*}$, $\nu \ge 0$ (cf. Theorem 5.1), and $\lambda > \lambda^* > \lambda'_{2,\kappa_1} \vee \lambda_{p,\kappa^*}$,

$$\int |\lambda G_\lambda^{(m)} f - \lambda G_\lambda f| \, d\nu \le \int G_{\lambda^*} |\lambda(\lambda - L) G_\lambda^{(m)} f - \lambda f| \, d\nu.$$

Therefore, the weak convergence of $(\lambda G_\lambda^{(m)} f)_{m \in \mathbb{N}}$ to $\lambda G_\lambda f$ in $W C_{p,\kappa^*}$ is, in fact, uniformly in $\lambda \in [\lambda^*, \infty)$. Furthermore, by (4.3), $(\lambda - \lambda_{p,\kappa^*}) \|G_\lambda f\|_{p,\kappa^*} \le \|f\|_{p,\kappa^*}$, since $P_N \to \mathrm{id}_{X_p}$ strongly as $N \to \infty$ by [40], Section 2c16. Because $\mathcal{D}$ is dense in $W C_{p,\kappa^*}$, it follows that

$$\|G_\lambda\|_{W C_{p,\kappa^*} \to W C_{p,\kappa^*}} \le (\lambda - \lambda_{p,\kappa^*})^{-1}$$

by continuity. Note that, for $u \in \mathcal{D}_{\kappa_1}$, $\lambda, \mu > \lambda'_{2,\kappa_1} \vee \lambda_{p,\kappa^*}$, one has $u - G_\mu(\lambda - L)u = (\mu - \lambda) G_\mu u$ since $G_\mu$ is the left inverse to $(\mu - L)$. Hence, for $f \in (\lambda - L)(\mathcal{D}_{\kappa_1})$, we have, by substituting $u := G_\lambda f$, $G_\lambda f - G_\mu f = (\mu - \lambda) G_\mu G_\lambda f$, which is the resolvent identity. Since $(\lambda - L)(\mathcal{D}_{\kappa_1})$ is dense in $W_1 C_{p,\kappa^*}$ for $\lambda > \lambda'_{2,\kappa_1}$, we conclude that $G_\lambda$, $\lambda > \lambda'_{2,\kappa_1} \vee \lambda_{p,\kappa^*}$ is a pseudo-resolvent on $W_1 C_{p,\kappa^*}$, quasi-contractive in $W C_{p,\kappa^*}$.

Now we are left to prove that $G_\lambda$ is strongly continuous on $WC_{p,\kappa^*}$. Then the last assertion will follow by the Hille–Yoshida theorem. Let $f \in \mathcal{D}$ and let $N \in \mathbb{N}$ be such that $f = f \circ P_N$. Then, for all $x \in X_p$, $m \geq N$, $\lambda \geq \lambda_* > \lambda'_{2,\kappa_1} \vee \lambda_{p,\kappa^*}$,

$$|\lambda G_\lambda f(x) - f(x)| \leq |\lambda G_\lambda f - \lambda G_\lambda^{(m)} f|(x) + |\lambda G_\lambda^{(m)} f(P_N x) - f(P_N x)|.$$

As we have seen above, the first term in the right-hand side vanishes as $m \to \infty$ uniformly in $\lambda \in [\lambda_*, \infty)$. The second term in the right-hand side vanishes as $\lambda \to \infty$ for each $m \geq N$, by Corollary 4.2. Since $(\lambda - \lambda_{p,\kappa^*}) G_\lambda$ is quasi-contractive on $WC_{p,\kappa^*}$, it follows that $\lambda G_\lambda f \to f$ weakly in $WC_{p,\kappa^*}$ as $\lambda \to \infty$, by Theorem 5.1. Hence, by Lemma 6.5, $G_\lambda$ is strongly continuous on the closure of $\mathcal{D}$ in $WC_{p,\kappa^*}$. However, by Corollary 5.3, this closure is the whole space $WC_{p,\kappa^*}$. $\square$

REMARK 6.6. Since by (5.3) condition (F2d) holds with $p' \in [p, \infty) \cap Q_{\mathrm{reg}}$, $\kappa' \geq \kappa$, if it holds with $p \in [2, \infty)$, $\kappa \in (0, \infty)$, the above theorem (and, correspondingly, any of the results below) holds for any $\kappa^* \in (\kappa, \kappa_0)$ and with $p$ replaced by any $p' \in [p, \infty) \cap Q_{\mathrm{reg}}$. We note that the corresponding resolvents, hence, also the semigroups, are consistent when applied to functions in $\mathcal{D}$. In particular, the resolvents and semigroups of kernels constructed in the following proposition coincide for any $\kappa^* \in (\kappa, \kappa_0)$ and $p' \in [p, \infty) \cap Q_{\mathrm{reg}}$.

Next we shall prove that both $G_\lambda$ and $P_t$ in Theorem 6.4 above are given by kernels on $X_p$ uniquely determined by $L$ under a mild "growth condition."

PROPOSITION 6.7 (Existence of kernels). *Consider the situation of Theorem* 6.4, *let* $\lambda > \lambda_{p,\kappa^*} \vee \lambda'_{2,\kappa_1}$ *and* $t > 0$, *and let* $G_\lambda$ *and* $P_t$ *be as constructed there. Then:*

(i) *There exists a kernel* $g_\lambda(x, dy)$ *from* $X_p$ *to* $H_0^1$ *such that*

$$g_\lambda f(x) := \int f(y) g_\lambda(x, dy) = G_\lambda f(x) \qquad \text{for all } f \in W_1 C_{p,\kappa^*}, x \in X_p,$$

*which is extended by zero to a kernel from* $X_p$ *to* $X_p$. *Furthermore,* $\lambda g_\lambda 1 = 1$, $g_{\lambda'}$, $\lambda' > \lambda_{p,\kappa^*} \vee \lambda'_{2,\kappa_1}$, *is a resolvent of kernels and*

$$g_\lambda \Theta_{p,\kappa^*}(x) \leq \frac{1}{m_{p,\kappa^*}} V_{p,\kappa^*}(x) \qquad \text{for all } x \in X_p,$$

*with* $m_{p,\kappa^*}$ *as in* (F2a).

(ii) *There exists a kernel* $p_t(x, dy)$ *from* $X_p$ *to* $X_p$ *such that*

$$p_t f(x) := \int f(y) p_t(x, dy) = P_t f(x) \qquad \text{for all } f \in WC_{p,\kappa^*}, x \in X_p.$$

*Furthermore, $p_t 1 = 1$ [i.e., $p_t(x, dy)$ is Markovian], $p_\tau$, $\tau > 0$, is a measurable semigroup and*

$$p_t V_{p,\kappa^*}(x) \leq e^{\lambda_{p,\kappa^*} t} V_{p,\kappa^*}(x) \qquad \text{for all } x \in X_p.$$

(iii) *We have*

(6.4)
$$g_\lambda f(x) = \int_0^\infty e^{-\lambda \tau} p_\tau f(x) \, d\tau$$

$$\text{for all } f \in \mathcal{B}_b(X_p) \cup \mathcal{B}^+(X_p), x \in X_p.$$

[*We extend $g_\lambda$ for all $\lambda \in (0, \infty)$ using* (6.4) *as a definition.*]

(iv) *Let $x \in X_p$. Then*

$$\int_0^t p_\tau(x, X_p \setminus H_0^1) \, d\tau = 0.$$

(v) *For $x \in X_p$,*

$$\int_0^t p_\tau \Theta_{p,\kappa^*}(x) \, d\tau < \infty,$$

*so*

$$\int_0^t p_\tau |f|(x) \, d\tau < \infty \qquad \text{for all } f \in W_1 C_{p,\kappa^*}.$$

*In particular, if $u \in \mathcal{D}_{\kappa_1}$, then $\tau \mapsto p_\tau(|Lu|)(x)$ is in $L^1(0, t)$. Furthermore,*

(6.5) $\quad p_t u(x) - u(x) = \int_0^t p_\tau(Lu)(x) \, d\tau \qquad \text{for all } u \in \mathcal{D}_{\kappa_1}, x \in X_p.$

PROOF. (i) and (ii) are immediate consequences of Theorem 6.4, Corollary 5.2 and standard monotone class arguments. Equation (6.4) in (iii) holds by Theorem 6.4 for $f \in WC_{p,\kappa^*}$. Hence, (iii) follows by a monotone class argument. Now let us prove (iv). For all $f \in \mathcal{B}^+(X_p)$, by (iii), we have

(6.6) $\quad \int_0^t p_\tau f(x) \, d\tau \leq e^{\lambda t} \int_0^\infty e^{-\lambda \tau} p_\tau f(x) \, d\tau = e^{\lambda t} g_\lambda f(x), \qquad x \in X_p.$

Hence, (iv) follows with $f := \mathbb{1}_{X_p \setminus H_0^1}$ since $g_\lambda(x, X_p \setminus H_0^1) = 0$ for all $x \in X_p$. To prove (v), we just apply (6.6) to $f := \Theta_{p,\kappa^*}$ and the first two parts of the assertion follow by (i) and (iv). Now let $u \in \mathcal{D}_{\kappa_1}(\subset W_1 C_{p,\kappa^*})$. Recall that, by Theorem 6.4, $\lambda u - Lu \in W_1 C_{p,\kappa^*}$, hence, $Lu \in W_1 C_{p,\kappa^*}$, so

$$\int_0^t p_\tau(|Lu|)(x) \, d\tau < \infty \qquad \text{for all } x \in X_p.$$

Finally, to prove (6.5), first note that, for $u \in \mathcal{D}_{\kappa_1}(\subset WC_{p,\kappa^*})$, we have $G_\lambda u \in D(\tilde{L})$, where $\tilde{L}$ is the generator of $P_t$ on $WC_{p,\kappa^*}$, and

(6.7) $\qquad \tilde{L}(G_\lambda u) = -u + \lambda G_\lambda u = G_\lambda(Lu),$

since $G_\lambda$ is the left inverse of $(\lambda - L)\colon \mathcal{D}_{\kappa_1} \to W_1 C_{p,\kappa^*}$. Therefore,

$$
\begin{aligned}
(6.8) \quad \int_0^t p_\tau(g_\lambda(Lu)) \, d\tau &= \int_0^t P_\tau G_\lambda(Lu) \, d\tau = \int_0^t P_\tau \tilde{L}(G_\lambda u) \, d\tau \\
&= P_t G_\lambda u - G_\lambda u = p_t(g_\lambda u) - g_\lambda u.
\end{aligned}
$$

But integrating by parts with respect to $d\tau$, we obtain, for all $x \in X_p$,

$$
\begin{aligned}
\int_0^t p_\tau&(Lu)(x) \, d\tau \\
&= e^{\lambda t} \int_0^t e^{-\lambda \tau} p_\tau(Lu)(x) \, d\tau - \lambda \int_0^t e^{\lambda r} \int_0^r e^{-\lambda \tau} p_\tau(Lu)(x) \, d\tau \, dr \\
&= e^{\lambda t} \left[ g_\lambda(Lu)(x) - \int_t^\infty e^{-\lambda \tau} p_\tau(Lu)(x) \, d\tau \right] \\
&\quad - \lambda \int_0^t e^{\lambda r} \left[ g_\lambda(Lu)(x) - \int_r^\infty e^{-\lambda \tau} p_\tau(Lu)(x) \, d\tau \right] dr \\
&= e^{\lambda t} g_\lambda(Lu)(x) - p_t(g_\lambda(Lu))(x) \\
&\quad - (e^{\lambda t} - 1) g_\lambda(Lu)(x) + \lambda \int_0^t p_r(g_\lambda(Lu))(x) \, dr \\
&= p_t u(x) - \lambda p_t(g_\lambda u)(x) - u(x) + \lambda g_\lambda u(x) \\
&\quad + \lambda p_t(g_\lambda u)(x) - \lambda g_\lambda u(x) \\
&= p_t u(x) - u(x),
\end{aligned}
$$

where in the second to last step we used (6.8) and that, by the second equality in (6.7),

$$
g_\lambda(Lu) = -u + \lambda g_\lambda u. \qquad \square
$$

Before we prove our uniqueness result, we need the following:

LEMMA 6.8. *Consider the situation of Theorem 6.4 and let $\lambda > \lambda'_{2,\kappa_1} \vee \lambda_{p,\kappa^*}$. Then $(\lambda - L)(\mathcal{D})$ is dense in $W_1 C_{p,\kappa^*}$.*

PROOF. Let $u \in \mathcal{D}_{\kappa_1}$ and $N \in \mathbb{N}$ be such that $u = u \circ P_N$. Choose $\varphi \in C^\infty(\mathbb{R})$ such that $\varphi' \leq 0$, $0 \leq \varphi \leq 1$, $\varphi = 1$ on $[0, 1]$ and $\varphi = 0$ on $(2, \infty)$. For $n \in \mathbb{N}$, let $\varphi_n(x) := \varphi(\frac{|P_N x|_2^2}{n^2})$, $x \in X$, $u_n := \varphi_n u$. Then $u_n \in \mathcal{D}$ and

$$
Lu_n = \varphi_n Lu + uL\varphi_n + 2(Du, A_N D\varphi_n).
$$

Note that, for $i, j = 1, \ldots, N$, there are $c_j, c_{ij} \in (0, \infty)$ such that

$$
|\partial_j \varphi_n| \leq \frac{c_j}{n} \mathbb{1}_{\{|P_N x|_2 < 2n\}}, \qquad |\partial_{ij}^2 \varphi_n| \leq \frac{c_{ij}}{n^2} \mathbb{1}_{\{|P_N x|_2 < 2n\}}.
$$

Then $0 \leq \varphi_n \uparrow 1$ as $n \to \infty$, $|A_N D \varphi_n| \leq \frac{\max c_j}{n}$, and $|L \varphi_n(x)| \leq \frac{c}{n}(|x'|_2 + |P_N F|_2) \leq \frac{2c}{n} \Theta_{p,\kappa^*}(x)$ for all $x \in H_0^1$ and some $c \in (0, \infty)$ independent of $x$ and $n$ by (F2c) and (F2d). So $u_n \to u$ and $L u_n \to L u$ pointwise on $H_0^1$ and bounded in $W_1 C_{p,\kappa^*}$. Hence, by Theorems 5.1 and 6.4, it follows that $(\lambda - L)(\mathcal{D})$ is weakly, hence, strongly, dense in $W_1 C_{p,\kappa^*}$. $\square$

PROPOSITION 6.9. *Consider the situation of Theorem* 6.4 *and let* $(p_t)_{t>0}$ *be as in Proposition* 6.7. *Let* $(q_t)_{t>0}$ *be a semigroup of kernels from* $X_p$ *to* $X_p$ *such that*

$$(6.9) \quad \int_0^\infty e^{-\lambda \tau} q_\tau \Theta_{p,\kappa^*}(x) \, d\tau < \infty \qquad \text{for some } \lambda \in (0, \infty) \text{ and all } x \in X_p,$$

*and*

$$(6.10) \quad q_t u(x) - u(x) = \int_0^t q_\tau(Lu)(x) \, d\tau \qquad \text{for all } x \in X_p, u \in \mathcal{D}.$$

[*Note that the same arguments as in the proof of Proposition* 6.7(iv) *show that* $\int_0^t q_\tau(x, X_p \setminus H_0^1) \, d\tau = 0$, $x \in X_p$, *hence, the right-hand side of* (6.10) *is well-defined.*] *Then* $q_t(x, dy) = p_t(x, dy)$ *for all* $x \in X_p$, $t > 0$.

PROOF. Let $u \in \mathcal{D}$, $x \in X_p$, $t > 0$, and $\lambda$ as in (6.9). Integrating by parts with respect to $d\tau$ and then using (6.10), we obtain

$$\int_0^t e^{-\lambda \tau} q_\tau(Lu)(x) \, d\tau$$
$$= \int_0^t \lambda e^{-\lambda s} \int_0^s q_\tau(Lu)(x) \, d\tau \, ds + e^{-\lambda t} \int_0^t q_\tau(Lu)(x) \, d\tau$$
$$= \int_0^t \lambda e^{-\lambda s} q_s(u)(x) \, ds - \int_0^t \lambda e^{-\lambda s} u(x) \, ds + e^{-\lambda t}(q_t(u)(x) - u(x)),$$

so,

$$\int_0^t e^{-\lambda s} q_s(\lambda u - Lu)(x) \, ds = u(x) - e^{-\lambda t} q_t(u)(x).$$

Since (6.9) holds also with $\lambda' > \lambda$ instead of $\lambda$, we can let $\lambda \nearrow \infty$ to obtain that the resolvent $g_\lambda^q := \int_0^\infty e^{-\lambda s} q_s \, ds$, $\lambda > 0$, of $(q_t)_{t>0}$ is the left inverse of $(\lambda - L)\!\restriction_\mathcal{D}$. Hence, $g_\lambda$ and $g_\lambda^q$ coincide on $(\lambda - L)\mathcal{D}$ which is dense in $W_1 C_{p,\kappa^*}$. But by (6.9) and Theorem 5.1, $g_\lambda^q(x, dy) \in (W_1 C_{p,\kappa^*})'$ [and so is $g_\lambda(x, dy)$] for all $x \in X_p$. Hence, $g_\lambda^q = g_\lambda$. Since $t \mapsto q_t u(x)$ by (6.10) is continuous for all $u \in \mathcal{D}$, $x \in X_p$, the assertion follows by the uniqueness of the Laplace transform and a monotone class argument. $\square$

Another consequence of Lemma 6.8 is the following characterization of the generator domain of the $C_0$-semigroup $P_t$ on $W C_{p,\kappa^*}$. The second part of the following corollary will be crucial to prove the weak sample path continuity of the corresponding Markov process in the next section.

COROLLARY 6.10. *Consider the situation of Theorem* 6.4. *Let* $\bar{L}$ *denote the generator of* $P_t$ *as a* $C_0$*-semigroup on* $WC_{p,\kappa^*}$.

(i) *Then* $v \in WC_{p,\kappa^*}$ *belongs to* $\mathrm{Dom}(\bar{L})$ *if and only if there exist* $f \in WC_{p,\kappa^*}$ *and* $(u_n) \subset \mathcal{D}$ *such that* $u_n \to v$ *and* $Lu_n \to f$ *strongly, equivalently, weakly, in* $W_1C_{p,\kappa^*}$ *as* $n \to \infty$, *that is,* $u_n \to v$ *and* $Lu_n \to f$ *pointwise on* $H_0^1$, *and* $\sup_n(\|u_n\|_{1,p,\kappa^*} + \|Lu_n\|_{1,p,\kappa^*}) < \infty$. *In this case,* $\bar{L}v = f$ *and* $u_n \to v$ *weakly in* $WC_{p,\kappa^*}$ *as* $n \to \infty$.

(ii) *If* $v \in \mathrm{Dom}(\bar{L})$ *and* $v$, $\bar{L}v$ *are bounded, then the sequence* $(u_n) \subset \mathcal{D}$ *from* (i) *can be chosen uniformly bounded.*

(iii) *Let* $\lambda > \lambda_{p,\kappa^*} \vee \lambda'_{2,\kappa_1}$ *and* $v \in D(L)$ *such that* $v$, $\bar{L}v$ *are bounded, and let* $x \in X_p$. *Then there exists a Borel-measurable map* $\bar{D}_{A^{1/2}}^x v : X_p \to X$ *such that, for any sequence* $(u_n) \subset \mathcal{D}_{\kappa_1}$ *such that* $u_n \to v$, $Lu_n \to \bar{L}v$ *weakly in* $W_1C_{p,\kappa^*}$ *as* $n \to \infty$ *with* $\sup_n \|u_n\|_\infty < \infty$, *we have*

$$\lim_{n\to\infty} g_\lambda(|\bar{D}_{A^{1/2}}^x v - A^{1/2}Du_n|_2)(x) = 0.$$

*Furthermore, for all* $\chi \in C^2(\mathbb{R})$ *and* $t > 0$,

$$p_t(\chi \circ v)(x) - (\chi \circ v)(x)$$
$$= \int_0^t p_\tau(\chi' \circ v\bar{L}v)(x)\,d\tau + \int_0^t p_\tau(\chi'' \circ v(\bar{D}_{A^{1/2}}^x v, \bar{D}_{A^{1/2}}^x v))(x)\,d\tau.$$

*If, in particular,* $v = g_\lambda f$ *for some* $f \in \mathcal{D}$, *then, in addition, for all* $\kappa' \in (0, \kappa_1]$,

$$|\bar{D}_{A^{1/2}}^x v|(y) \le \frac{1}{\lambda - \lambda'_{2,\kappa'}}(f)_{0,2,\kappa'}V_{\kappa'}(y) \qquad \text{for } g_\lambda(x, dy)\text{-a.e. } y \in X_p.$$

PROOF. (i) Note that $v \in \mathrm{Dom}(\bar{L})$ if and only if $v = G_\lambda g$ for some $g \in WC_{p,\kappa^*}$, $\lambda > \lambda_{p,\kappa^*} \vee \lambda'_{2,\kappa_1}$. Given such $v$, by Lemma 6.8, there exist $u_n \in \mathcal{D}$, $n \in \mathbb{N}$, such that $(\lambda - L)u_n \to g$ in $W_1C_{p,\kappa^*}$ as $n \to \infty$. Then $u_n = G_\lambda(\lambda - L)u_n \to G_\lambda g = v$ in $WC_{p,\kappa^*}$ by Theorem 6.4, consequently,

$$Lu_n \to \lambda v - g =: f \in WC_{p,\kappa^*},$$

as $n \to \infty$ in $WC_{p,\kappa^*}$, hence, by Corollary 5.6 in $W_1C_{p,\kappa^*}$. On the other hand, let $v, f \in WC_{p,\kappa^*}$ be such that, for some $(u_n) \subset \mathcal{D}$, $u_n \to v$ and $Lu_n \to f$ weakly in $W_1C_{p,\kappa^*}$. Then, for $\lambda > \lambda'_{2,\kappa_1} \vee \lambda_{p,\kappa^*}$,

$$v = \lim_n u_n = \lim_n G_\lambda(\lambda - L)u_n = G_\lambda(\lambda v - f),$$

weakly in $WC_{p,\kappa^*}$, since, by Theorem 6.4, the latter equality holds as a weak limit in $WC_{p,\kappa^*}$ (hence, as a weak limit in $W_1C_{p,\kappa^*}$ by Corollary 5.6).

(ii) By assumption, $g := \lambda v - \bar{L}v$ is bounded. By Corollary 5.3, there exist $g_n \in \mathcal{D}$, $n \in \mathbb{N}$, which we can choose such that $\sup_n |g_n| \leq \|g\|_\infty$, converging to $g$ in $WC_{p,\kappa^*}$. Let $\lambda > \lambda_{p,\kappa^*} \vee \lambda'_{2,\kappa_1}$ and consider $v_{n,m} := G_\lambda^{(m)} g_n$, $m \in \mathbb{N}$. Then $v_{n,m} \in \mathcal{D}_{\kappa_1}$ by Corollary 4.2, and by Theorem 6.4,

$$(6.11) \qquad \lim_{m \to \infty} v_{n,m} = G_\lambda g_n \qquad \text{weakly in } WC_{p,\kappa^*}, \text{ hence, weakly in } W_1 C_{p,\kappa^*},$$

and

$$(6.12) \qquad \lim_{m \to \infty} (\lambda - L) v_{n,m} = g_n \qquad \text{weakly in } W_1 C_{p,\kappa^*}.$$

Therefore,

$$(6.13) \qquad \lim_{m \to \infty} L v_{n,m} = -g_n + \lambda G_\lambda g_n \to -g + \lambda G_\lambda g = \bar{L}v$$

weakly in $W_1 C_{p,\kappa^*}$, as $n \to \infty$. Since $\lambda G_\lambda^{(m)}$ is Markov, $v_{n,m}$, $n, m \in \mathbb{N}$, is uniformly bounded. Consequently, the pair $(v, \bar{L}v)$ lies in the weak closure of the convex set

$$(6.14) \qquad \{(u, Lu) | u \in \mathcal{D}_{\kappa_1}, \|u\|_\infty \leq \|g\|_\infty\}$$

in $W_1 C_{p,\kappa^*} \times W_1 C_{p,\kappa^*}$, hence, also in its strong closure. Repeating the same arguments as in Lemma 6.8, it follows that, in (6.14), $\mathcal{D}_{\kappa_1}$ can be replaced by $\mathcal{D}$ and assertion (ii) follows.

(iii) If $(u_n) \subset \mathcal{D}$ is a sequence as in the assertion, then, since $(u_n - u_m)^2 \in \mathcal{D}$,

$$(\lambda - L)(u_n - u_m)^2 + 2|A^{1/2} D(u_n - u_m)|^2 = 2(u_n - u_m)(\lambda - L)(u_n - u_m).$$

Hence, applying $g_\lambda(x, dy)$, we obtain

$$(u_n - u_m)^2(x) + 2g_\lambda(|A^{1/2} D(u_n - u_m)|^2)(x)$$
$$= 2g_\lambda((u_n - u_m)(\lambda - L)(u_n - u_m))(x).$$

Hence, the first assertion follows by Theorem 6.4 and Proposition 6.7(i) by Lebesgue's dominated convergence theorem. Furthermore,

$$\int_0^t p_\tau(x, dy) \, d\tau \leq e^{t\lambda} g_\lambda(x, dy),$$

$\chi(u_n) \in \mathcal{D}$, and by (6.5),

$$p_t(\chi \circ u_n)(x) - (\chi \circ u_n)(x)$$
$$= \int_0^t p_\tau(\chi' \circ u_n L u_n)(x) \, d\tau + \int_0^t p_\tau(\chi'' \circ u_n |A^{1/2} D u_n|_2^2)(x) \, d\tau.$$

Hence, the second assertion again follows by dominated convergence, since $u_n \to u$ weakly in $WC_{p,\kappa^*}$ as $n \to \infty$ by the last assertion of (i). To prove the final part of (iii), define

$$u_n := G_\lambda^{(n)} f, \qquad n \in \mathbb{N}.$$

Then by Theorem 6.4, $(u_n)$ has all properties above so that $(A^{1/2} Du_n)$ approximates $\bar{D}_{A^{1/2}}^x v$ in the above sense. But by (4.4), with $q := p$, $\kappa := \kappa'$, and Lemma 3.6,

$$|Du_n|(y) \le \frac{1}{\lambda - \lambda_{2,\kappa'}'} (f)_{0,2,\kappa'} V(y) \qquad \text{for all } y \in X \ (\supset X_p).$$ □

Next we want further regularity properties. We emphasize that these results will not be used in the next section. We extend both $g_\lambda(x, dy)$, $p_t(x, dy)$ by zero to kernels from $X_p$ to $X$.

PROPOSITION 6.11. *Consider the situation of Theorem 6.4 and let* $g_\lambda$, $p_t$ *be as in Proposition 6.7. Let* $q \in Q_{\text{reg}} \cap [2, p]$ *and* $\kappa \in [\kappa_1, \kappa^*]$ *with* $\lambda_{q,\kappa}$, $\lambda_{q,\kappa}'$ *and* $\lambda_{q,\kappa}''$ *as in Corollary 4.2. Let* $\lambda > \lambda_{q,\kappa} \vee \lambda_{p,\kappa^*} \vee \lambda_{2,\kappa_1}'$.

(i) *Let* $f \in WC_{q,\kappa}$. *Then* $g_\lambda f$ *uniquely extends to a continuous function on* $X_q$, *again denoted by* $g_\lambda f$ *such that*

$$\|g_\lambda f\|_{q,\kappa} \le \frac{1}{\lambda - \lambda_{q,\kappa}} \|f\|_{q,\kappa}. \tag{6.15}$$

*If* $f \in \text{Lip}_{0,2,\kappa_1} \cap \mathcal{B}_b(X)$, *then* $g_\lambda f$ *extends uniquely to a continuous function on* $X$, *again denoted by* $g_\lambda f$ *such that* $g_\lambda f \in \text{Lip}_{0,2,\kappa_1} \cap \mathcal{B}_b(X)$ *and for* $\lambda > \lambda_{q,\kappa}'$ *satisfying* (6.15) *and*

$$(g_\lambda f)_{0,q,\kappa} \le \frac{1}{\lambda - \lambda_{q,\kappa}'} (f)_{0,q,\kappa}. \tag{6.16}$$

*If, in addition,* (F2e) *holds, then, for* $\lambda > \lambda_{q,\kappa}''$ *and* $f \in \text{Lip}_{1,2,\kappa_1} \cap \mathcal{B}_b(X)$,

$$(g_\lambda f)_{1,q,\kappa} \le \frac{1}{\lambda - \lambda_{q,\kappa}''} (f)_{1,q,\kappa}. \tag{6.17}$$

(ii) *Let* $t > 0$ *and* $f \in \text{Lip}_{0,2,\kappa_1} \cap \mathcal{B}_b(X) \cap W_{p,\kappa^*}(\supset \mathcal{D})$. *Then* $p_t f$ *uniquely extends to a continuous function on* $X$, *again denoted by* $p_t f$, *which is in* $\text{Lip}_{0,2,\kappa_1} \cap \mathcal{B}_b(X)$, *such that*

$$\|p_t f\|_{q,\kappa} \le e^{t\lambda_{q,\kappa}} \|f\|_{q,\kappa}, \tag{6.18}$$

$$(p_t f)_{0,q,\kappa} \le e^{t\lambda_{q,\kappa}'} (f)_{0,q,\kappa}. \tag{6.19}$$

*If, in addition,* (F2e) *holds, then, for* $f \in \text{Lip}_{1,2,\kappa_1} \cap \mathcal{B}_b(X)$,

$$(p_t f)_{1,q,\kappa} \le e^{t\lambda_{q,\kappa}''} (f)_{1,q,\kappa}. \tag{6.20}$$

REMARK 6.12. (i) Because of Remark 6.6, the restriction $q \le p$ and $\kappa \in [\kappa_1, \kappa^*]$ in the above proposition are irrelevant since, for given $q \in Q_{\text{reg}}$ and $\kappa \in (0, \kappa_0)$, we can always choose $p$, $\kappa_1$, $\kappa^*$ suitably.

(ii) If (F2e) holds, by similar techniques, as in the following proof of Proposition 6.11 and by the last part of Proposition 5.8, one can prove that $p_t$ from Proposition 6.7 can be extended to a semigroup of kernels from $X$ to $X$ such that

$$\lim_{t \to 0} p_t u(x) = u(x) \qquad \text{for all } u \in \mathcal{D}, x \in X.$$

Then the proof of the first part of Theorem 7.1 in the next section implies the existence of a corresponding cadlag Markov process on $X$. However, we do not know whether this process solves our desired martingale problem, since it is not clear whether identity (6.5) holds for the above extended semigroup for all $x \in X$. As is well known and will become clear in the proof of Theorem 7.1 below, (6.5) is crucial for the martingale problem.

(iii) We emphasize that, in Proposition 6.11, it is not claimed that the extensions of $g_\lambda$ and $p_t$ satisfy the resolvent equation, have the semigroup property respectively on the larger spaces $X_q$ or $X$. It is also far from being clear whether $\lim_{t \to 0} p_t u(x) = u(x)$ for $u \in \mathcal{D}$ and all $x \in X$. Furthermore, it is also not clear whether $g_\lambda f \in WC_{q,\kappa}$ if $f \in WC_{q,\kappa}$.

PROOF OF PROPOSITION 6.11. (i) Let $f \in \mathrm{Lip}_{0,2,\kappa_1} \cap \mathcal{B}_b(X)$. Hence, by (4.3) and (4.4) [together with (2.8)] applied with $q = 2$, $\kappa = \kappa_1$, it follows by Proposition 5.8 that $(G_\lambda^{(m)} f)_{m \in \mathbb{N}}$ has subsequences converging to functions in $\mathrm{Lip}_{0,2,\kappa_1}$. Since we know by Theorem 6.4 that $(G_\lambda^{(m)} f)_{m \in \mathbb{N}}$ converges to the continuous function $G_\lambda f$ [$= g_\lambda f$ by Proposition 6.7(i)] on $X_p$ and since $X_p$ is dense in $X$, we conclude that all these limits must coincide. Hence, $g_\lambda f$ has a continuous extension in $\mathrm{Lip}_{0,2,\kappa_1}$, which we denote by the same symbol. Since $P_N \to \mathrm{id}_{X_q}$ strongly on $X_q$ as $N \to \infty$, by (4.3), (4.4) and Lemma 5.7, we obtain (6.15) and, provided $\lambda > \lambda_{q,\kappa} \vee \lambda'_{q,\kappa}$, (6.16) for such $f$, since $\mathrm{Lip}_{0,2,\kappa_1} \subset \mathrm{Lip}_{0,q,\kappa}$. If, in addition, (F2e) holds, (4.7) and Lemma 5.7 imply (6.17), provided $f \in \mathrm{Lip}_{1,2,\kappa_1} \cap \mathcal{B}_b(X)$ and $\lambda > \lambda''_{q,\kappa}$. Considering (6.15) for $f \in \mathcal{D}$, since $\mathcal{D}$ is dense in $WC_{q,\kappa}$, (6.15) extends to all of $WC_{q,\kappa}$ by continuity. For $f \in WC_{q,\kappa}$, the resulting function, lets call it $\overline{g_\lambda f}$ on $X_q$, is equal to $g_\lambda f$ on $X_p$, since by Theorem 5.1, for $u_n \in \mathcal{D}$, $n \in \mathbb{N}$, with $u_n \to f$ as $n \to \infty$ in $WC_{q,\kappa}$, it follows that $g_\lambda u_n(x) \to g_\lambda f(x)$ as $n \to \infty$ for all $x \in X_p$. So, $\overline{g_\lambda f}$ coincides with $g_\lambda f$ on $X_p$ and $\overline{g_\lambda f}$ is the desired extension. Since $X_p$ is dense in $X_q \subset X$ continuously, it follows that, for $f \in WC_{q,\kappa} \cap \mathrm{Lip}_{0,2,\kappa_1} \cap \mathcal{B}_b(X)$, the two constructed extensions of $g_\lambda f$ coincide on $X_q$ by continuity. So, (i) is completely proved.

(ii) First, we recall that, by Theorem 6.4 and Proposition 6.7(ii), since $f \in WC_{p,\kappa^*}$, $p_t f \in W_{p,\kappa^*}$ and

$$(6.21) \qquad p_t f = \lim_{n \to \infty} \left( \frac{n}{t} g_{n/t} \right)^n f \qquad \text{in } WC_{p,\kappa^*},$$

in particular, pointwise on $X_p$. But by (i), for (large enough) $n \in \mathbb{N}$, $(g_{n/t})^n f$ have continuous extensions which belong to $\mathrm{Lip}_{0,2,\kappa_1} \cap \mathcal{B}_b(X)$ and satisfy (6.15), (6.16) and, provided (F2e) holds, also (6.17) with $\lambda$ replaced by $\frac{n}{t}$. So, by Proposition 5.8, Lemma 5.7 and the same arguments as in the proof of (i), the assertion follows, since by Euler's formula, for $\lambda_0 > 0$,

$$\lim_{n \to \infty} \left( \frac{n/t}{n/t - \lambda_0} \right)^n = e^{t\lambda_0}. \qquad \square$$

**7. Solution of the martingale problem and of SPDE (1.1).** This section is devoted to the proof of the following theorem which is more general than Theorem 2.3.

THEOREM 7.1. *Assume that* (A), (F2) *hold and let* $\kappa_0$ *be as in* (F2a), $\kappa \in (0, \kappa_0)$ *and* $p \in Q_{\mathrm{reg}}$ *as in* (F2d). *Let* $\kappa^* \in (\kappa, \kappa_0)$, $\kappa_1 \in (0, \kappa^* - \kappa]$ *and let* $\lambda_{p,\kappa^*}$ *be as in Corollary* 4.2 *(with* $\kappa^*$ *replacing* $\kappa$ *there). Let* $(p_t)_{t>0}$ *be as in Proposition* 6.7(ii).

(i) *There exists a conservative strong Markov process* $\mathbb{M} := (\Omega, \mathcal{F}, (\mathcal{F}_t)_{t \geq 0}, (x_t)_{t \geq 0}, (\mathbb{P}_x)_{x \in X_p})$ *on* $X_p$ *with continuous sample paths in the weak topology whose transition semigroup is given by* $(p_t)_{t>0}$, *that is,* $\mathbb{E}_x f(x_t) = p_t f(x)$, $x \in X_p$, $t > 0$, *for all* $f \in \mathcal{B}_b(X_p)$, *where* $\mathbb{E}_x$ *denotes expectation with respect to* $\mathbb{P}_x$. *In particular,*

$$\mathbb{E}_x \left[ \int_0^\infty e^{-\lambda_{p,\kappa^*}s} \Theta_{p,\kappa^*}(x_s) \, ds \right] < \infty \qquad \text{for all } x \in X_p.$$

(ii) *("Existence") The assertion of Theorem* 2.3(ii) *holds for* $\mathbb{M}$.

(iii) *("Uniqueness") The assertion of Theorem* 2.3(iii) *holds with* $\kappa$, $\lambda_\kappa$ *replaced by* $\kappa^*$, $\lambda_{p,\kappa^*}$ *respectively.*

(iv) *If there exist* $p' \in [p, \infty)$, $\kappa' \in [\kappa^*, \kappa_0)$ *such that*

$$(7.1) \qquad \sup_{y \in H_0^1} \Theta_{p',\kappa'}^{-1}(y) |(F(y), \eta_m)|^2 < \infty \qquad \text{for all } m \in \mathbb{N},$$

*then* $\mathbb{M}$ *from assertion* (i) *weakly solves SPDE* (1.1) *for* $x \in X_{p'}$ *as initial condition.*

REMARK 7.2. (i) Due to Theorem 7.1(iv), it suffices to show that (F1) implies (7.1) to prove Theorem 2.3(iv). It follows from (F1) that, for all $m \in \mathbb{N}$ and $y \in H_0^1$,

$$|(F(y), \eta_m)| \leq |(y, \eta_m'')| + |(\Psi(y), \eta_m')| + |(\Phi(y), \eta_m)|$$
$$\leq \sqrt{2}\pi^2 m^2 (|y|_1 + |\Psi(y)|_1 + |\Phi(y)|_1).$$

Proceeding exactly as in the proof of Lemma 4.9, we find that, for all $p' \in [2, \infty)$, $\kappa' \in (0, \infty)$ up to a constant (which is independent of $y$) which is dominated by

$$\Theta_{p', \kappa'}^{(q_2 - 2 + 2/q_1)/(p'+2)}(y) + \Theta_{p', \kappa'}^{1/(2(p'+2))}(y).$$

Here we also used Remark 2.1(i). Note that $(2(p'+2))^{-1} \leq \frac{1}{2}$ and $(q_2 - 2 + \frac{2}{q_1})/(p'+2) \leq \frac{1}{2}$ if and only if $p' \geq 2q_2 - 6 + \frac{4}{q_1}$. Hence, in the latter case, (7.1) holds and, therefore, $\mathbb{M}$ weakly solves (1.1), by Theorem 7.1(iv).

(ii) Since by Remark 6.6 we can always increase $p$ as long as it is in $Q_{\text{reg}}$, which is equal to $[2, \infty)$ if (F1) holds, Theorem 7.1, in particular, implies that, for $\tilde{p} \geq p$, $\tilde{p} \in Q_{\text{reg}}$, $X_{\tilde{p}}$ is an invariant subset for the process $\mathbb{M}$ and the sample paths are even weakly continuous in $X_{\tilde{p}}$.

PROOF OF THEOREM 7.1. (i) and (ii): We mostly follow the lines of the proof of [7], Theorem I.9.4.

Let $\Omega_0 := X_p^{[0, \infty)}$ equipped with the product Borel $\sigma$-algebra $\mathcal{M}$, $x_t(\omega) := \omega(t)$ for $t > 0$, $\omega \in \Omega$ and, for $t \geq 0$, let $\mathcal{M}_t^0$ be the $\sigma$-algebra generated by the functions $x_s^0$, $0 \leq s \leq t$. By Kolmogorov's theorem, for each $x \in X_p$, there exists a probability measure $\mathbb{P}_x$ on $(\Omega_0, \mathcal{M}^0)$ such that $\mathbb{M}_0 := (\Omega_0, \mathcal{M}^0, (\mathcal{M}_t^0)_{t \geq 0}, (x_t^0)_{t \geq 0}, (\mathbb{P}_x)_{x \in X_p})$ is a conservative time homogeneous Markov process with $\mathbb{P}_x\{x_0^0 = x\} = 1$ and $p_t$ as (probability) transition semigroup.

Now we show that, for all $x \in X_p$, the trajectory $x_t^0$ is locally bounded $\mathbb{P}_x$-a.s., that is,

$$(7.2) \qquad \mathbb{P}_x\left\{\sup_{t \in [0, T] \cap \mathbb{Q}} |x_t^0|_p < \infty \ \forall T > 0\right\} = 1 \qquad \forall x \in X_p.$$

Let $g := V_{p, \kappa^*}$. Then by Proposition 6.7(iii),

$$(7.3) \qquad e^{-\lambda_{p, \kappa^*} t} p_t g(x) \leq g(x) \qquad \text{for all } x \in X_p, \ t > 0.$$

Hence, for all $x \in X_p$, the family $e^{-\lambda_{p, \kappa^*} t} g(x_t^0)$ is a super-martingale over $(\Omega_0, \mathcal{M}^0, \mathcal{M}_t^0, \mathbb{P}_x)$ since, given $0 \leq s < t$ and $Q \in \mathcal{M}_s^0$, by the Markov property,

$$\mathbb{E}_x\{e^{-\lambda_{p, \kappa^*} t} g(x_t^0), Q\} = e^{-\lambda_{p, \kappa^*} s} \mathbb{E}_x\{e^{-\lambda_{p, \kappa^*}(t-s)} p_{t-s} g(x_s^0), Q\}$$
$$\leq \mathbb{E}_x\{e^{-\lambda_{p, \kappa^*} s} g(x_s^0), Q\}.$$

Then, by [7], Theorem 0.1.5(b)

$$\mathbb{P}_x\left\{\exists \lim_{\mathbb{Q} \ni s \uparrow t} |x_s^0|_p \text{ and } \lim_{\mathbb{Q} \ni s \downarrow t} |x_s^0|_p \ \forall t \geq 0\right\} = 1 \qquad \forall x \in X_p.$$

In particular, (7.2) holds.

Now we show that $x_t^0$ can be modified to become weakly cadlag on $X_p$, that is, that

$$(7.4) \qquad \mathbb{P}_x\Big\{\exists \text{ w-}\lim_{\mathbb{Q}\ni s\uparrow t} x_s^0 \text{ and w-}\lim_{\mathbb{Q}\ni s\downarrow t} x_s^0 \ \forall t \geq 0\Big\} = 1 \qquad \forall x \in X_p.$$

For a positive $f \in \mathcal{D}$ and $\lambda > 0$, we have $e^{-\lambda t} p_t g_\lambda f \leq g_\lambda f$ for all $x \in X_p$ and $t \geq 0$. Hence, by the preceding argument, the family $e^{-\lambda t} g_\lambda f(x_t^0)$ is a super-martingale over $(\Omega_0, \mathcal{M}^0, \mathcal{M}_t^0, \mathbb{P}_x)$ for all $x \in X_p$ and

$$\mathbb{P}_x\Big\{\exists \lim_{\mathbb{Q}\ni s\uparrow t} g_\lambda f(x_s^0) \text{ and } \lim_{\mathbb{Q}\ni s\downarrow t} g_\lambda f(x_s^0) \ \forall t \geq 0\Big\} = 1 \qquad \forall x \in X_p.$$

By Proposition 6.7(i) and Theorem 6.4, we know that $\lambda g_\lambda f \to f$ as $\lambda \to \infty$ uniformly on balls in $X_p$. Since $(x_t^0)_{t\in\mathbb{Q}}$ is locally bounded in $X_p$ $\mathbb{P}_x$-a.s. for all $x \in X_p$, we conclude that

$$\mathbb{P}_x\Big\{\exists \lim_{\mathbb{Q}\ni s\uparrow t} f(x_s^0) \text{ and } \lim_{\mathbb{Q}\ni s\downarrow t} f(x_s^0) \ \forall t \geq 0\Big\} = 1 \qquad \forall x \in X_p.$$

Now let $f$ run through the countable set

$$(7.5) \qquad \tilde{\mathcal{D}} := \{\cos(\eta_k, \cdot) + 1, \sin(\eta_k, \cdot) + 1 | k \in \mathbb{N}\} \subset \mathcal{D},$$

which separates the points in $X_p$. Then we get

$$\mathbb{P}_x\Big\{\exists \lim_{\mathbb{Q}\ni s\uparrow t} f(x_s^0) \text{ and } \lim_{\mathbb{Q}\ni s\downarrow t} f(x_s^0) \ \forall t \geq 0, f \in \tilde{\mathcal{D}}\Big\} = 1 \qquad \forall x \in X_p.$$

Now (7.4) follows from the fact that $(x_t^0)_{t\in\mathbb{Q}}$ is locally in $t$ weakly relatively compact in $X_p$ $\mathbb{P}_x$-a.s. for all $x \in X_p$.

Let now

$$\Omega := \Big\{\exists \text{ w-}\lim_{\mathbb{Q}\ni s\uparrow t} x_s^0 \text{ and w-}\lim_{\mathbb{Q}\ni s\downarrow t} x_s^0 \ \forall t \geq 0\Big\},$$

$$\mathcal{M} := \{Q \cap \Omega' | Q \in \mathcal{M}^0\},$$

$$\mathcal{M}_t := \{Q \cap \Omega' | Q \in \mathcal{M}_t^0\}, \qquad t \geq 0,$$

$$x_t := \text{w-}\lim_{\mathbb{Q}\ni s\downarrow t} x_s^0, \qquad\qquad t \geq 0.$$

Then for all $x \in X_p$ and $f \in \tilde{\mathcal{D}}$, $t > 0$,

$$\mathbb{E}_x[|f(x_t^0) - f(x_t)|^2] = \lim_{\substack{s\downarrow t \\ s\in\mathbb{Q}}} \mathbb{E}_x[|f(x_t^0) - f(x_s^0)|^2]$$

$$= \lim_{\substack{s\downarrow t \\ s\in\mathbb{Q}}} (p_t f^2(x) - 2p_t(f p_{s-t} f)(x) + p_s f^2(x))$$

$$= 0,$$

since by (6.5), $t \mapsto p_t f(x)$ is continuous. Hence, $\mathbb{P}_x[x_t^0 = x_t] = 1$. Therefore, $\mathbb{M} := (\Omega, \mathcal{M}, (\mathcal{M}_{t+})_{t \geq 0}, (x_t)_{t \geq 0}, (\mathbb{P}_x)_{x \in X_p}$ is a weakly cadlag Markov process with $\mathbb{P}_x\{x_0 = x\} = 1$ and $p_t$ as transition semigroup.

Below, $\mathcal{F}, \mathcal{F}_t$ shall denote the usual completions of $\mathcal{M}, \mathcal{M}_{t+}$. Then it follows from [7], Theorem I.8.11. and Proposition I.8.11, that $\mathbb{M} := (\Omega, \mathcal{F}, (\mathcal{F}_t)_{t \geq 0}, (x_t)_{t \geq 0}, (\mathbb{P}_x)_{x \in X_p})$ is a strong Markov cadlag process with $\mathbb{P}_x\{x_0 = x\} = 1$ and $p_t$ as transition semigroup.

To prove that $\mathbb{M}$ even has weakly continuous sample paths, we first need to show that it solves the martingale problem. So, fix $x \in X_p$ and $u \in \mathcal{D}_{\kappa_1}$. It follows by Proposition 6.7(v), that for all $t \geq 0$,

$$(7.6) \qquad |Lu|(x.) \in L^1(\Omega \times [0, t], \mathbb{P}_x \otimes ds).$$

Furthermore, by (6.5) and the Markov property, it then follows in the standard way that, under $\mathbb{P}_x$,

$$(7.7) \qquad u(x_t) - u(x) - \int_0^t Lu(x_s) \, ds, \qquad t \geq 0,$$

is an $(\mathcal{F}_t)_{t \geq 0}$-martingale starting at 0.

Now we show weak continuity of the sample paths. Fix $x \in X_p$ and $f \in \mathcal{D}$. Let $\lambda > \lambda_{p,\kappa^*} \vee \lambda'_{2,\kappa_1}$, $u := g_\lambda f$ $(\in D(\bar{L}) \subset W_{p,\kappa^*})$ and $u \in \text{Lip}_{0,2,\kappa'}$ for all $\kappa' \in (0, \infty)$. Then $u$ and $Lu$ are bounded, and trivially,

$$[u(x_t) - u(x_s)]^4 = [u^4(x_t) - u^4(x_s)] - 4[u^3(x_t) - u^3(x_s)]u(x_s)$$
$$+ 6[u^2(x_t) - u^2(x_s)]u^2(x_s) - 4[u(x_t) - u(x_s)]u^3(x_s).$$

Since the martingale property is stable under $L^1(\mathbb{P}_x)$-limits, by (7.7) and (the proof of) Corollary 6.10(iii), the following processes are right continuous martingales under $\mathbb{P}_x$:

$$u(x_t) - u(x_0) - \int_0^t \bar{L}u(x_\tau) \, d\tau,$$

$$u^2(x_t) - u^2(x_0) - \int_0^t (2u\bar{L}u)(x_\tau) + |\bar{D}_{A^{1/2}}^x u|_2^2(x_\tau) \, d\tau,$$

$$u^3(x_t) - u^3(x_0) - \int_0^t (3u^2\bar{L}u)(x_\tau) + (3u|\bar{D}_{A^{1/2}}^x u|_2^2)(x_\tau) \, d\tau,$$

$$u^4(x_t) - u^4(x_0) - \int_0^t (4u^3\bar{L}u)(x_\tau) + (6u^2|\bar{D}_{A^{1/2}}^x u|_2^2)(x_\tau) \, d\tau,$$

$t \geq 0$. Hence, we obtain, for $t \geq s$,

$$\mathbb{E}_x[u(x_t) - u(x_s)]^4$$
$$= 4\mathbb{E}_x \int_s^t \bar{L}u(x_\tau)[u(x_\tau) - u(x_s)]^3 \, d\tau$$

$$+ 6\mathbb{E}_x \int_s^t |\bar{D}^x_{A^{1/2}} u|_2^2(x_\tau)[u(x_\tau) - u(x_s)]^2 \, d\tau$$

$$\leq 4\|\bar{L}u\|_\infty (t-s)^{1/4} \left(\mathbb{E}_x \int_s^t [u(x_\tau) - u(x_s)]^4 \, d\tau\right)^{3/4}$$

$$+ 6e^{\lambda t/3} (g_\lambda(|\bar{D}^x_{A^{1/2}} u|^6)(x))^{1/3} \left(\mathbb{E}_x \int_s^t |u(x_\tau) - u(x_s)|^3 \, d\tau\right)^{2/3}.$$

But by Corollary 6.10(iii) with $\kappa' = \kappa_1/6$, we have, for all $y \in X_p$,

$$g_\lambda(|\bar{D}^x_{A^{1/2}} u|^6)(y) \leq \left(\frac{1}{\lambda - \lambda'_{2,\kappa_1/4}}\right)^6 (f)_{0,2,\kappa_1/4}^6 \, g_\lambda(V_{\kappa_1})(y),$$

and by the last part of Proposition 6.7(ii),

$$g_\lambda(V_{\kappa_1})(x) \leq g_\lambda(V_{p,\kappa^*})(x)$$
$$\leq (\lambda - \lambda_{p,\kappa})^{-1} V_{p,\kappa^*}(x).$$

Therefore, for $T \in [1, \infty)$, we can find a constant $C > 0$ independent of $s, t \in [0, T]$, $t \geq s$, such that

$$(7.8) \qquad \begin{aligned} &\mathbb{E}_x[u(x_t) - u(x_s)]^4 \\ &\leq C\left[(t-s)^{1/4}\left(\mathbb{E}_x \int_s^t [u(x_\tau) - u(x_s)]^4 \, d\tau\right)^{1/4} + (t-s)^{1/6}\right] y(t), \end{aligned}$$

where, for $s \geq 0$ fixed, we set

$$(7.9) \qquad y(t) := \left(\int_s^t \mathbb{E}_x[u(x_\tau) - u(x_s)]^4 \, d\tau\right)^{1/2}, \qquad t \in [s, T].$$

Hence, we obtain from (7.8) that, for $B_T := CT^{1/4}$,

$$y'(t) \leq \tfrac{1}{2} B_T y^{1/2}(t) + \tfrac{1}{2} C(t-s)^{1/6}, \qquad t \in [s, T]$$
$$y(s) = 0.$$

Hence, for $\varepsilon > 0$, $t \in (s, T]$,

$$y'(t) \leq \frac{\varepsilon}{4} y(t) + \frac{1}{4\varepsilon} B_T^2 + \frac{C}{2}(t-s)^{1/6},$$

so, multiplying by $\exp(-\frac{\varepsilon}{4}(t-s))$ and integrating, we obtain

$$y(t) \leq \left(\frac{1}{\varepsilon^2} B_T^2 + \frac{3C}{7}(t-s)^{7/6}\right) e^{\varepsilon(t-s)/4}.$$

Choosing $\varepsilon := 4(t-s)^{-1}$, we arrive at

$$y(t) \leq (B_T^2 T^{5/6} + 2C)(t-s)^{7/6}, \qquad t \in [s, T].$$

Substituting according to (7.9) into (7.8), by the Kolmogorov–Chentsov criterion, we conclude that $t \mapsto u(x_t)$ is continuous (since by construction $x_t = \lim_{\mathbb{Q} \ni s \uparrow t} x_s^0$). Now we take $u \in \tilde{\mathcal{D}}_1 := \bigcup_{n \in \mathbb{N}} n g_n(\tilde{\mathcal{D}})$ [cf. (7.5)]. Since $\tilde{\mathcal{D}}$ separates the points of $X_p$, so does $\tilde{\mathcal{D}}_1$. It follows that the weakly cadlag path $t \to x_t$ is, in fact, weakly continuous in $X_p$.

(iii) Uniqueness is now an immediate consequence of Proposition 6.9.

(iv) As in [2], Theorem 1, one derives that componentwise $(x_t)_{t \geq 0}$ under $\mathbb{P}_x$ weakly solves the stochastic equation (1.1) for all starting points $x \in X_{p'}$. This follows from Levy's characterization theorem (since $\langle \eta_k, \cdot \rangle \in \mathcal{D}_{\kappa_1}$ $\forall k \in \mathbb{N}$) and by the fact that the quadratic variation of the weakly continuous martingale in (7.7) is equal to

$$(7.10) \qquad \int_0^t (ADu, Du)(x_s) \, ds, \qquad t \geq 0.$$

The latter can be shown by a little lengthy calculation, but it is well known in finite-dimensional situations, at least if the coefficients are bounded. For the convenience of the reader, we include a proof in our infinite-dimensional case in the Appendix (cf. Lemma A.1). Hence, assertion (iv) is completely proved. $\square$

REMARK 7.3. In Theorem 7.1(iv) SPDE (1.1) is solved in the sense of Theorem 5.7 in [2], which means componentwise. To solve it in $X_{p'}$, one needs, of course, to make assumptions on the decay of the eigenvalues of $A$ to have that $(\sqrt{A} w_t)_{t \geq 0}$ takes values in $X_{p'}$. If this is the case, by the same method as in [2], one obtains a solution to the integrated version of (1.1) where the equality holds in $X_{p'}$ (cf. [2], Theorem 6.6).

## APPENDIX

LEMMA A.1. *Consider the situation of Theorem 7.1(iv) and let $u \in \mathcal{D}$. Assume, without loss of generality, that $p' = p$, $\kappa = \kappa^*$. Let $x \in X_p$, and define, for $t \geq 0$,*

$$M_t := \left( u(x_t) - u(x_0) - \int_0^t Lu(x_r) \, dr \right)^2 - \int_0^t \Gamma(u)(x_r) \, dr,$$

*where $\Gamma(u) := (ADu, Du)$. Then $(M_t)_{t \geq 0}$ is an $(\mathcal{F}_t)_{t \geq 0}$-martingale under $\mathbb{P}_x$.*

PROOF. Let $s \in [0, t)$. We note that, by (7.1), $(M_t)_{t \geq 0}$ is a $\mathbb{P}_x$-square integrable martingale, so all integrals below are well defined. We have

$$M_t - M_s$$
$$= \left( u(x_t) - u(x_0) - \int_0^t Lu(x_r) \, dr + u(x_s) - u(x_0) - \int_0^s Lu(x_r) \, dr \right)$$

$$\times \left( u(x_t) - u(x_s) - \int_s^t Lu(x_r)\, dr \right) - 2\int_s^t \Gamma(u)(x_r)\, dr$$

$$= \left( u(x_t) + u(x_s) - 2u(x_0) - 2\int_0^s Lu(x_r)\, dr - \int_s^t Lu(x_r)\, dr \right)$$

$$\times \left( u(x_t) - u(x_s) - \int_s^t Lu(x_r)\, dr \right) - \int_s^t \Gamma(u)(x_r)\, dr$$

$$= u^2(x_t) - u^2(x_s) - 2u(x_0)(u(x_t) - u(x_s))$$

$$- 2(u(x_t) - u(x_s))\int_0^s Lu(x_r)\, dr - (u(x_t) - u(x_s))\int_0^{t-s} Lu(x_{r+s})\, dr$$

$$- (u(x_t) + u(x_s))\int_0^{t-s} Lu(x_{r+s})\, dr + 2u(x_0)\int_0^{t-s} Lu(x_{r+s})\, dr$$

$$+ 2\int_0^s Lu(x_r)\, dr \int_0^{t-s} Lu(x_{r+s})\, dr + \left( \int_0^{t-s} Lu(x_{r+s})\, dr \right)^2$$

$$- \int_s^t \Gamma(u)(x_r)\, dr.$$

Now we apply $E_x[\cdot|\mathcal{F}_s]$ to this equality and get by the Markov property that $P_x$-a.s.

$$E_x[M_t - M_s|\mathcal{F}_s]$$

$$= p_{t-s}u^2(x_s) - u^2(x_s) - 2u(x)(p_{t-s}u(x_s) - u(x_s))$$

$$- 2(p_{t-s}u(x_s) - u(x_s))\int_0^s Lu(x_r)\, dr$$

$$- 2\int_0^{t-s} p_r(Lu\, p_{t-s-r}u)(x_s)\, dr + 2u(x)\int_0^{t-s} p_r(Lu)(x_s)\, dr$$

$$+ 2\int_0^s Lu(x_r)\, dr \int_0^{t-s} p_r(Lu)(x_s)\, dr$$

$$+ 2\int_0^{t-s}\int_0^{r'} E_{x_s}[Lu(x_r)Lu(x_{r'})]\, dr\, dr' - \int_0^{t-s} p_r(\Gamma(u))(x_s)\, dr.$$

Since on the right-hand side the second and fifth, and also the third and sixth term add up to zero by Theorem 7.1(ii), we obtain

$$E_x[M_t - M_s|\mathcal{F}_s] = p_{t-s}u^2(x_s) - u^2(x_s) - 2\int_0^{t-s} p_r(Lup_{t-s-r}u)(x_s)\, dr$$

$$+ 2\int_0^{t-s}\int_0^{r'} p_r(Lup_{r'-r}(Lu))(x_s)\, dr'\, dr$$

$$- \int_0^{t-s} p_r(L(u^2))(x_s) + 2\int_0^{t-s} p_r(u\, Lu)(x_s)\, dr.$$

Since on the right-hand side the first and fourth term add up to zero and the third is, by Fubini's theorem, equal to

$$2 \int_0^{t-s} p_r \left( Lu \int_r^{t-s} p_{r'-r}(Lu) \, dr' \right)(x_s) \, dr$$

$$= 2 \int_0^{t-s} p_r (Lu(p_{t-s-r}u - u))(x_s) \, dr,$$

we see that

$$E_x[M_t - M_s | \mathcal{F}_s] = 0, \qquad P_x\text{-a.s.} \qquad \square$$

Now we shall prove Theorem 2.4, even under the weaker condition (F2), but assuming, in addition [to (F2c)], that

(A.1) $\qquad \lim_{N \to \infty} (k, F_N) = F^{(k)} \qquad$ uniformly on $H_0^1$-balls for all $k \in \mathbb{N}$,

which by Proposition 4.1 also holds under assumption (F1). So, we consider the situation of Theorem 7.1(i) and adopt the notation from there. First we need a lemma which is a modification of [9], Theorem 4.1.

LEMMA A.2.    *Let $E$ be a finite-dimensional linear space, $A : E \to \mathcal{L}(E)$ be a Borel measurable function taking values in the set of symmetric nonnegative definite linear operators on $E$ and $B : E \to E$ be a Borel measurable vector field. Let*

$$L_{A,B} u := \operatorname{Tr} A D^2 u + (B, Du), \qquad u \in C^2(E).$$

*Let $\mu$ be a probability measure on $E$ such that $L_{A,B}^* \mu = 0$ in the sense that $|A|_{E \to E}, |B|_E \in L_{\text{loc}}^1(E, \mu)$ and, for all $u \in C_c^2(E)$,*

$$\int L_{A,B} u \, d\mu = 0.$$

*Let $V : E \to \mathbb{R}_+$ be a $C^2$-smooth function with compact level sets and $\Theta : E \to \mathbb{R}_+$ be a Borel measurable function. Assume that there exists $Q \in L^1(E, \mu)$ such that $L_{A,B} V \le Q - \Theta$. Then $\Theta \in L^1(E, \mu)$ and*

$$\int \Theta \, d\mu \le \int Q \, d\mu.$$

PROOF.    Let $\xi : \mathbb{R}_+ \to \mathbb{R}_+$ be a $C^2$-smooth concave function such that $\xi(r) = r$ for $r \in [0, 1]$, $\xi(r) = 2$ for $r \ge 3$ and $0 \le \xi' \le 1$. For $k \in \mathbb{N}$, let $\xi_k(r) := k\xi(\frac{r}{k})$. Then $\xi_k$ is a $C^2$-smooth function, $\xi_k(r) = 2k$ for $r \ge 3k$, $\xi_k'' \le 0$, $0 \le \xi_k'(r) \uparrow 1$ and $\xi_k(r) \to r$ for all $r > 0$ as $k \to \infty$. Let $u_k := \xi_k \circ V - 2k$ for $x \in E_N$. Then $u_k \in C_c^2(E)$ and

$$L_{A,B} u_k(x) = \xi_k' \circ V L_{A,B} V + \xi_k'' \circ V(DV, ADV) \le \xi_k' \circ V L_{A,B} V$$

since $A$ is nonnegative definite and $\xi_k'' \leq 0$.

Now, since $\int L_{A,B} u_k \, d\mu = 0$, $0 \leq \xi_k' \circ V \leq 1$, $\Theta \geq 0$, and $L_{A,B} V \leq Q - \Theta$, we obtain

$$\int \xi_k' \circ V \Theta \, d\mu \leq \int \xi_k' \circ V Q \, d\mu.$$

Then the assertion follows from Fatou's lemma. $\quad\square$

PROOF OF THEOREM 2.4 [only assuming (F2) instead of (F1)]. (i) It follows from (F2a) that, for

$$C := \lambda \sup\{V_{p,\kappa^*}(x) | x \in H_0^1, |x'|_2 \leq 2\lambda_{p,\kappa^*}/m_{p,\kappa^*}\},$$

which is finite since $H_0^1$-balls are compact on $X_p$,

(A.2) $\qquad L_N V_{p,\kappa^*} \leq C - \dfrac{m_{p,\kappa^*}}{2} \Theta_{p,\kappa^*} \qquad$ on $E_N$ for all $N \in \mathbb{N}$.

Let $N \in \mathbb{N}$. Obviously, $V_{p,\kappa^*}(x) \to \infty$ as $|x|_2 \to \infty$, $x \in E_N$. Since $\Theta_{p,\kappa^*}(x) \to \infty$ as $|x|_2 \to \infty$, $x \in E_N$, we conclude from (A.2) that $L_N V_{p,\kappa^*}(x) \to -\infty$ as $|x|_2 \to \infty$, $x \in E_N$. Hence, a generalization of Hasminskii's theorem [8], Corollary 1.3, implies that there exists a probability measure $\mu_N$ on $E_N$ such that $L_N^* \mu_N = 0$, that is, $\int L_N u \, d\mu_N = 0$ for all $u \in C_c^2(E_N)$. Below we shall consider $\mu_N$ as a probability measure on $X_p$ by setting $\mu_N(X_p \setminus E_N) = 0$. Then, by Lemma A.2, we conclude from (A.2) that

(A.3) $$\int_X \Theta_{p,\kappa^*} \, d\mu_N \leq C.$$

Since $\Theta_{p,\kappa^*}$ has compact level sets in $X_p$, the sequence $(\mu_N)$ is uniformly tight on $X_p$. So, it has a limit point $\mu$ (in the weak topology of measures) which is a probability measure on $X_p$. Passing to a subsequence if necessary, we may assume that $\mu_N \to \mu$ weakly. Then (2.20) follows from (A.3) since $\Theta_{p,\kappa^*}$ is lower semi-continuous. In particular, $\mu(X_p \setminus H_0^1) = 0$.

Now we prove (2.19). Let $k \in \mathbb{N}$. Then it follows by (F2c), (F2d) that $F_N^{(k)} := (\eta_k, F_N) \in W_1 C_{p,\kappa^*}$. In particular, $F_N^{(k)} \in L^1(\mu_N) \cap L^1(\mu)$ for all $N \in \mathbb{N}$, due to (A.3) and (2.20). Also, the maps $x \mapsto (x'', \eta_k)$ belong to $L^1(\mu_N) \cap L^1(\mu)$ for all $N \in \mathbb{N}$ since $|(x'', \eta_k)| \leq |\eta_k''|_\infty |x|_2$. Thus, it follows from the dominated convergence theorem that $\int L_N u \, d\mu_N = 0$ for all $u \in C_b^2(E_N)$. Let $u \in \mathcal{D}$. Then we have $\mathrm{Tr}\{A_N D^2 u(x)\} = \mathrm{Tr}\{A D^2 u(x)\}$ for large enough $N$. Since $\mu_N \to \mu$ weakly, it follows that

$$\int \mathrm{Tr}\{A_N D^2 u\} \, d\mu_N \to \int \mathrm{Tr}\{A D^2 u\} \, d\mu.$$

So, we are left to show that

$$\int (F_N^{(k)} + (x'', \eta_k)) \, \partial_k u \, d\mu_N \to \int (F^{(k)} + (x'', \eta_k)) \, \partial_k u \, d\mu \qquad \text{as } N \to \infty.$$

Since $F^{(k)} \in W_1 C_{p,\kappa^*}$, by Corollary 5.3, there exists a sequence $G_{k,n} \in \mathcal{D}$ such that $\|F^{(k)} - G_{k,n}\|_{1,p,\kappa^*} \to 0$ as $n \to \infty$. Then

$$\int_X F^{(k)} \, \partial_k u(d\mu_N - d\mu)$$
$$= \int_X G_{k,n} \, \partial_k u(d\mu_N - d\mu) + \int_X (F^{(k)} - G_{k,n}) \, \partial_k u(d\mu_N - d\mu).$$

Since $\mu_N \to \mu$ weakly, we conclude that the first term vanishes as $N \to \infty$ for all $n \in \mathbb{N}$. On the other hand, the second term vanishes as $n \to \infty$ uniformly in $N \in \mathbb{N}$ since, by (A.3),

$$\int_X |F^{(k)} - G_{k,n}|(d\mu_N + d\mu)$$
$$\leq \|F^{(k)} - G_{k,n}\|_{1,p,\kappa^*} \int \Theta_{p,\kappa^*}(d\mu_N + d\mu)$$
$$\leq \left( C + \int \Theta_{p,\kappa^*} \, d\mu \right) \|F^{(k)} - G_{k,n}\|_{1,p,\kappa^*}.$$

Since $(x'', \eta_k) = (x, \eta_k'')$, the same arguments as can be applied to above $(x'', \eta_k)$. Furthermore, by (F2c), (F2d) for $R > 0$,

$$\int |F_N^{(k)} - F^{(k)}| |\partial_k u| \, d\mu_N$$
$$\leq \int_{\{\Theta_{p,\kappa^*} \leq R\}} |F_N^{(k)} - F^{(k)}| |\partial_k u| \, d\mu_N$$
$$+ 2 \sup_{\{\Theta_{p,\kappa^*} \geq R\}} \omega(\Theta_{p,\kappa^*}) \|\partial_k u\|_\infty \int \Theta_{p,\kappa^*} \, d\mu_N.$$

By (A.1) and (A.3) first letting $N \to \infty$ and then $R \to \infty$, the left-hand side of the above inequality goes to zero. So, (2.19) follows and (i) is completely proved.

(ii) Let $\mu$ be as in (i), $u \in \mathcal{D}$, and $\lambda > \lambda_{p,\kappa^*} \vee \lambda'_{2,\kappa_1}$. Then by Proposition 6.7(i) and Theorem 6.4,

$$(\text{A.4}) \qquad \int \lambda g_\lambda((\lambda - L)u) \, d\mu = \lambda \int u \, d\mu = \int (\lambda - L)u \, d\mu,$$

where we used (2.19) in the last step. By Lemma 6.8, $(\lambda - L)(\mathcal{D})$ is dense in $W_1 C_{p,\kappa^*}$ and by Proposition 6.7(i), for $f \in W_1 C_{p,\kappa^*}$,

$$g_\lambda |f| \leq \|f\|_{1,p,\kappa^*} g_\lambda \Theta_{p,\kappa^*} \leq \|f\|_{1,p,\kappa^*} \frac{1}{m_{p,\kappa^*}} V_{p,\kappa^*}$$

and $\int |f| \, d\mu \leq \int \Theta_{p,\kappa^*} \, d\mu \, \|f\|_{1,p,\kappa^*}$. So, by (2.20) and Lebesgue's dominated convergence theorem, we conclude that (A.4) extends to any $f \in W_1 C_{p,\kappa^*}$

replacing $(\lambda - L)u$. Hence, by (6.4) and Fubini's theorem, for every $f \in \mathcal{D} \subset W_1 C_{p,\kappa^*}$,

$$\lambda \int_0^\infty e^{-\lambda t} \int p_t f \, d\mu \, dt = \int u \, d\mu = \lambda \int_0^\infty e^{-\lambda t} \int u \, d\mu \, dt.$$

Since $t \mapsto p_t f(x)$ is right continuous by (6.5) for all $x \in X_p$ and bounded, assertion (ii) follows by the uniqueness of the Laplace transform and a monotone class argument. $\quad \square$

REMARK A.3. One can check that if $u \in \mathcal{D}$, $u = 0$ $\mu$-a.e., then $Lu = 0$ $\mu$-a.e. (cf. [21], Lemma 3.1, where this is proved in a similar case). Hence, $(L, \mathcal{D})$ can be considered as a linear operator on $L^s(X, \mu)$, $s \in [1, \infty)$, where we extend $\mu$ by zero to all of $X$. By [25], Appendix B, Lemma 1.8, $(L, \mathcal{D})$ is dissipative on $L^s(X, \mu)$. Then by Lemma 6.8, we know that, for large enough $\lambda$, $(\lambda - L)(\mathcal{D})$ is dense in $W_1 C_{p,\kappa^*}$ which, in turn, is dense in $L^1(X, \mu)$. Hence, the closure of $(L, \mathcal{D})$ is maximal dissipative on $L^s(X, \mu)$, that is, strong uniqueness holds for $(L, \mathcal{D})$ on $L^s(X, \mu)$ for $s = 1$. In case (F1+) holds or $\Psi = 0$, similar arguments show that our results in Section 4 imply that this is true for all $s \in [1, \infty)$ as well. A more refined analysis, however, gives that this is, in fact, true merely under condition (F2). Details will be contained in a forthcoming paper. This generalizes the main result in [16] which was proved there for $s = 2$ in the special situation when $F$ satisfies (2.15) with $\Psi(x) = \frac{1}{2}x^2$, $x \in \mathbb{R}$, and $\Phi \equiv 0$, that is, in the case of the classical stochastic Burgers equation. For more details on the $L^1$-theory for the Kolmogorov operators of stochastic generalized Burgers equations, we refer to [49].

**Acknowledgments.** The authors would like to thank Giuseppe Da Prato for valuable discussions.

The results of this paper have been announced in [48] and presented in the "Seminar on Stochastic Analysis" in Bielefeld, as well as in invited talks at conferences in Vilnius, Beijing, Kyoto, Nagoya in June, August, September 2002, January 2003 respectively, and also both at the International Congress of Mathematical Physics in Lisbon in July 2003 and the conference in Levico Terme on SPDE in January 2004. The authors would like to thank the respective organizers for these very stimulating scientific events and their warm hospitality.

DEPARTMENTS OF MATHEMATICS AND STATISTICS
PURDUE UNIVERSITY
WEST LAFAYETTE, INDIANA 47906
USA
E-MAIL: roeckner@math.purdue.edu

DEPARTMENT OF MATHEMATICS
UNIVERSITY OF WALES SWANSEA
SINGLETON PARK, SWANSEA
SA2 8PP, WALES
UNITED KINGDOM
E-MAIL: z.sobol@swansea.ac.uk